\pdfoutput=1
\documentclass[11pt]{article}
\usepackage[left=1in,right=1in,top=1in,bottom=1in]{geometry}
\usepackage{times}
\usepackage{expl3}
\usepackage{cite}
\usepackage[table]{xcolor}
\usepackage{multirow}
\usepackage{stackengine} 
\usepackage{hhline}
\usepackage{lipsum}
\usepackage{titlesec}
\usepackage{wrapfig}
\usepackage{enumerate}
\usepackage{epsfig}
\usepackage{amsmath}
\usepackage{tabularx}
\usepackage{array}
\usepackage{booktabs}
\usepackage{enumitem}
\usepackage{bbm}
\usepackage{calc}
\usepackage{graphicx}
\usepackage{amsmath}
\usepackage[title]{appendix}
\usepackage{amssymb}
\usepackage{epstopdf}
\usepackage{boldline}
\usepackage{arydshln}
\usepackage{calligra}
\usepackage{bm}
\usepackage{url}
\usepackage{blindtext}
\usepackage{accents}

\newcommand{\define}{\stackrel{\mbox{\tiny def}}{=}}

\newtheorem{theorem}{Theorem}
\newtheorem{proposition}{Proposition}

\newtheorem{lemma}{Lemma}

\usepackage{mathtools}
\usepackage{epstopdf}
\usepackage{balance}
\usepackage{thmtools}
\usepackage{thm-restate}
\usepackage{hyperref}
\usepackage{cleveref}
\usepackage[mathscr]{euscript}

\usepackage[ruled,vlined]{algorithm2e}
\include{pythonlisting}

\newcommand{\ostar}{\mathbin{\mathpalette\make@circled\star}}

\makeatletter
\newcommand{\removelatexerror}{\let\@latex@error\@gobble}
\makeatother
\makeatletter
\newcommand*{\rom}[1]{\expandafter\@slowromancap\romannumeral #1@}
\makeatother

\ExplSyntaxOn
\newcommand\latinabbrev[1]{
  \peek_meaning:NTF . {
    #1\@}%
  { \peek_catcode:NTF a {
      #1.\@ }%
    {#1.\@}}}
\ExplSyntaxOff

\titleclass{\subsubsubsection}{straight}[\subsubsection]

\begin{document}
\vspace{1cm}
\title{Generalized Double Operator Integrals for Continuous Spectrum Operators}
\vspace{1.8cm}
\author{Shih-Yu~Chang
\thanks{Shih-Yu Chang is with the Department of Applied Data Science,
San Jose State University, San Jose, CA, U. S. A. (e-mail: {\tt
shihyu.chang@sjsu.edu}). 
           }}

\maketitle

\begin{abstract}
Continuous spectrum operators (CSOs), characterized by spectra comprising continuous intervals rather than discrete eigenvalues, are pivotal in quantum mechanics, wave propagation, and systems governed by partial differential equations. Traditional double operator integrals (DOIs)—central to analyzing operator functions and perturbations—have been limited to operators with finite or countable spectra, relying critically on self-adjointness. This work introduces a comprehensive framework for Generalized Double Operator Integrals (GDOIs), extending DOI theory to non-self-adjoint operators through the spectral structure of CSOs. By reinterpreting DOIs as instances of the spectral mapping theorem for CSOs, we establish GDOIs as a rigorous generalization, enabling their application to operators with continuous spectra. Key contributions include the development of GDOIs’ algebraic properties, perturbation formulas generalizing classical results, norm and Lipschitz-type inequalities, and continuity with respect to operator and function parameters. Applications to differentiating operator-valued functions demonstrate the framework’s utility in functional calculus. Furthermore, integrating recent spectral mapping theorems allows natural extension to hybrid spectrum operators, bridging operator theory with applied fields. This work significantly expands the analytical toolbox for systems with continuous spectral phenomena, offering new methodologies for quantum mechanics, control theory, and stochastic analysis, where non-self-adjoint and continuous spectral features are fundamental. The results unify and extend existing operator-theoretic techniques, fostering interdisciplinary advances in mathematics, physics, and engineering.
\end{abstract}


\section{Introduction}\label{sec: Introduction}

Continuous spectrum operators (CSOs) are operators—often unbounded and self-adjoint—whose spectrum contains no discrete eigenvalues but instead consists of a continuous set of values, typically corresponding to an interval or the entire real line. Unlike compact operators, which have countable point spectra (eigenvalues), CSOs do not have proper eigenfunctions in the Hilbert space, though they may admit generalized eigenfunctions. ﻿Mathematically, CSOs get up naturally in purposeful analysis and operator theory, specifically within the spectral evaluation of differential and integral operators~\cite{zhang2010continuous}. In physics, they're essential in quantum mechanics, where key observables consisting of location and momentum are modeled as CSOs with spectra similar to measurable physical quantities. The continuous spectrum in these cases reflects the reality that those observables can take a continuum of values~\cite{deift1983almost}. In engineering, CSOs appear in the study of wave propagation, control theory, and systems governed by partial differential equations (PDEs), where spectral methods help analyze system stability, resonance, and energy distribution. Their rich spectral structure necessitates advanced tools such as the spectral theorem, generalized Fourier transforms, and double operator integrals to enable rigorous analysis and applications~\cite{beauchard2010controllability}.

Double Operator Integrals (DOIs) are a central tool in modern operator theory, providing a framework to analyze functions of operators and their perturbations. When the involved operators have a \textbf{finite or countable spectrum}—such as compact self-adjoint operators or matrices—the DOI framework becomes particularly tractable and insightful. DOIs have far-reaching applications in mathematics, including spectral shift functions, trace formulae, and perturbation bounds. In the sciences and engineering, they are relevant in quantum mechanics (e.g., for energy-level shifts), control theory, and signal processing where operators often act in finite-dimensional or compact settings~\cite{skripka2019multilinear,birman2003double}. Recently, The auhors in~\cite{chang2022randomDTI,chang2022randomMOI} initiated new research directions involving the \textit{probabilistic analysis} of DOIs based on finite dimensional objects, tensors. In particular, Chang's work explores \textbf{tail bounds for random DOIs}, providing one of the first approaches to understanding how DOIs behave under randomness when the underlying spectra are countable. This opens pathways to connect operator theory with probabilistic concentration, random matrix theory, and stochastic analysis, particularly in settings involving structured noise or uncertainty in finite-rank operator systems.

In the advancement of operator theory, the extension of the double operator integral (DOI) framework beyond self-adjoint operators represents a significant and influential development. Traditionally, the DOI formulation requires the assumption that the involved operators are self-adjoint, a condition that permits the use of spectral measures defined on the real line~\cite{skripka2019multilinear,birman2003double,de2002double}. However, this restriction limits the applicability of DOI to a relatively narrow class of operators. To overcome this restriction, the introduction of continuous spectrum operators (CSOs) without self-adjoint assumotion has proven to be a pivotal innovation. By leveraging the analytical structure of CSOs, it becomes feasible to define and compute double operator integrals for a broader class of non-self-adjoint operators. This advancement facilitates the transition from the classical DOI framework to what is now referred to as the generalized double operator integral (GDOI). The principal contribution of this generalization lies in its capacity to extend the robust techniques and analytical tools of DOI to a wider class of operators, effectively replacing the rigid requirement of self-adjointness with the more flexible spectral framework offered by CSOs.

This paper is structured to systematically introduce and develop the theory of Generalized Double Operator Integrals (GDOI) for continuous spectrum operators (CSOs), offering a significant extension beyond the traditional self-adjoint setting. The central contributions include the formulation of GDOI for CSOs, the establishment of their algebraic and analytic properties, and the derivation of functional inequalities and continuity results relevant to this generalized framework. In  Section~\ref{sec:DOIs as Special Cases of the Spectral Mapping Theorem for Continuous Spectrum Operators} , we reinterpret classical double operator integrals as special cases of the spectral mapping theorem tailored for continuous spectrum operators. This viewpoint lays the theoretical foundation for generalizing DOI beyond the self-adjoint scenario. Section~\ref{sec:Generalized Double Operator Integrals for Continuous Spectrum Operators}  introduces the formal definition of GDOI for CSOs and investigates their fundamental algebraic properties. This development is the cornerstone of the paper, marking the transition from classical DOI to the generalized framework. In Section~\ref{sec:Perturbation Formula} , we derive a perturbation formula for GDOI involving continuous spectrum operators, which generalizes known results from the self-adjoint setting and highlights the flexibility of the GDOI formulation. Section~\ref{sec:Inequalities Related to GDOI for Continuous Spectrum Operators}  explores several norm inequalities, including Lipschitz-type estimates, further solidifying the analytical robustness of the GDOI framework. In  Section~\ref{sec:Continuity of GDOI for Continuous Spectrum Operators}, we establish continuity properties of GDOI with respect to both the function $\beta$ and the involved operators, ensuring that GDOI behave well under limits and perturbations. Finally,  in Section~\ref{sec:Differentiation of Functions with Arguments as Continuous Spectrum Operators}, we  apply the developed theory to derive formulas for the differentiation of operator-valued functions, demonstrating the utility of GDOI in functional calculus. Overall, the paper’s key contributions lie in extending the DOI theory from self-adjoint operators to a more general setting involving continuous spectrum operators, thereby introducing the generalized DOI (GDOI) framework and establishing its foundational properties, estimates, and applications.

\textbf{Remark 1:} The current work can be considered a companion to our previous study, in which we focused exclusively on finite-dimensional operators (matrices), as presented in ~\cite{chang2025GDOIFinite}.

\textbf{Remark 2:} Also, note that by utilizing the spectral mapping theorem for hybrid spectrum operators from Section 5 of ~\cite{chang2024operatorChar} and ~\cite{chang2025GDOIFinite}, the current work can be naturally extended to accommodate hybrid spectrum operators.

\section{DOIs as Special Cases of the Spectral Mapping Theorem for Continuous Spectrum Operators}\label{sec:DOIs as Special Cases of the Spectral Mapping Theorem for Continuous Spectrum Operators}

We begin with a short review of the conventional DOI definition with continuous spectrum operators. Given a function $\beta: \mathbb{R}^2 \rightarrow \mathbb{C}$, two self-adjoint operators $\bm{X}_1, \bm{X}_2$, and any operator $\bm{Y}$. From the spectral mapping theorem of continuous spectrum operators~\cite{chang2024operatorChar}, we have
\begin{eqnarray}\label{eq0:  conv DOI def}
\bm{X}_1&=&\int\limits_{\lambda_1\in \sigma(\bm{X}_1)}\lambda_1 d\bm{E}_{\bm{X}_1}(\lambda_1);\nonumber \\
\bm{X}_2&=&\int\limits_{\lambda_2\in \sigma(\bm{X}_2)}\lambda_2 d\bm{E}_{\bm{X}_2}(\lambda_2);\nonumber \\
\bm{Y}&=&\int\limits_{\lambda_3 \in \sigma(\bm{Y})}\lambda_3 d\bm{E}_{\bm{Y}}(\lambda_3)+
\int\limits_{\lambda_3 \in \sigma(\bm{Y})}\left(\bm{Y}-\lambda_3\bm{I}\right)d\bm{E}_{\bm{Y}}(\lambda_3),
\end{eqnarray}
where $\lambda_1, \lambda_2$ and $\lambda_3$ are eigenvalues of the operators $\bm{X}_1, \bm{X}_2$ and $\bm{Y}$, respectively; $\sigma(\bm{X}_1), \sigma(\bm{X}_2)$ and $\sigma(\bm{Y})$ are spectrums of the operators $\bm{X}_1, \bm{X}_2$ and $\bm{Y}$, respectively; $\bm{E}_{\bm{X}_1}(\lambda_1), \bm{E}_{\bm{X}_2}(\lambda_2)$ and $\bm{E}_{\bm{Y}}(\lambda_3)$ are spectrum measures with respect to eigenvalues $\lambda_1, \lambda_2$ and $\lambda_3$, respectively. 

The conventional DOI definition for continuous spectrum operators, e.g., DOI on Schatten classes, is given by~\cite{skripka2019multilinear}:
\begin{eqnarray}\label{eq1:  conv DOI def}
T_{\beta}^{\bm{X}_1,\bm{X}_2}(\bm{Y})\define\int\limits_{\lambda_1\in \sigma(\bm{X}_1)}\int\limits_{\lambda_2\in \sigma(\bm{X}_2)}\beta(\lambda_1,\lambda_2)d\bm{E}_{\bm{X}_1}(\lambda_1)\bm{Y}d\bm{E}_{\bm{X}_2}(\lambda_2).
\end{eqnarray}
From the decomposition of the operator $\bm{Y}$ given by Eq.~\eqref{eq0:  conv DOI def}, we can express Eq.~\eqref{eq1:  conv DOI def} via the decomposition of the operator $\bm{Y}$ as 
\begin{eqnarray}\label{eq2:  conv DOI def}
\lefteqn{T_{\beta}^{\bm{X}_1,\bm{X}_2}(\bm{Y})}\nonumber \\
&=&\int\limits_{\lambda_1\in \sigma(\bm{X}_1)}\int\limits_{\lambda_2\in \sigma(\bm{X}_2)}\beta(\lambda_1,\lambda_2)d\bm{E}_{\bm{X}_1}(\lambda_1)\left(\int\limits_{\lambda_3 \in \sigma(\bm{Y})}\lambda_3 d\bm{E}_{\bm{Y}}(\lambda_3) \right. \nonumber \\
&&
\left.+
\int\limits_{\lambda_3 \in \sigma(\bm{Y})}\left(\bm{Y}-\lambda_3\bm{I}\right)d\bm{E}_{\bm{Y}}(\lambda_3)\right)d\bm{E}_{\bm{X}_2}(\lambda_2)\nonumber \\
&=&\int\limits_{\lambda_1\in \sigma(\bm{X}_1)}\int\limits_{\lambda_2\in \sigma(\bm{X}_2)}\int\limits_{\lambda_3\in \sigma(\bm{X}_3)}\beta(\lambda_1,\lambda_2)\lambda_3 d\bm{E}_{\bm{X}_1}(\lambda_1)d\bm{E}_{\bm{X}_3}(\lambda_3)d\bm{E}_{\bm{X}_2}(\lambda_2)\nonumber \\
&&+\int\limits_{\lambda_1\in \sigma(\bm{X}_1)}\int\limits_{\lambda_2\in \sigma(\bm{X}_2)}\int\limits_{\lambda_3\in \sigma(\bm{X}_3)}\beta(\lambda_1,\lambda_2)d\bm{E}_{\bm{X}_1}(\lambda_1)\left(\bm{Y}-\lambda_3\bm{I}\right)d\bm{E}_{\bm{Y}}(\lambda_3)d\bm{E}_{\bm{X}_2}(\lambda_2)
\end{eqnarray}

We begin by recalling Theorem 11 from~\cite{chang2024operatorChar}. Prior to stating the theorem, we introduce several notational conventions that are pertinent to its formulation. 

Given $r$ positive integers $q_1, q_2, \ldots, q_r$, we define $\alpha_{\kappa}(q_1, \ldots, q_r)$ to be the selection of these $r$ arguments $q_1, \ldots, q_r$ into $\kappa$ arguments, i.e., we have
\begin{equation}
\alpha_{\kappa}(q_1, \ldots, q_r) = \{ q_{\iota_1}, q_{\iota_2}, \ldots, q_{\iota_\kappa} \}.
\end{equation}

We use $\mathrm{Ind}(\alpha_{\kappa}(q_1, \ldots, q_r))$ to obtain the indices of those $\kappa$ positive integers $\{ q_{\iota_1}, q_{\iota_2}, \ldots, q_{\iota_\kappa} \}$, i.e., we have
\begin{equation}
\mathrm{Ind}(\alpha_{\kappa}(q_1, \ldots, q_r)) = \{ \iota_1, \iota_2, \ldots, \iota_\kappa \}.
\end{equation}

We use $\alpha_{\kappa}(q_1, \ldots, q_r) = 1$ to represent $q_{\iota_1} = 1, q_{\iota_2} = 1, \ldots, q_{\iota_\kappa} = 1$. 

We also use
\[
m_{\lambda_{\mbox{Ind}(\alpha_{\kappa}(q_1,\ldots,q_r))}}-1
\]
to represent
\[
m_{\lambda_{\iota_1}}-1,\ldots,m_{\lambda_{\iota_\kappa}}-1,
\]
where $m_{\lambda_{\iota_j}}$ is the order of the nilpotent  $\left(\bm{X}_{\iota_j}-\lambda_{\iota_j}\bm{I}\right)d\bm{E}_{\bm{X}_{\iota_j}}(\lambda_{\iota_j})$, i.e., 
\[
\left[\left(\bm{X}_{\iota_j}-\lambda_{\iota_j}\bm{I}\right)d\bm{E}_{\bm{X}_{\iota_j}}(\lambda_{\iota_j})\right]^\ell= \bm{0}, \quad \text{for } \ell \geq m_{\lambda_{\iota_j}}  \text{ and } j = 1, 2, \ldots, \kappa.
\]

Theorem 11 from~\cite{chang2024operatorChar} is given below. 
\begin{theorem}\label{thm: Spectral Mapping Theorem for r Variables inf}
Given an analytic function $f(z_1,z_2,\ldots,z_r)$ within the domain for $|z_l| < R_l$, and the operator $\bm{X}_l$ decomposed by:
\begin{eqnarray}\label{eq1-1: thm: Spectral Mapping Theorem for r Variables inf}
\bm{X}_l&=&\int\limits_{\lambda_l \in \sigma(\bm{X}_l)}\lambda_l d\bm{E}_{\bm{X}_l}(\lambda_l)+
\int\limits_{\lambda_l \in \sigma(\bm{X}_l)}\left(\bm{X}_l-\lambda_l\bm{I}\right)d\bm{E}_{\bm{X}_l}(\lambda_l),
\end{eqnarray}
where $\left\vert\lambda_{l}\right\vert<R_l$ for $l=1,2,\ldots,r$.

Then, we have
\begin{eqnarray}\label{eq2: thm: Spectral Mapping Theorem for kappa Variables inf}
\lefteqn{f(\bm{X}_1,\ldots,\bm{X}_r)=}\nonumber \\
&&\int\limits_{\lambda_1 \in \sigma(\bm{X}_1)}\cdots\int\limits_{\lambda_r \in \sigma(\bm{X}_r)}
f(\lambda_1,\ldots,\lambda_r)d\bm{E}_{\bm{X}_1}(\lambda_1)\cdots d\bm{E}_{\bm{X}_r}(\lambda_r) \nonumber \\
&&+\int\limits_{\lambda_1 \in \sigma(\bm{X}_1)}\cdots\int\limits_{\lambda_r \in \sigma(\bm{X}_r)}\sum\limits_{\kappa=1}^{r-1}\sum\limits_{\alpha_\kappa(q_1,\ldots,q_r)}\Bigg(\sum\limits_{\alpha_{\kappa}(q_1,\ldots,q_r)=1}^{m_{\lambda_{\mbox{Ind}(\alpha_{\kappa}(q_1,\ldots,q_r))}}-1}\nonumber \\
&&~~~~~ \frac{f^{\alpha_{\kappa}(q_1,\ldots,q_r)}(\lambda_1,\ldots,\lambda_r)}{q_{\iota_1}!q_{\iota_2}!\ldots q_{\iota_\kappa}!}\times \prod\limits_{\substack{\beta' =\mbox{Ind}(\alpha_{\kappa}(q_1,\ldots,q_r)), \bm{Y}=\left(\bm{X}_{\beta'} - \lambda_{\beta'}\bm{I}\right)^{q_{\beta'}}d\bm{E}_{\bm{X}_{\beta'}}(\lambda_{\beta'}) \\ \beta' \neq \mbox{Ind}(\alpha_{\kappa}(q_1,\ldots,q_r)), \bm{Y}=d\bm{E}_{\bm{X}_{\beta'}}(\lambda_{\beta'})}
}^{r} \bm{Y}\Bigg) 
\nonumber \\
&&+\int\limits_{\lambda_1 \in \sigma(\bm{X}_1)}\cdots\int\limits_{\lambda_r \in \sigma(\bm{X}_r)}\sum\limits_{q_1=\ldots=q_r=1}^{m_{\lambda_1}-1,\ldots,m_{\lambda_r}-1}
\frac{f^{(q_1,\ldots,q_r)}(\lambda_1,\ldots,\lambda_r)}{q_1!\cdots q_r!}\nonumber \\
&&\times \left(\bm{X}_1 - \lambda_1\bm{I}\right)^{q_1}d\bm{E}_{\bm{X}_1}(\lambda_1) \left(\bm{X}_2 - \lambda_2\bm{I}\right)^{q_2}d\bm{E}_{\bm{X}_2}(\lambda_2)\cdots \left(\bm{X}_r - \lambda_r\bm{I}\right)^{q_r}d\bm{E}_{\bm{X}_r}(\lambda_r),
\end{eqnarray}
where we have
\begin{itemize}
\item $\sum\limits_{\alpha_\kappa(q_1,\ldots,q_r)}$ runs over all $\kappa$ selections of $q_1,\ldots,q_r$ by $\alpha_\kappa(q_1,\ldots,q_r)$;
\item $m_{\lambda_{\mbox{Ind}(\alpha_{\kappa}(q_1,\ldots,q_r))}}-1 =$ $m_{\lambda_{\iota_1}}-1$,$\ldots,m_{\lambda_{\iota_\kappa}}-1$;
\item $f^{\alpha_{\kappa}(q_1,\ldots,q_r)}(\lambda_1,\ldots,\lambda_r)$ represents the partial derivatives with respect to variables with indices $\iota_1,\iota_2,\ldots,\iota_\kappa$ and the orders of derivatives given by $q_{\iota_1},q_{\iota_2},\ldots,q_{\iota_\kappa}$.
\end{itemize}
\end{theorem}

By setting $r = 3$ and $f(z_1, z_2, z_3) = \beta(z_1, z_3) z_2$, where $z_1 = \lambda_{1}$, $z_2 = \lambda_{3}$, and $z_3 = \lambda_{2}$, in Theorem~\ref{thm: Spectral Mapping Theorem for r Variables inf}, Eq.~\eqref{eq2: thm: Spectral Mapping Theorem for kappa Variables inf} can be reduced to Eq.~\eqref{eq2:  conv DOI def} because
\begin{eqnarray}
\frac{\partial f(z_1, z_2, z_3)}{\partial z_2} &=& \beta(z_1, z_3); \nonumber \\
\frac{\partial^\ell f(z_1, z_2, z_3)}{\partial^\ell z_2} &=& 0,
\end{eqnarray}
where $\ell > 1$. Therefore, we have
\begin{eqnarray} \label{eq2-1:  conv DOI def}
T_{\beta}^{\bm{X}_1, \bm{X}_2}(\bm{Y}) &=& f(\bm{X}_1, \bm{Y}, \bm{X}_2).
\end{eqnarray}

If we set $r = 3$ and $f(z_1, z_2, z_3) =  z_1\beta(z_2, z_3)$, where $z_1 = \lambda_{3}$, $z_2 = \lambda_{1}$, and $z_3 = \lambda_{2}$, in Theorem~\ref{thm: Spectral Mapping Theorem for r Variables inf}, we obtain a first variation of the conventional DOI by changing the position of the variable matrix $\bm{Y}$ in the DOI as
\begin{eqnarray} \label{eq3:  conv DOI def}
T_{\beta}^{\prime, \bm{X}_1, \bm{X}_2}(\bm{Y}) \define \int\limits_{\lambda_1 \in \sigma(\bm{X}_1)}\int\limits_{\lambda_2 \in \sigma(\bm{X}_2)}\beta(\lambda_{1}, \lambda_{2}) \bm{Y}d\bm{E}_{\bm{X}_1}(\lambda_1)d\bm{E}_{\bm{X}_2}(\lambda_2),
\end{eqnarray}
because
\begin{eqnarray}
\frac{\partial f(z_1, z_2, z_3)}{\partial z_1} &=& \beta(z_{2},z_{3}); \nonumber \\
\frac{\partial^\ell f(z_1, z_2, z_3)}{\partial^\ell z_1} &=& 0,
\end{eqnarray}
where $\ell > 1$. Therefore, we have
\begin{eqnarray} \label{eq3-1:  conv DOI def}
T_{\beta}^{\prime, \bm{X}_1, \bm{X}_2}(\bm{Y}) &=& f(\bm{Y}, \bm{X}_1, \bm{X}_2).
\end{eqnarray}

Moreover, if we set $r = 3$ and $f(z_1, z_2, z_3) = \beta(z_1, z_2) z_3$, where $z_1 = \lambda_{1}$, $z_2 = \lambda_{2}$, and $z_3 = \lambda_{3}$, in Theorem~\ref{thm: Spectral Mapping Theorem for r Variables inf}, we obtain a first variation of the conventional DOI by changing the position of the variable matrix $\bm{Y}$ in the DOI as
\begin{eqnarray} \label{eq3:  conv DOI def}
T_{\beta}^{'', \bm{X}_1, \bm{X}_2}(\bm{Y}) \define \int\limits_{\lambda_1 \in \sigma(\bm{X}_1)}\int\limits_{\lambda_2 \in \sigma(\bm{X}_2)}\beta(\lambda_{1}, \lambda_{2}) d\bm{E}_{\bm{X}_1}(\lambda_1)d\bm{E}_{\bm{X}_2}(\lambda_2)\bm{Y},
\end{eqnarray}
because
\begin{eqnarray}
\frac{\partial f(z_1, z_2, z_3)}{\partial z_3} &=& \beta(z_{1}, z_{2}); \nonumber \\
\frac{\partial^\ell f(z_1, z_2, z_3)}{\partial^\ell z_3} &=& 0,
\end{eqnarray}
where $\ell > 1$. Therefore, we have
\begin{eqnarray} \label{eq4-1:  conv DOI def}
T_{\beta}^{'', \bm{X}_1, \bm{X}_2}(\bm{Y}) &=& f(\bm{X}_1, \bm{X}_2, \bm{Y}).
\end{eqnarray}

From Eq.~\eqref{eq2-1: conv DOI def}, Eq.~\eqref{eq3-1: conv DOI def}, and Eq.~\eqref{eq4-1: conv DOI def}, it follows that the application of the multivariable operator spectral mapping theorem, as discussed in~\cite{chang2024operatorChar}, enables an extension of the conventional DOI definition by considering different positional arrangements of the input matrix $\bm{Y}$ relative to $\bm{X}_1$ and $\bm{X}_2$.

\section{Generalized Double Operator Integrals for Continuous Spectrum Operators}\label{sec:Generalized Double Operator Integrals for Continuous Spectrum Operators}

The generalized double operator integrals (GDOI) for continuous spectrum operators will be defined in Section~\ref{sec:GDOI for Continuous Spectrum Operators}. Besides, the algebraic properties of GDOI will be explored in Section~\ref{sec:Algebraic Properties}. 

\subsection{GDOI for Continuous Spectrum Operators}\label{sec:GDOI for Continuous Spectrum Operators}

Let us recall Thoerem 10 in~\cite{chang2024operatorChar}, which is given below.
\begin{theorem}\label{thm: Spectral Mapping Theorem for Two Variables inf}
Given an analytic function $f(z_1,z_2)$ within the domain for $|z_1| < R_1$ and $|z_2| < R_2$, the first operator $\bm{X}_1$ decomposed by:
\begin{eqnarray}\label{eq1-1: thm: Spectral Mapping Theorem for Two Variables inf}
\bm{X}_1&=&\int\limits_{\lambda_1 \in \sigma(\bm{X}_1)}\lambda_1 d\bm{E}_{\bm{X}_1}(\lambda_1)+
\int\limits_{\lambda_1 \in \sigma(\bm{X}_1)}\left(\bm{X}_1-\lambda_1\bm{I}\right)d\bm{E}_{\bm{X}_1}(\lambda_1),
\end{eqnarray}
where $\left\vert\lambda_{1}\right\vert<R_1$, and second operator $\bm{X}_2$ decomposed by:
\begin{eqnarray}\label{eq1-2: thm: Spectral Mapping Theorem for Two Variables inf}
\bm{X}_2&=&\int\limits_{\lambda_2 \in \sigma(\bm{X}_2)}\lambda_2 d\bm{E}_{\bm{X}_2}(\lambda_2)+
\int\limits_{\lambda_2 \in \sigma(\bm{X}_2)}\left(\bm{X}_2-\lambda_2\bm{I}\right)d\bm{E}_{\bm{X}_2}(\lambda_2),
\end{eqnarray}
where $\left\vert\lambda_{2}\right\vert<R_2$.

Then, we have
\begin{eqnarray}\label{eq2: thm: Spectral Mapping Theorem for Two Variables inf}
f(\bm{X}_1, \bm{X}_2)&=&\int\limits_{\lambda_1 \in \sigma(\bm{X}_1)}\int\limits_{\lambda_2 \in \sigma(\bm{X}_2)}f(\lambda_{1}, \lambda_{2})d\bm{E}_{\bm{X}_1}(\lambda_1)d\bm{E}_{\bm{X}_2}(\lambda_2) \nonumber \\
&&+\int\limits_{\lambda_1 \in \sigma(\bm{X}_1)}\int\limits_{\lambda_2 \in \sigma(\bm{X}_2)}\sum_{q_2=1}^{m_{\lambda_2}-1}\frac{f^{(-,q_2)}(\lambda_{1},\lambda_{2})}{q_2!}d\bm{E}_{\bm{X}_1}(\lambda_1)\left(\bm{X}_2-\lambda_2\bm{I}\right)^{q_2}d\bm{E}_{\bm{X}_2}(\lambda_2) \nonumber \\
&&+\int\limits_{\lambda_1 \in \sigma(\bm{X}_1)}\int\limits_{\lambda_2 \in \sigma(\bm{X}_2)}\sum_{q_1=1}^{m_{\lambda_1}-1}\frac{f^{(q_1,-)}(\lambda_{1},\lambda_{2})}{q_1!}\left(\bm{X}_1-\lambda_1\bm{I}\right)^{q_1}d\bm{E}_{\bm{X}_1}(\lambda_1)d\bm{E}_{\bm{X}_2}(\lambda_2)\nonumber \\
&&+\int\limits_{\lambda_1 \in \sigma(\bm{X}_1)}\int\limits_{\lambda_2 \in \sigma(\bm{X}_2)}\sum_{q_1=1}^{m_{\lambda_1}-1}\sum_{q_2=1}^{m_{\lambda_2}-1}\nonumber \\
&&\frac{f^{(q_1,q_2)}(\lambda_{1},\lambda_{2})}{q_1!q_2!}\left(\bm{X}_1-\lambda_1\bm{I}\right)^{q_1}d\bm{E}_{\bm{X}_1}(\lambda_1)\left(\bm{X}_2-\lambda_2\bm{I}\right)^{q_2}d\bm{E}_{\bm{X}_2}(\lambda_2).
\end{eqnarray}
\end{theorem}

According to Theorem~\ref{thm: Spectral Mapping Theorem for Two Variables inf}, the GDOI for continuous spectrum operators, denoted by $T_{\beta}^{\bm{X}_1, \bm{X}_2}(\bm{Y})$, can be defined as 
\begin{eqnarray}\label{eq1:  GDOI def}
T_{\beta}^{\bm{X}_1, \bm{X}_2}(\bm{Y})&\define&\int\limits_{\lambda_1 \in \sigma(\bm{X}_1)}\int\limits_{\lambda_2 \in \sigma(\bm{X}_2)}\beta(\lambda_{1}, \lambda_{2})d\bm{E}_{\bm{X}_1}(\lambda_1)\bm{Y}d\bm{E}_{\bm{X}_2}(\lambda_2) \nonumber \\
&&+\int\limits_{\lambda_1 \in \sigma(\bm{X}_1)}\int\limits_{\lambda_2 \in \sigma(\bm{X}_2)}\sum_{q_2=1}^{m_{\lambda_2}-1}\frac{f^{(-,q_2)}(\lambda_{1},\lambda_{2})}{q_2!}d\bm{E}_{\bm{X}_1}(\lambda_1)\bm{Y}\left(\bm{X}_2-\lambda_2\bm{I}\right)^{q_2}d\bm{E}_{\bm{X}_2}(\lambda_2) \nonumber \\
&&+\int\limits_{\lambda_1 \in \sigma(\bm{X}_1)}\int\limits_{\lambda_2 \in \sigma(\bm{X}_2)}\sum_{q_1=1}^{m_{\lambda_1}-1}\frac{f^{(q_1,-)}(\lambda_{1},\lambda_{2})}{q_1!}\left(\bm{X}_1-\lambda_1\bm{I}\right)^{q_1}d\bm{E}_{\bm{X}_1}(\lambda_1)\bm{Y}d\bm{E}_{\bm{X}_2}(\lambda_2)\nonumber \\
&&+\int\limits_{\lambda_1 \in \sigma(\bm{X}_1)}\int\limits_{\lambda_2 \in \sigma(\bm{X}_2)}\sum_{q_1=1}^{m_{\lambda_1}-1}\sum_{q_2=1}^{m_{\lambda_2}-1}\nonumber \\
&&\frac{f^{(q_1,q_2)}(\lambda_{1},\lambda_{2})}{q_1!q_2!}\left(\bm{X}_1-\lambda_1\bm{I}\right)^{q_1}d\bm{E}_{\bm{X}_1}(\lambda_1)\bm{Y}\left(\bm{X}_2-\lambda_2\bm{I}\right)^{q_2}d\bm{E}_{\bm{X}_2}(\lambda_2).
\end{eqnarray}
From Eq.~\eqref{eq1: GDOI def}, when the operators $\bm{X}_1$ and $\bm{X}_2$ are self-adjoint, the operator $T_{\beta}^{\bm{X}_1,\bm{X}_2}(\bm{Y})$ reduces to the conventional DOI form given in Eq.~\eqref{eq1: conv DOI def}.

\subsection{Algebraic Properties}\label{sec:Algebraic Properties}

We will explore the algebraic properties for the GDOI defined by Eq.~\eqref{eq1:  GDOI def}. Because we have the following decompositions for the operators $\bm{X}_1$ and $\bm{X}_2$:
\begin{eqnarray}\label{eq:X1 decomp p and n parts}
\bm{X}_1&=&\int\limits_{\lambda_1 \in \sigma(\bm{X}_1)}\lambda_1 d\bm{E}_{\bm{X}_1}(\lambda_1)+
\int\limits_{\lambda_1 \in \sigma(\bm{X}_1)}\left(\bm{X}_1-\lambda_1\bm{I}\right)d\bm{E}_{\bm{X}_1}(\lambda_1)\nonumber \\
&\define&\bm{X}_{1,\bm{P}}+\bm{X}_{1,\bm{N}},
\end{eqnarray}
where $\left\vert\lambda_{1}\right\vert<R_1$, and the second operator $\bm{X}_2$ is decomposed by:
\begin{eqnarray}\label{eq:X2 decomp p and n parts}
\bm{X}_2&=&\int\limits_{\lambda_2 \in \sigma(\bm{X}_2)}\lambda_2 d\bm{E}_{\bm{X}_2}(\lambda_2)+
\int\limits_{\lambda_2 \in \sigma(\bm{X}_2)}\left(\bm{X}_2-\lambda_2\bm{I}\right)d\bm{E}_{\bm{X}_2}(\lambda_2)\nonumber \\
&\define&\bm{X}_{2,\bm{P}}+\bm{X}_{2,\bm{N}},
\end{eqnarray}
where $\left\vert\lambda_{2}\right\vert<R_2$. 

From decompositions given by Eq.~\eqref{eq:X1 decomp p and n parts} and Eq.~\eqref{eq:X2 decomp p and n parts}, we have the following decomposition proposition based on the GDOI definition given by Eq.~\eqref{eq1:  GDOI def}. 
\begin{proposition}
Given operators $\bm{X}_1$ and $\bm{X}_2$, which are decomposed as Eq.~\eqref{eq:X1 decomp p and n parts} and Eq.~\eqref{eq:X2 decomp p and n parts}, respectively, then, we have
\begin{eqnarray}
T_{\beta}^{\bm{X}_1,\bm{X}_2}(\bm{Y})&=&T_{\beta}^{\bm{X}_{1,\bm{P}},\bm{X}_{2,\bm{P}}}(\bm{Y})+T_{\beta}^{\bm{X}_{1,\bm{P}},\bm{X}_{2,\bm{N}}}(\bm{Y})+T_{\beta}^{\bm{X}_{1,\bm{N}},\bm{X}_{2,\bm{P}}}(\bm{Y})+T_{\beta}^{\bm{X}_{1,\bm{N}},\bm{X}_{2,\bm{N}}}(\bm{Y}).
\end{eqnarray}
\end{proposition}

In the following, we concentrate on the algebraic properties related to the function $\beta$. We start with Lemma~\ref{lma: ind of Proj and Nilp}, which concerns the independence of the product between the projection part and the nilpotent part of operators. Note that we use the subscript $\bm{P}$ and the subscript $\bm{N}$ in Eq.~\eqref{eq1-1: thm: Spectral Mapping Theorem for Two Variables inf} and Eq.~\eqref{eq1-2: thm: Spectral Mapping Theorem for Two Variables inf} to represent the projection part and the nilpotent part, respectively. 

\begin{lemma}\label{lma: ind of Proj and Nilp}
Let $\bm{X}_1$ and $\bm{X}_2$ are two operators with spectral decomposition given by Eq.~\eqref{eq1-1: thm: Spectral Mapping Theorem for Two Variables inf} and Eq.~\eqref{eq1-2: thm: Spectral Mapping Theorem for Two Variables inf}, respectively. We have the linear independence of the following four categories of operators:
\begin{align}\label{eq1: lma: ind of Proj and Nilp}
1.\quad & d\bm{E}_{\bm{X}_1}(\lambda_1)\bm{Y}d\bm{E}_{\bm{X}_2}(\lambda_2); \nonumber \\
2.\quad & d\bm{E}_{\bm{X}_1}(\lambda_1) \bm{Y}\left(\bm{X}_2-\lambda_2\bm{I}\right)^{q_2}d\bm{E}_{\bm{X}_2}(\lambda_2); \nonumber \\  
3.\quad & \left(\bm{X}_1-\lambda_1\bm{I}\right)^{q_1}d\bm{E}_{\bm{X}_1}(\lambda_1)\bm{Y}d\bm{E}_{\bm{X}_2}(\lambda_2); \nonumber \\ 
4.\quad & \left(\bm{X}_1-\lambda_1\bm{I}\right)^{q_1}d\bm{E}_{\bm{X}_1}(\lambda_1) \bm{Y}\left(\bm{X}_2-\lambda_2\bm{I}\right)^{q_2}d\bm{E}_{\bm{X}_2}(\lambda_2).
\end{align}

where $\bm{Y} \neq \bm{0}$, $1 \leq q_1 < m_{\lambda_1}$ and $1 \leq q_2 < m_{\lambda_2}$.
\end{lemma} 
\textbf{Proof:}
We define the set $\mathcal{S}_{\bm{X}_1, \bm{X}_2, \bm{Y}}$ of operators as
\begin{eqnarray}\label{eq2: lma: ind of Proj and Nilp}
\mathcal{S}_{\bm{X}_1, \bm{X}_2, \bm{Y}}&\define&\Bigg\{
\int\limits_{\lambda_1 \in \sigma(\bm{X}_1)}\int\limits_{\lambda_2 \in \sigma(\bm{X}_2)}c^{PP}_{\lambda_1,\lambda_2}d\bm{E}_{\bm{X}_1}(\lambda_1)\bm{Y}d\bm{E}_{\bm{X}_2}(\lambda_2)
\nonumber \\
&&+\int\limits_{\lambda_1 \in \sigma(\bm{X}_1)}\int\limits_{\lambda_2 \in \sigma(\bm{X}_2)}\sum_{q_2=1}^{m_{\lambda_2}-1}c^{PN}_{\lambda_1,\lambda_2,q_2}d\bm{E}_{\bm{X}_1}(\lambda_1)\bm{Y}\left(\bm{X}_2-\lambda_2\bm{I}\right)^{q_2}d\bm{E}_{\bm{X}_2}(\lambda_2)\nonumber \\
&&+ \int\limits_{\lambda_1 \in \sigma(\bm{X}_1)}\int\limits_{\lambda_2 \in \sigma(\bm{X}_2)}\sum_{q_1=1}^{m_{\lambda_1}-1}c^{NP}_{\lambda_1,\lambda_2,q_1}\left(\bm{X}_1-\lambda_1\bm{I}\right)^{q_1}d\bm{E}_{\bm{X}_1}(\lambda_1)\bm{Y}d\bm{E}_{\bm{X}_2}(\lambda_2)\nonumber \\
&&+\int\limits_{\lambda_1 \in \sigma(\bm{X}_1)}\int\limits_{\lambda_2 \in \sigma(\bm{X}_2)}\sum_{q_1=1}^{m_{\lambda_1}-1}\sum_{q_2=1}^{m_{\lambda_2}-1}c^{NN}_{\lambda_1,\lambda_2,q_1,q_2}\left(\bm{X}_1-\lambda_1\bm{I}\right)^{q_1}d\bm{E}_{\bm{X}_1}(\lambda_1)\nonumber \\
&&~~~~\times \bm{Y}\left(\bm{X}_2-\lambda_2\bm{I}\right)^{q_2}d\bm{E}_{\bm{X}_2}(\lambda_2)\Bigg\},
\end{eqnarray}

where $c^{PP}_{\lambda_1,\lambda_2}$, $c^{PN}_{\lambda_1,\lambda_2,q_2}$, $c^{NP}_{\lambda_1,\lambda_2,q_1}$ and $c^{NN}_{\lambda_1,\lambda_2,q_1,q_2}$ are complex scalers.  

For $j=1,2$, we have
\begin{eqnarray}
d\bm{E}_{\bm{X}_j}(\lambda_j)d\bm{E}_{\bm{X}_j}(\lambda_j')&=&d\bm{E}_{\bm{X}_j}(\lambda_j)\delta(\lambda_j,\lambda_j'), \nonumber 
\end{eqnarray}
\begin{eqnarray}
(\bm{X}_j-\lambda_j\bm{I})d\bm{E}_{\bm{X}_j}(\lambda_j)(\bm{X}_j-\lambda_j'\bm{I})d\bm{E}_{\bm{X}_j}(\lambda_j')=(\bm{X}_j-\lambda_j\bm{I})^2d\bm{E}_{\bm{X}_j}(\lambda_j)\delta(\lambda_j,\lambda_j'),\nonumber 
\end{eqnarray} 
\begin{eqnarray}\label{eq3: lma: ind of Proj and Nilp}
d\bm{E}_{\bm{X}_j}(\lambda_j')(\bm{X}_j-\lambda_j\bm{I})d\bm{E}_{\bm{X}_j}(\lambda_j)&=&(\bm{X}_j-\lambda_j\bm{I})d\bm{E}_{\bm{X}_j}(\lambda)d\bm{E}_{\bm{X}_j}(\lambda_j') \nonumber \\
&=&(\bm{X}_j-\lambda_j\bm{I})d\bm{E}_{\bm{X}_j}(\lambda_j)\delta(\lambda_j,\lambda_j').
\end{eqnarray} 
therefore, \(d\bm{E}_{\bm{X}_j}(\lambda_j)\) and \(\left(\bm{X}_j-\lambda_j\bm{I}\right)d\bm{E}_{\bm{X}_j}(\lambda_j)\) act on different generalized eigenspaces, their behaviors are distinct. Then, we have
\begin{itemize}
\item The term \(d\bm{E}_{\bm{X}_1}(\lambda_1)\bm{Y}d\bm{E}_{\bm{X}_2}(\lambda_2)\) isolates components where both left and right transformations remain in the eigenspaces.
\item The term \(d\bm{E}_{\bm{X}_1}(\lambda_1) \bm{Y}\left(\bm{X}_2-\lambda_2\bm{I}\right)^{q_2}d\bm{E}_{\bm{X}_2}(\lambda_2)\) involves a right multiplication by a nilpotent matrix, affecting only part of the generalized eigenspace.
\item Similarly, \( \left(\bm{X}_1-\lambda_1\bm{I}\right)^{q_1}d\bm{E}_{\bm{X}_1}(\lambda_1)\bm{Y}d\bm{E}_{\bm{X}_2}(\lambda_2)\) applies nilpotent transformations on the left.
\item The term \(\left(\bm{X}_1-\lambda_1\bm{I}\right)^{q_1}d\bm{E}_{\bm{X}_1}(\lambda_1) \bm{Y}\left(\bm{X}_2-\lambda_2\bm{I}\right)^{q_2}d\bm{E}_{\bm{X}_2}(\lambda_2)\) applies nilpotent transformations on both sides.
\end{itemize}
Since nilpotent operators act non-trivially within their respective Jordan blocks, while projection operators restrict transformations to specific eigenspaces, these terms span distinct transformation spaces when $\bm{Y} \neq \bm{0}$. Therefore, they cannot be expressed as linear combinations of one another, which establishes that the four categories of operators given in Eq.~\eqref{eq1: lma: ind of Proj and Nilp} are linearly independent.
$\hfill\Box$

If we are given a domain $\mathcal{D} \in \mathbb{C}^2$ with $(\lambda_{1}, \lambda_{2})\in \mathcal{D}$ with a bi-variable analytic function defined over $\mathcal{D}$, namely $\beta(z_1,z_2)$,  we use $\mathcal{S}_{\bm{X}_1,\bm{X}_2,\bm{Y}}(\mathcal{D})$ to represent the following set:
\begin{eqnarray}\label{eq0: lma: GDOI linear homomorphism}
\mathcal{S}_{\bm{X}_1, \bm{X}_2, \bm{Y}}(\mathcal{D})&\define&\Bigg\{
\int\limits_{\lambda_1 \in \sigma(\bm{X}_1)}\int\limits_{\lambda_2 \in \sigma(\bm{X}_2)}c^{PP}_{\lambda_1,\lambda_2}d\bm{E}_{\bm{X}_1}(\lambda_1)\bm{Y}d\bm{E}_{\bm{X}_2}(\lambda_2)
\nonumber \\
&&+\int\limits_{\lambda_1 \in \sigma(\bm{X}_1)}\int\limits_{\lambda_2 \in \sigma(\bm{X}_2)}\sum_{q_2=1}^{m_{\lambda_2}-1}c^{PN}_{\lambda_1,\lambda_2,q_2}d\bm{E}_{\bm{X}_1}(\lambda_1)\bm{Y}\left(\bm{X}_2-\lambda_2\bm{I}\right)^{q_2}d\bm{E}_{\bm{X}_2}(\lambda_2)\nonumber \\
&&+ \int\limits_{\lambda_1 \in \sigma(\bm{X}_1)}\int\limits_{\lambda_2 \in \sigma(\bm{X}_2)}\sum_{q_1=1}^{m_{\lambda_1}-1}c^{NP}_{\lambda_1,\lambda_2,q_1}\left(\bm{X}_1-\lambda_1\bm{I}\right)^{q_1}d\bm{E}_{\bm{X}_1}(\lambda_1)\bm{Y}d\bm{E}_{\bm{X}_2}(\lambda_2)\nonumber \\
&&+\int\limits_{\lambda_1 \in \sigma(\bm{X}_1)}\int\limits_{\lambda_2 \in \sigma(\bm{X}_2)}\sum_{q_1=1}^{m_{\lambda_1}-1}\sum_{q_2=1}^{m_{\lambda_2}-1}c^{NN}_{\lambda_1,\lambda_2,q_1,q_2}\left(\bm{X}_1-\lambda_1\bm{I}\right)^{q_1}d\bm{E}_{\bm{X}_1}(\lambda_1)\nonumber \\
&&~~~~\times \bm{Y}\left(\bm{X}_2-\lambda_2\bm{I}\right)^{q_2}d\bm{E}_{\bm{X}_2}(\lambda_2)\Bigg\},
\end{eqnarray}
such that $c^{PP}_{\lambda_1,\lambda_2}$, $c^{PN}_{\lambda_1,\lambda_2,q_2}$, $c^{NP}_{\lambda_1,\lambda_2,q_1}$ and $c^{NN}_{\lambda_1,\lambda_2,q_1,q_2}$ satisfy the following:
 \begin{eqnarray}\label{eq0-1: lma: GDOI linear homomorphism}
c^{PP}_{\lambda_1,\lambda_2}&=&\beta(\lambda_{1},\lambda_{2});\nonumber \\
c^{PN}_{\lambda_1,\lambda_2,q_2}&=&\frac{\beta^{(-,q_2)}(\lambda_{1},\lambda_{2})}{q_2!};\nonumber \\
c^{NP}_{\lambda_1,\lambda_2,q_1}&=&\frac{\beta^{(q_1,-)}(\lambda_{1},\lambda_{2})}{q_1!};\nonumber \\
c^{NN}_{\lambda_1,\lambda_2,q_2}&=&\frac{\beta^{(q_1,q_2)}(\lambda_{1},\lambda_{2})}{q_1!q_2!}.
\end{eqnarray}

The composition between  $T^{\bm{X}_1,\bm{X}_2}_\beta$ and $T^{\bm{X}_1,\bm{X}_2}_\gamma$ is represented by the symbol $\circ$, and is defined by 
\begin{eqnarray}
T^{\bm{X}_1,\bm{X}_2}_\beta \circ T^{\bm{X}_1,\bm{X}_2}_\gamma \define T^{\bm{X}_1,\bm{X}_2}_\beta(T^{\bm{X}_1,\bm{X}_2}_\gamma).
\end{eqnarray}

From the composition relation, we have the following Lemma~\ref{lma: GDOI linear homomorphism} about linear homomorphism property of the operator $T^{\bm{X}_1,\bm{X}_2}_\beta$. 
\begin{lemma}\label{lma: GDOI linear homomorphism}
Let $\bm{X}_1$ and $\bm{X}_2$ be operators and consider the operator $T^{\bm{X}_1,\bm{X}_2}_\beta$, which is defined by Eq.~\eqref{eq1:  GDOI def}. Given $\beta(z_1,z_2)$ as a bi-variable analytic function over the domain $\mathcal{D}$, the mapping $f: \beta \rightarrow T^{\bm{X}_1,\bm{X}_2}_{\beta}$ is a linear homomorphism. 
\end{lemma}
\textbf{Proof:}
To prove that the mapping $f$ is a linear homomorphism, the following properties have to be estblished:
\begin{eqnarray}\label{eq1-1: lma: GDOI linear homomorphism}
T^{\bm{X}_1,\bm{X}_2}_{c_1\beta_1+c_2\beta_2}&=&c_1 T^{\bm{X}_1,\bm{X}_2}_{\beta_1}+ c_2 T^{\bm{X}_1,\bm{X}_2}_{\beta_2},
\end{eqnarray}
where $c_1$ and $c_2$ are two scalers; and 
\begin{eqnarray}\label{eq1-2: lma: GDOI linear homomorphism}
T^{\bm{X}_1,\bm{X}_2}_{\beta\gamma}&=&T^{\bm{X}_1,\bm{X}_2}_{\beta}\circ T^{\bm{X}_1,\bm{X}_2}_{\gamma}.
\end{eqnarray}

\textbf{Proof of Eq.~\eqref{eq1-1: lma: GDOI linear homomorphism}}

Because we have
\begin{eqnarray}\label{eq2: lma: GDOI linear homomorphism}
\lefteqn{T_{c_1\beta_1+c_2\beta_2}^{\bm{X}_1,\bm{X}_2}(\bm{Y})}\nonumber \\
&=_1&\int\limits_{\lambda_1 \in \sigma(\bm{X}_1)}\int\limits_{\lambda_2 \in \sigma(\bm{X}_2)}[c_1\beta_1(\lambda_{1}, \lambda_{2})+c_2\beta_2(\lambda_{1}, \lambda_{2})]d\bm{E}_{\bm{X}_1}(\lambda_1)\bm{Y}d\bm{E}_{\bm{X}_2}(\lambda_2) \nonumber \\
&&+\int\limits_{\lambda_1 \in \sigma(\bm{X}_1)}\int\limits_{\lambda_2 \in \sigma(\bm{X}_2)}\sum_{q_2=1}^{m_{\lambda_2}-1}\frac{[c_1\beta_1^{(-,q_2)}(\lambda_{1},\lambda_{2})+c_2\beta_2^{(-,q_2)}(\lambda_{1},\lambda_{2})]}{q_2!}\nonumber \\
&& \times d\bm{E}_{\bm{X}_1}(\lambda_1)\bm{Y}\left(\bm{X}_2-\lambda_2\bm{I}\right)^{q_2}d\bm{E}_{\bm{X}_2}(\lambda_2) \nonumber \\
&&+\int\limits_{\lambda_1 \in \sigma(\bm{X}_1)}\int\limits_{\lambda_2 \in \sigma(\bm{X}_2)}\sum_{q_1=1}^{m_{\lambda_1}-1}\frac{[c_1\beta_1^{(q_1,-)}(\lambda_{1},\lambda_{2}) + c_2 \beta_2^{(q_1,-)}(\lambda_{1},\lambda_{2})]}{q_1!}\nonumber \\
&& \times\left(\bm{X}_1-\lambda_1\bm{I}\right)^{q_1}d\bm{E}_{\bm{X}_1}(\lambda_1)\bm{Y}d\bm{E}_{\bm{X}_2}(\lambda_2)  \nonumber \\
&&+\int\limits_{\lambda_1 \in \sigma(\bm{X}_1)}\int\limits_{\lambda_2 \in \sigma(\bm{X}_2)}\sum_{q_1=1}^{m_{\lambda_1}-1}\sum_{q_2=1}^{m_{\lambda_2}-1}\frac{[c_1\beta_1^{(q_1,q_2)}(\lambda_{1},\lambda_{2}) + c_2\beta_2^{(q_1,q_2)}(\lambda_{1},\lambda_{2})]}{q_1!q_2!}\nonumber \\
&& \times\left(\bm{X}_1-\lambda_1\bm{I}\right)^{q_1}d\bm{E}_{\bm{X}_1}(\lambda_1)\bm{Y}\left(\bm{X}_2-\lambda_2\bm{I}\right)^{q_2}d\bm{E}_{\bm{X}_2}(\lambda_2)\nonumber \\
&=&c_1 T^{\bm{X}_1,\bm{X}_2}_{\beta_1}(\bm{Y})+ c_2 T^{\bm{X}_1,\bm{X}_2}_{\beta_2}(\bm{Y}),
\end{eqnarray}
where we apply the linearity of the partial derivative in $=_1$. Then, Eq.~\eqref{eq1-1: lma: GDOI linear homomorphism} is established.

\textbf{Proof of Eq.~\eqref{eq1-2: lma: GDOI linear homomorphism}}

The term $T^{\bm{X}_1,\bm{X}_2}_{\beta\gamma}$ can be expressed by
\begin{eqnarray}\label{eq3: lma: GDOI linear homomorphism}
\lefteqn{T_{\beta\gamma}^{\bm{X}_1,\bm{X}_2}(\bm{Y})}\nonumber \\
&=_1&\underbrace{\int\limits_{\lambda_1 \in \sigma(\bm{X}_1)}\int\limits_{\lambda_2 \in \sigma(\bm{X}_2)}[\beta(\lambda_{1}, \lambda_{2})\gamma(\lambda_{1}, \lambda_{2})]d\bm{E}_{\bm{X}_1}(\lambda_1)\bm{Y}d\bm{E}_{\bm{X}_2}(\lambda_2)}_{\mbox{Part 1}} \nonumber \\
&&+\underbrace{\int\limits_{\lambda_1 \in \sigma(\bm{X}_1)}\int\limits_{\lambda_2 \in \sigma(\bm{X}_2)}\sum_{q_2=1}^{m_{\lambda_2}-1}\frac{[\beta(\lambda_{1}, \lambda_{2})\gamma(\lambda_{1}, \lambda_{2})]^{(-,q_2)}}{q_2!}d\bm{E}_{\bm{X}_1}(\lambda_1)\bm{Y}\left(\bm{X}_2-\lambda_2\bm{I}\right)^{q_2}d\bm{E}_{\bm{X}_2}(\lambda_2) }_{\mbox{Part 2}} \nonumber \\
&&+\underbrace{\int\limits_{\lambda_1 \in \sigma(\bm{X}_1)}\int\limits_{\lambda_2 \in \sigma(\bm{X}_2)}\sum_{q_1=1}^{m_{\lambda_1}-1}\frac{[\beta(\lambda_{1}, \lambda_{2})\gamma(\lambda_{1}, \lambda_{2})]^{(q_1,-)}}{q_1!}\left(\bm{X}_1-\lambda_1\bm{I}\right)^{q_1}d\bm{E}_{\bm{X}_1}(\lambda_1)\bm{Y}d\bm{E}_{\bm{X}_2}(\lambda_2)}_{\mbox{Part 3}} \nonumber \\
&&+\int\limits_{\lambda_1 \in \sigma(\bm{X}_1)}\int\limits_{\lambda_2 \in \sigma(\bm{X}_2)}\sum_{q_1=1}^{m_{\lambda_1}-1}\sum_{q_2=1}^{m_{\lambda_2}-1}\frac{[\beta(\lambda_{1}, \lambda_{2})\gamma(\lambda_{1}, \lambda_{2})]^{(q_1,q_2)}}{q_1!q_2!}\nonumber\\
&& \underbrace{ \times \left(\bm{X}_1-\lambda_1\bm{I}\right)^{q_1}d\bm{E}_{\bm{X}_1}(\lambda_1)\bm{Y}\left(\bm{X}_2-\lambda_2\bm{I}\right)^{q_2}d\bm{E}_{\bm{X}_2}(\lambda_2)}_{\mbox{Part 4}} .
\end{eqnarray}

On the other hand, we can express $T^{\bm{X}_1,\bm{X}_2}_{\beta}\circ T^{\bm{X}_1,\bm{X}_2}_{\gamma}$ as
\begin{eqnarray}\label{eq4: lma: GDOI linear homomorphism}
 \lefteqn{T^{\bm{X}_1,\bm{X}_2}_{\beta}\circ T^{\bm{X}_1,\bm{X}_2}_{\gamma}=T^{\bm{X}_1,\bm{X}_2}_{\beta}( T^{\bm{X}_1,\bm{X}_2}_{\gamma})}\nonumber \\
&=&\int\limits_{\lambda_1 \in \sigma(\bm{X}_1)}\int\limits_{\lambda_2 \in \sigma(\bm{X}_2)}\beta(\lambda_{1}, \lambda_{2})d\bm{E}_{\bm{X}_1}(\lambda_1)[T^{\bm{X}_1,\bm{X}_2}_{\gamma}(\bm{Y})]d\bm{E}_{\bm{X}_2}(\lambda_2) \nonumber \\
&&+\int\limits_{\lambda_1 \in \sigma(\bm{X}_1)}\int\limits_{\lambda_2 \in \sigma(\bm{X}_2)}\sum_{q_2=1}^{m_{\lambda_2}-1}\frac{\beta^{(-,q_2)}(\lambda_{1},\lambda_{2})}{q_2!}d\bm{E}_{\bm{X}_1}(\lambda_1)[T^{\bm{X}_1,\bm{X}_2}_{\gamma}(\bm{Y})]\left(\bm{X}_2-\lambda_2\bm{I}\right)^{q_2}d\bm{E}_{\bm{X}_2}(\lambda_2) \nonumber \\
&&+\int\limits_{\lambda_1 \in \sigma(\bm{X}_1)}\int\limits_{\lambda_2 \in \sigma(\bm{X}_2)}\sum_{q_1=1}^{m_{\lambda_1}-1}\frac{\beta^{(q_1,-)}(\lambda_{1},\lambda_{2})}{q_1!}\left(\bm{X}_1-\lambda_1\bm{I}\right)^{q_1}d\bm{E}_{\bm{X}_1}(\lambda_1)[T^{\bm{X}_1,\bm{X}_2}_{\gamma}(\bm{Y})]d\bm{E}_{\bm{X}_2}(\lambda_2)  \nonumber \\
&&+\int\limits_{\lambda_1 \in \sigma(\bm{X}_1)}\int\limits_{\lambda_2 \in \sigma(\bm{X}_2)}\sum_{q_1=1}^{m_{\lambda_1}-1}\sum_{q_2=1}^{m_{\lambda_2}-1}\frac{\beta^{(q_1,q_2)}(\lambda_{1},\lambda_{2})}{q_1!q_2!}\nonumber \\
&& \times \left(\bm{X}_1-\lambda_1\bm{I}\right)^{q_1}d\bm{E}_{\bm{X}_1}(\lambda_1)[T^{\bm{X}_1,\bm{X}_2}_{\gamma}(\bm{Y})]\left(\bm{X}_2-\lambda_2\bm{I}\right)^{q_2}d\bm{E}_{\bm{X}_2}(\lambda_2).
\end{eqnarray}
 
We arrange the expansion of R.H.S. of Eq.~\eqref{eq3: lma: GDOI linear homomorphism} according to 
\begin{align}
1.\quad & d\bm{E}_{\bm{X}_1}(\lambda_1)\bm{Y}d\bm{E}_{\bm{X}_2}(\lambda_2), \nonumber \\
2.\quad & d\bm{E}_{\bm{X}_1}(\lambda_1)\bm{Y}\left(\bm{X}_2-\lambda_2\bm{I}\right)^{q_2}d\bm{E}_{\bm{X}_2}(\lambda_2),\nonumber \\
3.\quad &\left(\bm{X}_1-\lambda_1\bm{I}\right)^{q_1}d\bm{E}_{\bm{X}_1}(\lambda_1)\bm{Y}d\bm{E}_{\bm{X}_1}(\lambda_1),\nonumber \\
4.\quad &\left(\bm{X}_1-\lambda_1\bm{I}\right)^{q_1}d\bm{E}_{\bm{X}_1}(\lambda_1)\bm{Y}\left(\bm{X}_2-\lambda_2\bm{I}\right)^{q_2}d\bm{E}_{\bm{X}_2}(\lambda_2),
\end{align}
 to obtain the following four expressions :
\begin{align}\label{eq4-1: lma: GDOI linear homomorphism}
\int\limits_{\lambda_1 \in \sigma(\bm{X}_1)}\int\limits_{\lambda_2 \in \sigma(\bm{X}_2)}[\beta(\lambda_{1}, \lambda_{2})\gamma(\lambda_{1}, \lambda_{2})]d\bm{E}_{\bm{X}_1}(\lambda_1)\bm{Y}d\bm{E}_{\bm{X}_2}(\lambda_2), 
\end{align}
\begin{align}\label{eq4-2: lma: GDOI linear homomorphism}
\int\limits_{\lambda_1 \in \sigma(\bm{X}_1)}\int\limits_{\lambda_2 \in \sigma(\bm{X}_2)}\sum_{q_2=1}^{m_{\lambda_2}-1}\frac{\beta(\lambda_{1},\lambda_{2})\gamma^{(-,q_2)}(\lambda_{1},\lambda_{2})}{q_2!}\nonumber \\ \times d\bm{E}_{\bm{X}_1}(\lambda_1)\bm{Y}\left(\bm{X}_2-\lambda_2\bm{I}\right)^{q_2}d\bm{E}_{\bm{X}_2}(\lambda_2)\nonumber \\
+\int\limits_{\lambda_1 \in \sigma(\bm{X}_1)}\int\limits_{\lambda_2 \in \sigma(\bm{X}_2)}\sum_{q_2=1}^{m_{\lambda_2}-1}\frac{\beta^{(-,q_2)}(\lambda_{1},\lambda_{2})\gamma(\lambda_{1},\lambda_{2})}{q_2!}\nonumber \\ \times d\bm{E}_{\bm{X}_1}(\lambda_1)\bm{Y}\left(\bm{X}_2-\lambda_2\bm{I}\right)^{q_2}d\bm{E}_{\bm{X}_2}(\lambda_2)\nonumber \\
+\int\limits_{\lambda_1 \in \sigma(\bm{X}_1)}\int\limits_{\lambda_2 \in \sigma(\bm{X}_2)}\sum_{q'_2=1, q''_2=1}^{q'_2+q''_2 = m_{\lambda_2}-1}\frac{\beta^{(-,q'_2)}(\lambda_{1},\lambda_{2})}{q'_2!}\frac{\gamma^{(-,q''_2)}(\lambda_{1},\lambda_{2})}{q''_2!} \nonumber 
\\ \times d\bm{E}_{\bm{X}_1}(\lambda_1)\bm{Y}\left(\bm{X}_2-\lambda_2\bm{I}\right)^{q'_2+q''_2}d\bm{E}_{\bm{X}_2}(\lambda_2),
\end{align}
\begin{eqnarray}\label{eq4-3: lma: GDOI linear homomorphism}
\int\limits_{\lambda_1 \in \sigma(\bm{X}_1)}\int\limits_{\lambda_2 \in \sigma(\bm{X}_2)}\sum_{q_1=1}^{m_{\lambda_1}-1}\frac{\beta(\lambda_{1},\lambda_{2})\gamma^{(q_1,-)}(\lambda_{1},\lambda_{2})}{q_1!}\nonumber \\
\times \left(\bm{X}_1-\lambda_1\bm{I}\right)^{q_1}d\bm{E}_{\bm{X}_1}(\lambda_1)\bm{Y}d\bm{E}_{\bm{X}_2}(\lambda_2)\nonumber \\
+\int\limits_{\lambda_1 \in \sigma(\bm{X}_1)}\int\limits_{\lambda_2 \in \sigma(\bm{X}_2)}\sum_{q_1=1}^{m_{\lambda_1}-1}\frac{\beta^{(q_1,-)}(\lambda_{1},\lambda_{2})\gamma(\lambda_{1},\lambda_{2})}{q_1!}\nonumber \\
\times \left(\bm{X}_1-\lambda_1\bm{I}\right)^{q_1}d\bm{E}_{\bm{X}_1}(\lambda_1)\bm{Y}d\bm{E}_{\bm{X}_2}(\lambda_2)\nonumber \\
+\int\limits_{\lambda_1 \in \sigma(\bm{X}_1)}\int\limits_{\lambda_2 \in \sigma(\bm{X}_2)}\sum_{q'_1=1, q''_1=1}^{q'_1+q''_1 = m_{\lambda_1}-1}\frac{\beta^{(q'_1,-)}(\lambda_{1},\lambda_{2})}{q'_1!}\frac{\gamma^{(q''_1,-)}(\lambda_{1},\lambda_{2})}{q''_1!}\nonumber \\
\times\left(\bm{X}_1-\lambda_1\bm{I}\right)^{q'_1+q''_1}d\bm{E}_{\bm{X}_1}(\lambda_1)\bm{Y}d\bm{E}_{\bm{X}_2}(\lambda_2),    
\end{eqnarray}
\begin{eqnarray}\label{eq4-4: lma: GDOI linear homomorphism}
\sum\limits_{q_1=1}^{m_{\lambda_1}}\sum\limits_{q_2=1}^{m_{\lambda_2}}\frac{\beta^{(q_1,q_2)}(\lambda_{1},\lambda_{2})\gamma(\lambda_{1},\lambda_{2})}{q_1! q_2!}\left(\bm{X}_1-\lambda_1\bm{I}\right)^{q_1}d\bm{E}_{\bm{X}_1}(\lambda_1)\bm{Y}\left(\bm{X}_2-\lambda_2\bm{I}\right)^{q_2}d\bm{E}_{\bm{X}_2}(\lambda_2)\nonumber \\
+\sum\limits_{q_1=1}^{m_{\lambda_1}}\sum\limits_{q_2=1}^{m_{\lambda_2}}\frac{\beta^{(q_1,-)}(\lambda_{1},\lambda_{2})\gamma^{(-, q_2)}(\lambda_{1},\lambda_{2})}{q_1! q_2!}\left(\bm{X}_1-\lambda_1\bm{I}\right)^{q_1}d\bm{E}_{\bm{X}_1}(\lambda_1)\bm{Y}\left(\bm{X}_2-\lambda_2\bm{I}\right)^{q_2}d\bm{E}_{\bm{X}_2}(\lambda_2)\nonumber \\
+\sum\limits_{q_1=1}^{m_{\lambda_1}}\sum\limits_{q_2=1}^{m_{\lambda_2}}\frac{\beta^{(-,q_2)}(\lambda_{1},\lambda_{2})\gamma^{(q_1,-)}(\lambda_{1},\lambda_{2})}{q_1! q_2!}\left(\bm{X}_1-\lambda_1\bm{I}\right)^{q_1}d\bm{E}_{\bm{X}_1}(\lambda_1)\bm{Y}\left(\bm{X}_2-\lambda_2\bm{I}\right)^{q_2}d\bm{E}_{\bm{X}_2}(\lambda_2)\nonumber \\
+\sum\limits_{q_1=1}^{m_{\lambda_1}}\sum\limits_{q_2=1}^{m_{\lambda_2}}\frac{\beta(\lambda_{1},\lambda_{2})\gamma^{(q_1,q_2)}(\lambda_{1},\lambda_{2})}{q_1! q_2!}\left(\bm{X}_1-\lambda_1\bm{I}\right)^{q_1}d\bm{E}_{\bm{X}_1}(\lambda_1)\bm{Y}\left(\bm{X}_2-\lambda_2\bm{I}\right)^{q_2}d\bm{E}_{\bm{X}_2}(\lambda_2)\nonumber \\
+\sum\limits_{q_1=1}^{m_{\lambda_1}}\sum\limits_{q'_2=1, q''_2=1}^{q'_2 + q''_2=m_{\lambda_2}-1}\frac{\beta^{(q_1, q'_2)}(\lambda_{1},\lambda_{2})\gamma^{(-, q''_2)}(\lambda_{1},\lambda_{2})}{q_1! q'_2! q''_2!}\nonumber \\
\times \left(\bm{X}_1-\lambda_1\bm{I}\right)^{q_1}d\bm{E}_{\bm{X}_1}(\lambda_1)\bm{Y}\left(\bm{X}_2-\lambda_2\bm{I}\right)^{q'_2+q''_2}d\bm{E}_{\bm{X}_2}(\lambda_2)\nonumber \\
+\sum\limits_{q_1=1}^{m_{\lambda_1}}\sum\limits_{q'_2=1, q''_2=1}^{q'_2 + q''_2=m_{\lambda_2}-1}\frac{\beta^{(-, q'_2)}(\lambda_{1},\lambda_{2})\gamma^{(q_1, q''_2)}(\lambda_{1},\lambda_{2}) }{q_1! q'_2! q''_2!}\nonumber \\
\times \left(\bm{X}_1-\lambda_1\bm{I}\right)^{q_1}d\bm{E}_{\bm{X}_1}(\lambda_1)\bm{Y}\left(\bm{X}_2-\lambda_2\bm{I}\right)^{q'_2+q''_2}d\bm{E}_{\bm{X}_2}(\lambda_2)\nonumber \\
+\sum\limits_{q'_1=1, q''_1=1}^{q'_1 + q''_1=m_{\lambda_1}-1}\sum\limits_{q_2=1}^{m_{\lambda_2}}\frac{\beta^{(q'_1, q_2)}(\lambda_{1},\lambda_{2})\gamma^{(q''_1, -)}(\lambda_{1},\lambda_{2})}{q'_1! q_2! q''_1!}\nonumber \\ 
\times \left(\bm{X}_1-\lambda_1\bm{I}\right)^{q'_1+q''_1}d\bm{E}_{\bm{X}_1}(\lambda_1)\bm{Y}\left(\bm{X}_2-\lambda_2\bm{I}\right)^{q_2}d\bm{E}_{\bm{X}_2}(\lambda_2)\nonumber \\
+\sum\limits_{q'_1=1, q''_1=1}^{q'_1 + q''_1=m_{\lambda_1}-1}\sum\limits_{q_2=1}^{m_{\lambda_2}}\frac{\beta^{(q'_1, -)}(\lambda_{1},\lambda_{2})\gamma^{(q''_1,q_2)}(\lambda_{1},\lambda_{2})}{q'_1! q_2! q''_1!}\nonumber \\
\times \left(\bm{X}_1-\lambda_1\bm{I}\right)^{q'_1+q''_1}d\bm{E}_{\bm{X}_1}(\lambda_1)\bm{Y}\left(\bm{X}_2-\lambda_2\bm{I}\right)^{q_2}d\bm{E}_{\bm{X}_2}(\lambda_2)\nonumber \\
+\sum\limits_{q'_1=1, q''_1=1}^{q'_1 + q''_1=m_{\lambda_1}-1}\sum\limits_{q'_2=1, q''_2=1}^{q'_2 + q''_2=m_{\lambda_2}-1}\frac{\beta^{(q'_1, q'_2)}(\lambda_{1},\lambda_{2})\gamma^{(q''_1,q''_2)}(\lambda_{1},\lambda_{2})}{q'_1! q'_2! q''_1!q''_2!}\nonumber \\
\left(\bm{X}_1-\lambda_1\bm{I}\right)^{q'_1+q''_1}d\bm{E}_{\bm{X}_1}(\lambda_1)\bm{Y}\left(\bm{X}_2-\lambda_2\bm{I}\right)^{q'_2+q''_2}d\bm{E}_{\bm{X}_2}(\lambda_2).
\end{eqnarray}

The Leibniz rule for partial derivatives of a product of two differentiable functions \( f(x, y) \) and \( g(x, y) \) is:
\begin{eqnarray}\label{eq5: lma: GDOI linear homomorphism}
(f(x, y)g(x, y))^{(q_1, q_2)} = \sum_{k_1=0}^{q_1} \sum_{k_2=0}^{q_2} \binom{q_1}{k_1} \binom{q_2}{k_2} f^{(k_1, k_2)}(x, y) g^{(q_1 - k_1, q_2 - k_2)}(x, y),
\end{eqnarray}
where \( f^{(k_1, k_2)} \) denotes the partial derivative of \( f(x, y) \) with respect to \( x \) for \( k_1 \) times and with respect to \( y \) for \( k_2 \) times. 

It is evident that Part 1 in Eq.~\eqref{eq3: lma: GDOI linear homomorphism} is identical to Eq.~\eqref{eq4-1: lma: GDOI linear homomorphism}. As for the remaining Parts 2, 3, and 4 in Eq.~\eqref{eq3: lma: GDOI linear homomorphism}, they correspond exactly to Eq.~\eqref{eq4-2: lma: GDOI linear homomorphism}, Eq.~\eqref{eq4-3: lma: GDOI linear homomorphism}, and Eq.~\eqref{eq4-4: lma: GDOI linear homomorphism}, respectively, by applying the Leibniz rule for partial derivatives stated in Eq.~\eqref{eq5: lma: GDOI linear homomorphism}. Consequently, Eq.~\eqref{eq1-2: lma: GDOI linear homomorphism} holds true as well.
$\hfill\Box$

Below, we will present Lemma~\ref{lma:Analytic Function Identity via Derivative Agreement at Finite Points}, which will be used to show the injective property of GDOI.  

\begin{lemma}[Analytic Function Identity via Derivative Agreement at Finite Points]\label{lma:Analytic Function Identity via Derivative Agreement at Finite Points}
Let \(\beta, \gamma \in \mathcal{O}(\mathcal{D})\) be analytic bivaraite functions defined on a connected domain \(\mathcal{D} \subseteq \mathbb{C}^2\). Suppose there exists a finite set of distinct points  
\[
\mathcal{P} = \{ (\lambda_1^{(i)}, \lambda_2^{(i)}) \}_{i=1}^N \subset \sigma(\bm{X}_1) \times  \sigma(\bm{X}_2) \subseteq \mathcal{D}
\]  
and corresponding non-negative integers \((p_i, q_i)\) for \(i = 1, \ldots, N\), such that:
\[
\frac{\partial^{p+q} \beta}{\partial \lambda_1^p \partial \lambda_2^q}(\lambda_1^{(i)}, \lambda_2^{(i)}) = \frac{\partial^{p+q} \gamma}{\partial \lambda_1^p \partial \lambda_2^q}(\lambda_1^{(i)}, \lambda_2^{(i)}), \quad \forall \ 0 \leq p \leq p_i, \ 0 \leq q \leq q_i, \ \forall i.
\]  

Then, defining \(h := \beta - \gamma\), we have:
\[
h \in \mathcal{N}_{(p_1, q_1), \dots, (p_N, q_N)},
\]
where \(\mathcal{N}_{(p_1, q_1), \dots, (p_N, q_N)}\) is the space of analytic functions vanishing up to order \((p_i, q_i)\) at each \((\lambda_1^{(i)}, \lambda_2^{(i)})\).

Furthermore, if for each \(i\), the limits \(p_i, q_i \to \infty\), then:
\[
\beta \equiv \gamma \quad \text{on } \mathcal{D}.
\]
\end{lemma}
\textbf{Proof:}
The first step of this proof is to construct the vanishing function. Let \(h := \beta - \gamma\). Since \(\beta, \gamma \in \mathcal{O}(\mathcal{D})\), their difference \(h \in \mathcal{O}(\mathcal{D})\) is also analytic.

From the hypothesis:
\[
\frac{\partial^{p+q} h}{\partial \lambda_1^p \partial \lambda_2^q}(\lambda_1^{(i)}, \lambda_2^{(i)}) = 0, \quad \forall \ 0 \leq p \leq p_i, \ 0 \leq q \leq q_i, \ \forall i.
\]

This means \(h\) and all its partial derivatives up to order \((p_i, q_i)\) vanish at each point \((\lambda_1^{(i)}, \lambda_2^{(i)})\). Therefore:
\[
h \in \mathcal{N}_{(p_1, q_1), \dots, (p_N, q_N)},
\]
i.e., \(h\) belongs to the vanishing subspace defined by the filtration at those points.

The next step is to show the uniqueness via analytic continuation. Suppose now that \(p_i, q_i \to \infty\) for all \(i\). Then all mixed partial derivatives of \(h\) vanish at \((\lambda_1^{(i)}, \lambda_2^{(i)})\) of all orders. From complex analysis, a classical result states: \\

If a function \(f \in \mathcal{O}(\mathcal{D})\) has all partial derivatives vanish at a point \((z_1, z_2) \in \mathcal{D}\), then \(f \equiv 0\) in a neighborhood of that point. Analyticity implies that \(f\) vanishes identically on all of \(\mathcal{D}\) (by identity theorem for several complex variables).

Applying this result to \(h\), and noting that \(\mathcal{D}\) is connected, we conclude:
\[
h = \beta - \gamma \equiv 0 \quad \text{on } \mathcal{D} \quad \Rightarrow \quad \beta \equiv \gamma.
\]
$\hfill\Box$

The main objective of this section is to present Theorem~\ref{thm: GDOI algebraic properties}, which characterizes the algebraic structure of the operator \( T^{\bm{X}_1,\bm{X}_2}_\beta \).
\begin{theorem}\label{thm: GDOI algebraic properties}
Let \( \bm{X}_1 \) and \( \bm{X}_2 \) be arbitrary operators whose spectra lie within the connected domain \( \mathcal{D} \). Consider the operator \( T^{\bm{X}_1,\bm{X}_2}_\beta \) defined in Eq.~\eqref{eq1:  GDOI def}, where \( \beta(z_1, z_2) \) is a bi-variable analytic function such that all bi-variable functions $\beta$ defined over \( \mathcal{D} \) satisfy Lemma~\ref{lma:Analytic Function Identity via Derivative Agreement at Finite Points} conditions. Then the mapping
\[
f: \beta \mapsto T^{\bm{X}_1,\bm{X}_2}_\beta
\]
is a linear isomorphism provided that \( T^{\bm{X}_1,\bm{X}_2}_\beta(\bm{Y}) \in \mathcal{S}_{\bm{X}_1, \bm{X}_2, \bm{Y}}(\mathcal{D}) \).
\end{theorem}
\textbf{Proof:}
From Lemma~\ref{lma: GDOI linear homomorphism}, we know that the mapping is a linear homomorphism. This theorem is proved if the injective and surjective are satisfied by the mapping $f$.  

To show that the mapping $f$ is injective, we first assume that $T^{\bm{X}_1,\bm{X}_2}_{\beta}=T^{\bm{X}_1,\bm{X}_2}_{\gamma}$, i.e., $\bm{0} = T^{\bm{X}_1,\bm{X}_2}_{\beta}-T^{\bm{X}_1,\bm{X}_2}_{\gamma}$. Then, we have
\begin{eqnarray}
\lefteqn{\bm{0}=\int\limits_{\lambda_1 \in \sigma(\bm{X}_1)}\int\limits_{\lambda_2 \in \sigma(\bm{X}_2)}(\beta(\lambda_{1}, \lambda_{2})-\gamma(\lambda_{1}, \lambda_{2}))}\nonumber \\
&&~~\times d\bm{E}_{\bm{X}_1}(\lambda_1)\bm{Y}d\bm{E}_{\bm{X}_2}(\lambda_2)\nonumber \\
&&+\int\limits_{\lambda_1 \in \sigma(\bm{X}_1)}\int\limits_{\lambda_2 \in \sigma(\bm{X}_2)}\sum_{q_2=1}^{m_{\lambda_2}-1}\frac{\beta^{(-,q_2)}(\lambda_{1},\lambda_{2}) - \gamma^{(-,q_2)}(\lambda_{1},\lambda_{2})}{q_2!}\nonumber\\  
&&~~\times d\bm{E}_{\bm{X}_1}(\lambda_1)\bm{Y}\left(\bm{X}_2-\lambda_2\bm{I}\right)^{q_2}d\bm{E}_{\bm{X}_2}(\lambda_2)\nonumber \\
&&+\int\limits_{\lambda_1 \in \sigma(\bm{X}_1)}\int\limits_{\lambda_2 \in \sigma(\bm{X}_2)}\sum_{q_1=1}^{m_{\lambda_1}-1}\frac{\beta^{(q_1,-)}(\lambda_{1},\lambda_{2})- \gamma^{(q_1,-)}(\lambda_{1},\lambda_{2})}{q_1!}\nonumber\\  
&&~~\times  \left(\bm{X}_1-\lambda_1\bm{I}\right)^{q_1}d\bm{E}_{\bm{X}_1}(\lambda_1)\bm{Y}d\bm{E}_{\bm{X}_1}(\lambda_1)\nonumber \\
&&+\int\limits_{\lambda_1 \in \sigma(\bm{X}_1)}\int\limits_{\lambda_2 \in \sigma(\bm{X}_2)}\sum_{q_1=1}^{m_{\lambda_1}-1}\sum_{q_2=1}^{m_{\lambda_2}-1}\frac{\beta^{(q_1,q_2)}(\lambda_{1},\lambda_{2}) - \gamma^{(q_1,q_2)}(\lambda_{1},\lambda_{2})}{q_1!q_2!}\nonumber\\  
&&~~\times \left(\bm{X}_1-\lambda_1\bm{I}\right)^{q_1}d\bm{E}_{\bm{X}_1}(\lambda_1)\bm{Y}\left(\bm{X}_2-\lambda_2\bm{I}\right)^{q_2}d\bm{E}_{\bm{X}_2}(\lambda_2)
\end{eqnarray}
From Lemma~\ref{lma: ind of Proj and Nilp} and Lemma~\ref{lma:Analytic Function Identity via Derivative Agreement at Finite Points}, we have $\beta=\gamma$ in the domain $\mathcal{D}$. 

The surjective of the mapping $f$ is easy to check from the definition of $\mathcal{S}_{\bm{X}_1, \bm{X}_2, \bm{Y}}(\mathcal{D})$ given by Eq.~\eqref{eq0: lma: GDOI linear homomorphism} and Eq.~\eqref{eq0-1: lma: GDOI linear homomorphism}. 
$\hfill\Box$

\section{Perturbation Formula}\label{sec:Perturbation Formula}

We have the following Theorem~\ref{thm:pert formula} to determine the difference between $f(\bm{X}_1)$ and $f(\bm{X}_2)$ via GDOI.
\begin{theorem}\label{thm:pert formula}
Given an analytic function $f(z)$ within the domain for $|z| < R$, the first operator $\bm{X}_1$ can be expressed by
\begin{eqnarray}\label{eq1-1:thm:pert formula}
\bm{X}_1&=&\int\limits_{\lambda_1 \in \sigma(\bm{X}_1)}\lambda_1 d\bm{E}_{\bm{X}_1}(\lambda_1)+
\int\limits_{\lambda_1 \in \sigma(\bm{X}_1)}\left(\bm{X}_1-\lambda_1\bm{I}\right)d\bm{E}_{\bm{X}_1}(\lambda_1),
\end{eqnarray}
where $\left\vert\lambda_{1}\right\vert<R$, and second  operator $\bm{X}_2$ can be expressed by
\begin{eqnarray}\label{eq1-2:thm:pert formula}
\bm{X}_2&=&\int\limits_{\lambda_2 \in \sigma(\bm{X}_2)}\lambda_2 d\bm{E}_{\bm{X}_2}(\lambda_2)+
\int\limits_{\lambda_2 \in \sigma(\bm{X}_2)}\left(\bm{X}_2-\lambda_2\bm{I}\right)d\bm{E}_{\bm{X}_2}(\lambda_2)\nonumber \\
\end{eqnarray}
where $\left\vert\lambda_{2}\right\vert<R$. We also assume that $\lambda_{1} \neq \lambda_{2}$. Then, we have
\begin{eqnarray}\label{eq2:thm:pert formula}
f(\bm{X}_1)\bm{Y}-\bm{Y}f(\bm{X}_2)&=&T_{\frac{f(x_1) - f(x_2)}{x_1 -x_2}}^{\bm{X}_{1},\bm{X}_{2}}(\bm{X}_{1}\bm{Y}-\bm{Y}\bm{X}_{2}). 
\end{eqnarray}
Moreover, by setting $\bm{Y}=\bm{I}$, we have
\begin{eqnarray}\label{eq1:cor:pert formula I}
f(\bm{X}_1)-f(\bm{X}_2)&=&T_{\frac{f(x_1) - f(x_2)}{x_1 -x_2}}^{\bm{X}_{1},\bm{X}_{2}}(\bm{X}_{1}-\bm{X}_{2}). 
\end{eqnarray}
\end{theorem}
\textbf{Proof:}
We first define two projection functions $\pi_1(x_1,x_2)$ and $\pi_2(x_1,x_2)$ as
\begin{eqnarray}\label{eq3:thm:pert formula}
\pi_1(x_1,x_2)=x_1,~~\pi_2(x_1,x_2)&=&x_2.
\end{eqnarray}

Then, from GDOI definition given by Eq.~\eqref{eq1:  GDOI def}, we have
\begin{eqnarray}\label{eq4-1:thm:pert formula}
T_{\pi_1}^{\bm{X}_1,\bm{X}_2}(\bm{Y})&=&\int\limits_{\lambda_1 \in \sigma(\bm{X}_1)}\int\limits_{\lambda_2 \in \sigma(\bm{X}_2)}\lambda_{1}d\bm{E}_{\bm{X}_1}(\lambda_1)\bm{Y}d\bm{E}_{\bm{X}_2}(\lambda_2) \nonumber \\
&&+\int\limits_{\lambda_1 \in \sigma(\bm{X}_1)}\int\limits_{\lambda_2 \in \sigma(\bm{X}_2)}\left(\bm{X}_1-\lambda_1\bm{I}\right)d\bm{E}_{\bm{X}_1}(\lambda_1)\bm{Y}d\bm{E}_{\bm{X}_2}(\lambda_2) \nonumber \\
&=&\bm{X}_1 \bm{Y}.
\end{eqnarray}
Similarly, we have
\begin{eqnarray}\label{eq4-2:thm:pert formula}
T_{\pi_2}^{\bm{X}_1,\bm{X}_2}(\bm{Y})&=&\bm{Y}\bm{X}_2.
\end{eqnarray}
Besides, we also have 
\begin{eqnarray}\label{eq5-1:thm:pert formula}
T_{f \circ \pi_1}^{\bm{X}_1,\bm{X}_2}(\bm{Y})&=&\int\limits_{\lambda_1 \in \sigma(\bm{X}_1)}\int\limits_{\lambda_2 \in \sigma(\bm{X}_2)}f(\lambda_{1})d\bm{E}_{\bm{X}_1}(\lambda_1)\bm{Y}d\bm{E}_{\bm{X}_2}(\lambda_2)  \nonumber \\
&&+\int\limits_{\lambda_1 \in \sigma(\bm{X}_1)}\int\limits_{\lambda_2 \in \sigma(\bm{X}_2)}\sum_{q_1=1}^{m_{\lambda_1}-1}\frac{f^{(q_1)}(\lambda_{1})}{q_1!}\left(\bm{X}_1-\lambda_1\bm{I}\right)^{q_1}d\bm{E}_{\bm{X}_1}(\lambda_1)\bm{Y}d\bm{E}_{\bm{X}_2}(\lambda_2) \nonumber \\
&=&f(\bm{X}_1)\bm{Y}.
\end{eqnarray}
Similarly, we have
\begin{eqnarray}\label{eq5-2:thm:pert formula}
T_{f \circ \pi_2}^{\bm{X}_1,\bm{X}_2}(\bm{Y})&=&\bm{Y}f(\bm{X}_2).
\end{eqnarray}

Then, we have
\begin{eqnarray}\label{eq6:thm:pert formula}
f(\bm{X}_1)\bm{Y} - \bm{Y}f(\bm{X}_2)&=_1&T_{f \circ \pi_1}^{\bm{X}_1,\bm{X}_2}(\bm{Y}) - T_{f \circ \pi_2}^{\bm{X}_1,\bm{X}_2}(\bm{Y})\nonumber \\
&=&T_{f \circ \pi_1 - f \circ \pi_2}^{\bm{X}_1,\bm{X}_2}(\bm{Y})\nonumber \\
&=&T_{\frac{f(x_1) - f(x_2)}{x_1 - x_2} (\pi_1 - \pi_2)}^{\bm{X}_1,\bm{X}_2}(\bm{Y})\nonumber \\
&=&T_{\frac{f(x_1) - f(x_2)}{x_1 - x_2}\pi_1}^{\bm{X}_1,\bm{X}_2}(\bm{Y}) - T_{\frac{f(x_1) - f(x_2)}{x_1 - x_2}\pi_2}^{\bm{X}_1,\bm{X}_2}(\bm{Y})\nonumber \\
&=&T_{\frac{f(x_1) - f(x_2)}{x_1 - x_2}}^{\bm{X}_1,\bm{X}_2}(T_{\pi_1}^{\bm{X}_1,\bm{X}_2}(\bm{Y})) - T_{\frac{f(x_1) - f(x_2)}{x_1 - x_2}}^{\bm{X}_1,\bm{X}_2}(T_{\pi_2}^{\bm{X}_1,\bm{X}_2}(\bm{Y}))\nonumber \\
&=_2&T_{\frac{f(x_1) - f(x_2)}{x_1 - x_2}}^{\bm{X}_1,\bm{X}_2}(\bm{X}_1\bm{Y}-\bm{Y}\bm{X}_2),
\end{eqnarray}
where we apply Eq.~\eqref{eq5-1:thm:pert formula} and Eq.~\eqref{eq5-2:thm:pert formula} in $=_1$, we apply Eq.~\eqref{eq4-1:thm:pert formula} and Eq.~\eqref{eq4-2:thm:pert formula} in $=_2$ and other equalities come from Lemma~\ref{lma: GDOI linear homomorphism}.

Finally, Eq.~\eqref{eq1:cor:pert formula I} is obtained by setting $\bm{Y}$ as an identiy operator in Eq.~\eqref{eq2:thm:pert formula}.
$\hfill\Box$

We present the following Theorem~\ref{thm:pert formula diff DOI}, which characterizes the difference between the GDOI and the conventional DOI for self-adjoint operators, i.e., 
\[
T_{\frac{f(x_1) - f(x_2)}{x_1 - x_2}}^{\bm{X}_{1,\bm{P}}, \bm{X}_{2,\bm{P}}}(\bm{X}_{1,\bm{P}} - \bm{X}_{2,\bm{P}}),
\]
as described in Theorem~\ref{thm:pert formula}.

\begin{theorem}\label{thm:pert formula diff DOI}
Let $f(z)$ be an analytic function defined on the domain $|z| < R$. Suppose $\bm{X}_1$ is an $m \times m$ matrix with $K_1$ distinct eigenvalues $\lambda_{1}$ for $k_1 = 1, 2, \ldots, K_1$, such that
\begin{eqnarray}\label{eq1-1:thm:pert formula diff DOI}
\bm{X}_1&=&\int\limits_{\lambda_1 \in \sigma(\bm{X}_1)}\lambda_1 d\bm{E}_{\bm{X}_1}(\lambda_1)+
\int\limits_{\lambda_1 \in \sigma(\bm{X}_1)}\left(\bm{X}_1-\lambda_1\bm{I}\right)d\bm{E}_{\bm{X}_1}(\lambda_1)\nonumber \\
&\define&\bm{X}_{1,\bm{P}}+\bm{X}_{1,\bm{N}},
\end{eqnarray}
where $\left|\lambda_{1}\right| < R$. Similarly, let $\bm{X}_2$ be another operator, such that
\begin{eqnarray}\label{eq1-2:thm:pert formula diff DOI}
\bm{X}_2&=&\int\limits_{\lambda_2 \in \sigma(\bm{X}_2)}\lambda_2 d\bm{E}_{\bm{X}_2}(\lambda_2)+
\int\limits_{\lambda_2 \in \sigma(\bm{X}_2)}\left(\bm{X}_2-\lambda_2\bm{I}\right)d\bm{E}_{\bm{X}_2}(\lambda_2)\nonumber \\
&\define&\bm{X}_{2,\bm{P}}+\bm{X}_{2,\bm{N}},
\end{eqnarray}
with $\left|\lambda_{2}\right| < R$. Moreover, we assume that $\lambda_{1} \neq \lambda_{2}$. Then, the following identity holds:
\begin{eqnarray}\label{eq2:thm:pert formula diff DOI}
T_{\frac{f(x_1) - f(x_2)}{x_1 - x_2}}^{\bm{X}_1, \bm{X}_2}(\bm{X}_1 - \bm{X}_2) &=&
T_{\frac{f(x_1) - f(x_2)}{x_1 - x_2}}^{\bm{X}_{1,\bm{P}}, \bm{X}_{2,\bm{P}}}(\bm{X}_{1,\bm{P}} - \bm{X}_{2,\bm{P}}) \nonumber \\
&& + \left[T_{f(x_1)}^{\bm{X}_{1,\bm{N}}, \bm{X}_{2,\bm{P}}}(\bm{I}) - T_{f(x_2)}^{\bm{X}_{1,\bm{P}}, \bm{X}_{2,\bm{N}}}(\bm{I})\right].
\end{eqnarray}
\end{theorem}
\textbf{Proof:}
From Theorem 9 in~\cite{chang2024operatorChar}, we have 
\begin{eqnarray}\label{eq3-1:thm:pert formula diff DOI}
f(\bm{X}_1)&=&\int\limits_{\lambda_1 \in \sigma(\bm{X}_1)}f(\lambda_1)d\bm{E}_{\bm{X}_1}(\lambda_1)+\int\limits_{\lambda \in \sigma(\bm{X}_1)}\sum\limits_{q_1=1}^{m_{\lambda_1}-1}\frac{f^{(q_1)}(\lambda_1)}{q_1!}(\bm{X}_1-\lambda_1\bm{I})^{q_1} d\bm{E}_{\bm{X}_1}(\lambda_1),
\end{eqnarray}
and
\begin{eqnarray}\label{eq3-2:thm:pert formula diff DOI}
f(\bm{X}_2)&=&\int\limits_{\lambda_2 \in \sigma(\bm{X}_2)}f(\lambda_2)d\bm{E}_{\bm{X}_2}(\lambda_2)+\int\limits_{\lambda \in \sigma(\bm{X}_2)}\sum\limits_{q_2=1}^{m_{\lambda_2}-1}\frac{f^{(q_2)}(\lambda_2)}{q_2!}(\bm{X}_2-\lambda_2\bm{I})^{q_2} d\bm{E}_{\bm{X}_2}(\lambda_2),
\end{eqnarray}
Then, we can express $f(\bm{X}_1)-f(\bm{X}_2)$ as
\begin{eqnarray}\label{eq4:thm:pert formula diff DOI}
\lefteqn{f(\bm{X}_1)-f(\bm{X}_2)}\nonumber \\
&=&\underbrace{\left[\int\limits_{\lambda_1 \in \sigma(\bm{X}_1)}f(\lambda_{1})d\bm{E}_{\bm{X}_1}(\lambda_1)-\int\limits_{\lambda_2 \in \sigma(\bm{X}_2)}f(\lambda_{2})d\bm{E}_{\bm{X}_2}(\lambda_2)\right]}_{\mbox{Part I}}\nonumber \\
&&+\left[\int\limits_{\lambda_1 \in \sigma(\bm{X}_1)}\sum\limits_{q_1=1}^{m_{\lambda_1}-1}\frac{f^{(q_1)}(\lambda_{1})}{q_1!}(\bm{X}_1-\lambda_1\bm{I})^{q_1}d\bm{E}_{\bm{X}_1}(\lambda_1) \right. \nonumber \\
&& \left.-\int\limits_{\lambda_2 \in \sigma(\bm{X}_2)}\sum\limits_{q_2=1}^{m_{\lambda_2}-1}\frac{f^{(q_2)}(\lambda_{2})}{q_2!}(\bm{X}_2-\lambda_2\bm{I})^{q_2}d\bm{E}_{\bm{X}_2}(\lambda_2)\right]_{\mbox{Part II}}
\end{eqnarray}

From GDOI definition given by Eq.~\eqref{eq1:  GDOI def}, we have
\begin{eqnarray}\label{eq5:thm:pert formula diff DOI}
T_{\frac{f(x_1) - f(x_2)}{x_1 -x_2}}^{\bm{X}_{1,\bm{P}},\bm{X}_{2,\bm{P}}}(\bm{X}_{1,\bm{P}}-\bm{X}_{2,\bm{P}})&=&\int\limits_{\lambda_1 \in \sigma(\bm{X}_1)}\int\limits_{\lambda_2 \in \sigma(\bm{X}_2)}\frac{f(\lambda_{1}) - f(\lambda_{2})}{\lambda_{1} -\lambda_{2}}d\bm{E}_{\bm{X}_1}(\lambda_1)(\bm{X}_{1,\bm{P}}-\bm{X}_{2,\bm{P}})d\bm{E}_{\bm{X}_2}(\lambda_2) \nonumber \\
&=&\int\limits_{\lambda_1 \in \sigma(\bm{X}_1)}\int\limits_{\lambda_2 \in \sigma(\bm{X}_2)}\frac{f(\lambda_{1}) - f(\lambda_{2})}{\lambda_{1} -\lambda_{2}}d\bm{E}_{\bm{X}_1}(\lambda_1)\bm{X}_{1,\bm{P}}d\bm{E}_{\bm{X}_2}(\lambda_2) \nonumber \\
&&- \int\limits_{\lambda_1 \in \sigma(\bm{X}_1)}\int\limits_{\lambda_2 \in \sigma(\bm{X}_2)}\frac{f(\lambda_{1}) - f(\lambda_{2})}{\lambda_{1} -\lambda_{2}}d\bm{E}_{\bm{X}_1}(\lambda_1)\bm{X}_{2,\bm{P}}d\bm{E}_{\bm{X}_2}(\lambda_2) \nonumber \\
&=_1&\int\limits_{\lambda_1 \in \sigma(\bm{X}_1)}\int\limits_{\lambda_2 \in \sigma(\bm{X}_2)}\lambda_{1}\frac{f(\lambda_{1}) - f(\lambda_{2})}{\lambda_{1} -\lambda_{2}}d\bm{E}_{\bm{X}_1}(\lambda_1)d\bm{E}_{\bm{X}_2}(\lambda_2) \nonumber \\
&&-\int\limits_{\lambda_1 \in \sigma(\bm{X}_1)}\int\limits_{\lambda_2 \in \sigma(\bm{X}_2)}\lambda_{2}\frac{f(\lambda_{1}) - f(\lambda_{2})}{\lambda_{1} -\lambda_{2}}d\bm{E}_{\bm{X}_1}(\lambda_1)d\bm{E}_{\bm{X}_2}(\lambda_2)  \nonumber \\
&=&\int\limits_{\lambda_1 \in \sigma(\bm{X}_1)}\int\limits_{\lambda_2 \in \sigma(\bm{X}_2)}f(\lambda_{1})d\bm{E}_{\bm{X}_1}(\lambda_1)d\bm{E}_{\bm{X}_2}(\lambda_2) \nonumber \\
&&-\int\limits_{\lambda_1 \in \sigma(\bm{X}_1)}\int\limits_{\lambda_2 \in \sigma(\bm{X}_2)}f(\lambda_{2})d\bm{E}_{\bm{X}_1}(\lambda_1)d\bm{E}_{\bm{X}_2}(\lambda_2)\nonumber \\
&=_2&\underbrace{\left[\int\limits_{\lambda_1 \in \sigma(\bm{X}_1)}f(\lambda_{1})d\bm{E}_{\bm{X}_1}(\lambda_1)- \int\limits_{\lambda_2 \in \sigma(\bm{X}_2)}f(\lambda_{2})d\bm{E}_{\bm{X}_2}(\lambda_2)\right]}_{\mbox{Part I}},
\end{eqnarray}
where we apply $\bm{X}_{1,\bm{P}}$ and $\bm{X}_{2,\bm{P}}$ definitions in $=_1$ with relations provided by Eq.~\eqref{eq3: lma: ind of Proj and Nilp}, and apply relations $\int\limits_{\lambda_1 \in \sigma(\bm{X}_1)}d\bm{E}_{\bm{X}_1}(\lambda_1)=\int\limits_{\lambda_2 \in \sigma(\bm{X}_2)}d\bm{E}_{\bm{X}_2}(\lambda_2) = \bm{I}$ in $=_2$. 

On the other hand, by applying the GDOI definition given by Eq.~\eqref{eq1:  GDOI def}, we also have
\begin{eqnarray}\label{eq6:thm:pert formula diff DOI}
\lefteqn{\left[T_{f(x_1)}^{\bm{X}_{1,\bm{N}},\bm{X}_{2,\bm{P}}}(\bm{I})-T_{f(x_2)}^{\bm{X}_{1,\bm{P}},\bm{X}_{2,\bm{N}}}(\bm{I})\right]}\nonumber \\
&=&\int\limits_{\lambda_1 \in \sigma(\bm{X}_1)}\int\limits_{\lambda_2 \in \sigma(\bm{X}_2)}\sum\limits_{q_1=1}^{m_{\lambda_1}-1}\frac{f^{(q_1)}(\lambda_{1})}{q_1!}(\bm{X}_1-\lambda_1\bm{I})^{q_1}d\bm{E}_{\bm{X}_1}(\lambda_1) \bm{I}d\bm{E}_{\bm{X}_2}(\lambda_2)-\nonumber \\
&&\int\limits_{\lambda_1 \in \sigma(\bm{X}_1)}\int\limits_{\lambda_2 \in \sigma(\bm{X}_2)}\sum\limits_{q_2=1}^{m_{\lambda_2}-1}\frac{f^{(q_2)}(\lambda_{2})}{q_2!}d\bm{E}_{\bm{X}_1}(\lambda_1)\bm{I}(\bm{X}_2-\lambda_2\bm{I})^{q_2}d\bm{E}_{\bm{X}_2}(\lambda_2)\nonumber \\
&=_1&\left[\int\limits_{\lambda_1 \in \sigma(\bm{X}_1)}\sum\limits_{q_1=1}^{m_{\lambda_1}-1}\frac{f^{(q_1)}(\lambda_{1})}{q_1!}(\bm{X}_1-\lambda_1\bm{I})^{q_1}d\bm{E}_{\bm{X}_1}(\lambda_1) \right. \nonumber \\
&& \left.-\int\limits_{\lambda_2 \in \sigma(\bm{X}_2)}\sum\limits_{q_2=1}^{m_{\lambda_2}-1}\frac{f^{(q_2)}(\lambda_{2})}{q_2!}(\bm{X}_2-\lambda_2\bm{I})^{q_2}d\bm{E}_{\bm{X}_2}(\lambda_2)\right]_{\mbox{Part II}},
\end{eqnarray}
where we apply $\int\limits_{\lambda_1 \in \sigma(\bm{X}_1)}d\bm{E}_{\bm{X}_1}(\lambda_1)=\int\limits_{\lambda_2 \in \sigma(\bm{X}_2)}d\bm{E}_{\bm{X}_2}(\lambda_2) = \bm{I}$  in $=_1$.

This theorem is proved by combining Eq.~\eqref{eq4:thm:pert formula diff DOI}, Eq.~\eqref{eq5:thm:pert formula diff DOI} and Eq.~\eqref{eq6:thm:pert formula diff DOI}.
$\hfill\Box$

Comparing Theorem~\ref{thm:pert formula diff DOI} with the conventional perturbation formula under Hermitian operators assumptions of $\bm{X}_1$ and $\bm{X}_2$, we have one extra term, which is $\left[T_{f(x_1)}^{\bm{X}_{1,\bm{N}},\bm{X}_{2,\bm{P}}}(\bm{I})-T_{f(x_2)}^{\bm{X}_{1,\bm{P}},\bm{X}_{2,\bm{N}}}(\bm{I})\right]$. In the remaining part of this subsection, we will explore properties of the following extra term for GDOI with CSOs: 
\begin{eqnarray}\label{eq:extra nilpotent term in pertf}
\mu(\bm{X}_1, \bm{X}_2, f)&\define&\left[T_{f(x_1)}^{\bm{X}_{1,\bm{N}},\bm{X}_{2,\bm{P}}}(\bm{I})-T_{f(x_2)}^{\bm{X}_{1,\bm{P}},\bm{X}_{2,\bm{N}}}(\bm{I})\right].
\end{eqnarray}

Based on the proof of  Lemma 4 in~\cite{chang2025GDOIFinite}, we have the following Proposition~\ref{prop: GDOI div DOI} to quantify the deviation between the GDOI and the DOI due to the non-self-adjoint nature of operators $\bm{X}_1$ and $\bm{X}_2$.
\begin{proposition}\label{prop: GDOI div DOI}
Let $\mu(\bm{X}_1, \bm{X}_2, f)$ be an operator obtained by Eq.~\eqref{eq:extra nilpotent term in pertf}. The operator $\mu(\bm{X}_1, \bm{X}_2, f)$ can be categorized as either non-nilpotent or nilpotent. Define $\ell_1$ as the measure of non-zero eigenvalues of $\mu(\bm{X}_1, \bm{X}_2, f)$. The following conditions hold:

1. If $\ell_1 > 0$, the operator $\mu(\bm{X}_1, \bm{X}_2, f)$ is non-nilpotent.\\

2. If $\ell_1 = 0$, the operator $\mu(\bm{X}_1, \bm{X}_2, f)$ is nilpotent. \\

For the nilpotent case, let $\ell_2$ represent the nilpotent degree of $\mu(\bm{X}_1, \bm{X}_2, f)$, meaning $\ell_2$ is the smallest integer $k$ such that
\[
\mu^{k}(\bm{X}_1, \bm{X}_2, f) = \bm{0}.
\]

For two operators $\mu(\bm{X}_1, \bm{X}_2, f)$ and $\mu(\bm{X}'_1, \bm{X}'_2, f')$, if $\ell_1 = \ell'_1$ and $\ell_2 = \ell'_2$, the operator norm of the operators $\mu(\bm{X}_1, \bm{X}_2, f)$ and $\mu(\bm{X}'_1, \bm{X}'_2, f')$ can be used to distinguish between them by the total-ordering structure via triple $(\ell_1, \ell_2, r)$. Moreover, if $(\ell_1, \ell_2, r) = (0,0,0)$, the GDOI reduces to the DOI, indicating that both operators $\bm{X}_1$ and $\bm{X}_2$ are self-adjoint.
\end{proposition}

We pose two design questions concerning the conditions under which the operator $\mu(\bm{X}_1, \bm{X}_2, f)$ becomes \emph{nilpotent}:

\begin{equation}
\begin{aligned}
\text{(1)} \quad &\text{Given the operators} \bm{X}_1 \text{ and } \bm{X}_2, \text{ what analytic function } f \\
&\text{makes the operator } \mu(\bm{X}_1, \bm{X}_2, f) \text{ nilpotent?}
\end{aligned}
\label{eq:design_Q1}
\end{equation}

\begin{equation}
\begin{aligned}
\text{(2)} \quad &\text{Given the analytic function } f, \text{ what operators} \bm{X}_1 \text{ and } \bm{X}_2 \\
&\text{make the operator } \mu(\bm{X}_1, \bm{X}_2, f) \text{ nilpotent?}
\end{aligned}
\label{eq:design_Q2}
\end{equation}

The next Proposition~\ref{prop: commutative make mu nilpotent} shows that the operator $\mu(\bm{X}_1, \bm{X}_2, f)$ is a nilpotent operator if nilpotent components in $\bm{X}_1$ and $\bm{X}_2$ commute each other. 

\begin{proposition}\label{prop: commutative make mu nilpotent}
As the setting in Theorem~\ref{thm:pert formula}, if $(\bm{X}_1-\lambda_1\bm{I})d\bm{E}_{\bm{X}_1}(\lambda_1)$ and $(\bm{X}_2-\lambda_2\bm{I})d\bm{E}_{\bm{X}_2}(\lambda_2)$ commute each other, the operator $\mu(\bm{X}_1, \bm{X}_2, f)$ is nilpotent. The nilpotent index of the operator $\mu(\bm{X}_1, \bm{X}_2, f)$ is at most $\max\limits_{\lambda_1 \in \sigma(\bm{X}_1)}[m_{\lambda_1}]+\max\limits_{\lambda_2 \in \sigma(\bm{X}_2)}[m_{\lambda_2}]$.
\end{proposition}
\textbf{Proof:}
We first show that $\bm{A}_1$ and $\bm{A}_2$ are commutative nilpotent operators, then, for any scalars \(a, b \in \mathbb{C}\), the linear combination \(\bm{A} = a\bm{A}_1 + b\bm{A}_2\) is also nilpotent.

Suppose \(\bm{A}_1^{k_1} = 0\) and \(\bm{A}_2^{k_2} = 0\) where \(k_1\) and \(k_2\) are the nilpotent indices of \(\bm{A}_1\) and \(\bm{A}_2\), respectively. Consider the linear combination \(\bm{A} = a\bm{A}_1 + b\bm{A}_2\). We analyze \(\bm{A}^{k_1 + k_2}\) using the binomial theorem:

\[
\bm{A}^{k_1 + k_2} = (a\bm{A}_1 + b\bm{A}_2)^{k_1 + k_2} = \sum_{j=0}^{k_1 + k_2} \binom{k_1 + k_2}{j} a^j b^{k_1 + k_2 - j} \bm{A}_1^j \bm{A}_2^{k_1 + k_2 - j}
\]

Note that if \(j \geq k_1\), then \(\bm{A}_1^j = \bm{0}\). Also, if \(k_1 + k_2 - j \geq k_2\), then \(\bm{A}_2^{k_1 + k_2 - j} = 0\). Since at least one of these conditions is always true for every term in the sum, all terms in the expansion are zero. Therefore:

\[
\bm{A}^{k_1 + k_2} = \bm{0}
\]

Since there exists a finite power of \(\bm{A}\) that is zero, \(\bm{A}\) is nilpotent. The nilpotent index of \(\bm{A}\) is at most \(k_1 + k_2\). Because we also have 
\begin{eqnarray}
\mu(\bm{X}_1, \bm{X}_2, f)&=&\int\limits_{\lambda_1 \in \sigma(\bm{X}_1)}\sum\limits_{q_1=1}^{m_{\lambda_1}-1}\frac{f^{(q_1)}(\lambda_{1})}{q_1!}(\bm{X}_1-\lambda_1\bm{I})^{q_1}d\bm{E}_{\bm{X}_1}(\lambda_1)\nonumber \\
&&-\int\limits_{\lambda_2 \in \sigma(\bm{X}_2)}\sum\limits_{q_2=1}^{m_{\lambda_2}-1}\frac{f^{(q_2)}(\lambda_{2})}{q_2!}(\bm{X}_2-\lambda_2\bm{I})^{q_2}d\bm{E}_{\bm{X}_2}(\lambda_2),
\end{eqnarray}
then, the nilpotent index for the first term in the above equation is $\max\limits_{\lambda_1  \in \sigma(\bm{X}_1)}[m_{\lambda_1}]$ and the nilpotent index for the second term in the above equation is $\max\limits_{\lambda_2  \in \sigma(\bm{X}_2)}[m_{\lambda_2}]$. From the above fact just proved, the nilpotent index of the operator $\mu(\bm{X}_1, \bm{X}_2, f)$ is at most $\max\limits_{\lambda_1 \in \sigma(\bm{X}_1)}[m_{\lambda_1}]+\max\limits_{\lambda_2 \in \sigma(\bm{X}_2)}[m_{\lambda_2}]$.
$\hfill\Box$

\section{Inequalities Related to GDOI for Continuous Spectrum Operators}\label{sec:Inequalities Related to GDOI for Continuous Spectrum Operators}

Norm inequalities associated with GDOI for continuous spectrum operators are presented in Section~\ref{sec:Norm Inequalities}. Lipschitz estimates for continuous spectrum operators are derived in Section~\ref{sec:Lipschitz Estimations}.

\subsection{Norm Inequalities}\label{sec:Norm Inequalities}

We need the following Lemma~\ref{lma:conv triangle for Frob norm} to provide the lower bound for the GDOI norm. Note that Lemma~\ref{lma:conv triangle for Frob norm} is Lemma 3 in our previous work~\cite{chang2025GDOIFinite} and the proof is similar to there by replacing the matrix norm to the operator norm here. 
\begin{lemma}\label{lma:conv triangle for Frob norm}
Given $n$ operators $\bm{A}_1, \bm{A}_2, \ldots, \bm{A}_n$ such that $\left\Vert\bm{A}_1\right\Vert\geq\left\Vert\bm{A}_2\right\Vert\geq\ldots\geq\left\Vert\bm{A}_n\right\Vert$, then, we have
\begin{eqnarray}\label{eq1:lma:conv triangle for Frob norm}
\left\Vert\sum\limits_{i=1}^n \bm{A}_i \right\Vert \geq \max\left[0, \left\Vert\bm{A}_1\right\Vert - \sum\limits_{i=2}^n \left\Vert\bm{A}_i\right\Vert\right]
\end{eqnarray}
\end{lemma}

Theorem~\ref{thm: GDOI norm} below is provided to give the upper bound and the lower bound for the norm  associated with GDOI for continuous spectrum operators. 
\begin{theorem}\label{thm: GDOI norm}
We have the upper bound for the operator norm of $T_{\beta}^{\bm{X}_1,\bm{X}_2}(\bm{Y})$, which is given by
\begin{eqnarray}\label{eq1: thm: GDOI norm}
\lefteqn{\left\Vert T_{\beta}^{\bm{X}_1,\bm{X}_2}(\bm{Y}) \right\Vert\leq}\nonumber \\
&&\left[\max\limits_{\lambda_1 \in \sigma_{\bm{X}_1},\lambda_2 \in \sigma_{\bm{X}_2}} \left\vert\beta(\lambda_1,\lambda_2)\right\vert\right]\left\Vert\bm{Y}\right\Vert \nonumber \\
&&+\left\{\int\limits_{\lambda_2 \in \sigma(\bm{X}_2)}\sum_{q_2=1}^{m_{\lambda_2}-1}\left[\max\limits_{\lambda_1 \in \sigma_{\bm{X}_1},\lambda_2 \in \sigma_{\bm{X}_2}}\left\vert\frac{\beta^{(-,q_2)}(\lambda_{1},\lambda_{2})}{q_2!}\right\vert\right]\left\Vert(\bm{X}_2-\lambda_2\bm{I})^{q_2}\right\Vert d\bm{E}_{\bm{X}_2}(\lambda_2)\right\}\left\Vert\bm{Y}\right\Vert\nonumber \\
&&+\left\{\int\limits_{\lambda_1 \in \sigma(\bm{X}_1)}\sum_{q_1=1}^{m_{\lambda_1}-1}\left[\max\limits_{\lambda_1 \in \sigma_{\bm{X}_1},\lambda_2 \in \sigma_{\bm{X}_2}}\left\vert\frac{\beta^{(q_1,-)}(\lambda_{1},\lambda_{2})}{q_1!}\right\vert\right]\left\Vert(\bm{X}_1-\lambda_1\bm{I})^{q_1} \right\Vert d\bm{E}_{\bm{X}_1}(\lambda_1)\right\}\left\Vert\bm{Y}\right\Vert  \nonumber \\
&&+ \Bigg\{ \int\limits_{\lambda_1 \in \sigma(\bm{X}_1)}\int\limits_{\lambda_2 \in \sigma(\bm{X}_2)}\sum_{q_1=1}^{m_{\lambda_1}-1}\sum_{q_2=1}^{m_{\lambda_2}-1}\left[\max\limits_{\lambda_1 \in \sigma_{\bm{X}_1},\lambda_2 \in \sigma_{\bm{X}_2}}\left\vert\frac{\beta^{(q_1,q_2)}(\lambda_{1},\lambda_{2})}{q_1!q_2!}\right\vert\right]  \nonumber\\
&&~~~ \times  \left\Vert(\bm{X}_1-\lambda_1\bm{I})^{q_1} \right\Vert \left\Vert(\bm{X}_2-\lambda_2\bm{I})^{q_2} \right\Vert d\bm{E}_{\bm{X}_1}(\lambda_1) d\bm{E}_{\bm{X}_2}(\lambda_2)\Bigg\}\left\Vert\bm{Y}\right\Vert.
\end{eqnarray}
where $\sigma_{\bm{X}_1}$ and $\sigma_{\bm{X}_2}$ are spectrums of the operator $\bm{X}_1$ and the operator $\bm{X}_2$, respectively.

On the other hand, let us define the following operators
\begin{eqnarray}\label{eq1-1: thm: GDOI norm}
\bm{A}_1&\define&\int\limits_{\lambda_1 \in \sigma(\bm{X}_1)}\int\limits_{\lambda_2 \in \sigma(\bm{X}_2)}\beta(\lambda_{1}, \lambda_{2}) d\bm{E}_{\bm{X}_1}(\lambda_1)\bm{Y}d\bm{E}_{\bm{X}_2}(\lambda_2)\nonumber \\
\bm{A}_2&\define&\int\limits_{\lambda_1 \in \sigma(\bm{X}_1)}\int\limits_{\lambda_2 \in \sigma(\bm{X}_2)}\sum_{q_2=1}^{m_{\lambda_2}-1}\frac{\beta^{(-,q_2)}(\lambda_{1},\lambda_{2})}{q_2!}d\bm{E}_{\bm{X}_1}(\lambda_1)\bm{Y}(\bm{X}_2-\lambda_2\bm{I})^{q_2}d\bm{E}_{\bm{X}_2}(\lambda_2)\nonumber \\
\bm{A}_3&\define&\int\limits_{\lambda_1 \in \sigma(\bm{X}_1)}\int\limits_{\lambda_2 \in \sigma(\bm{X}_2)}\sum_{q_1=1}^{m_{\lambda_1}-1}\frac{\beta^{(q_1,-)}(\lambda_{1},\lambda_{2})}{q_1!}(\bm{X}_1-\lambda_1\bm{I})^{q_1}d\bm{E}_{\bm{X}_1}(\lambda_1)\bm{Y}d\bm{E}_{\bm{X}_2}(\lambda_2)\nonumber \\
\bm{A}_4&\define&\int\limits_{\lambda_1 \in \sigma(\bm{X}_1)}\int\limits_{\lambda_2 \in \sigma(\bm{X}_2)}\sum_{q_1=1}^{m_{\lambda_1}-1}\sum_{q_2=1}^{m_{\lambda_2}-1}\frac{\beta^{(q_1,q_2)}(\lambda_{1},\lambda_{2})}{q_1!q_2!}\nonumber \\
&&\times (\bm{X}_1-\lambda_1\bm{I})^{q_1}d\bm{E}_{\bm{X}_1}(\lambda_1)\bm{Y}(\bm{X}_2-\lambda_2\bm{I})^{q_2}d\bm{E}_{\bm{X}_2}(\lambda_2).
\end{eqnarray}
then, we have the lower bound for the operator norm of $T_{\beta}^{\bm{X}_1,\bm{X}_2}(\bm{Y})$, which is given by
\begin{eqnarray}\label{eq2: thm: GDOI norm}
\left\Vert T_{\beta}^{\bm{X}_1,\bm{X}_2}(\bm{Y}) \right\Vert\geq \max\left[0, \left\Vert\bm{A}_{\sigma(1)}\right\Vert - \sum\limits_{i=2}^4 \left\Vert\bm{A}_{\sigma(i)}\right\Vert\right],
\end{eqnarray}
where $\sigma$ is the permutation of operators $\bm{A}_i$ for $i=1,2,3,4$ such that $\left\Vert\bm{A}_{\sigma(1)}\right\Vert \geq \left\Vert\bm{A}_{\sigma(2)}\right\Vert \geq\left\Vert\bm{A}_{\sigma(3)}\right\Vert \geq\left\Vert\bm{A}_{\sigma(4)}\right\Vert \geq$. 

Further, if we have $\left[\min\limits_{\lambda_1 \in \sigma_{\bm{X}_1},\lambda_2 \in \sigma_{\bm{X}_2}} \beta(\lambda_1,\lambda_2)\right]\left\Vert\bm{Y}\right\Vert\geq \left\Vert\bm{A}_{2}\right\Vert+\left\Vert\bm{A}_{3}\right\Vert + \left\Vert\bm{A}_{4}\right\Vert$,  the lower bound for the operator norm of $T_{\beta}^{\bm{X}_1,\bm{X}_2}(\bm{Y})$ can be expressed by
\begin{eqnarray}\label{eq3: thm: GDOI norm}
\left\Vert T_{\beta}^{\bm{X}_1,\bm{X}_2}(\bm{Y}) \right\Vert\geq \left[\min\limits_{\lambda_1 \in \sigma_{\bm{X}_1},\lambda_2 \in \sigma_{\bm{X}_2}} \beta(\lambda_1,\lambda_2)\right]\left\Vert\bm{Y}\right\Vert - \sum\limits_{i=2}^4 \left\Vert\bm{A}_{i}\right\Vert.
\end{eqnarray}
\end{theorem}
\textbf{Proof:}
From the definition of $T_{\beta}^{\bm{X}_1,\bm{X}_2}(\bm{Y})$ and the triangle inequality, we have 
\begin{eqnarray}\label{eq4: thm: GDOI norm}
\left\Vert T_{\beta}^{\bm{X}_1,\bm{X}_2}(\bm{Y}) \right\Vert&\leq&\underbrace{\left\Vert \int\limits_{\lambda_1 \in \sigma(\bm{X}_1)}\int\limits_{\lambda_2 \in \sigma(\bm{X}_2)}\beta(\lambda_{1}, \lambda_{2})d\bm{E}_{\bm{X}_1}(\lambda_1)\bm{Y}d\bm{E}_{\bm{X}_2}(\lambda_2)\right\Vert}_{\mbox{Part I}}\nonumber \\
&&+\underbrace{\left\Vert \int\limits_{\lambda_1 \in \sigma(\bm{X}_1)}\int\limits_{\lambda_2 \in \sigma(\bm{X}_2)}\sum_{q_2=1}^{m_{\lambda_2}-1}\frac{\beta^{(-,q_2)}(\lambda_{1},\lambda_{2})}{q_2!}d\bm{E}_{\bm{X}_1}(\lambda_1)\bm{Y}(\bm{X}_2-\lambda_2\bm{I})^{q_2}d\bm{E}_{\bm{X}_2}(\lambda_2) \right\Vert}_{\mbox{Part II}}\nonumber \\
&&+\underbrace{\left\Vert\int\limits_{\lambda_1 \in \sigma(\bm{X}_1)}\int\limits_{\lambda_2 \in \sigma(\bm{X}_2)}\sum_{q_1=1}^{m_{\lambda_1}-1}\frac{\beta^{(q_1,-)}(\lambda_{1},\lambda_{2})}{q_1!}(\bm{X}_1-\lambda_1\bm{I})^{q_1}d\bm{E}_{\bm{X}_1}(\lambda_1)\bm{Y}d\bm{E}_{\bm{X}_2}(\lambda_2) \right\Vert}_{\mbox{Part III}}\nonumber \\
&&+\left\Vert\int\limits_{\lambda_1 \in \sigma(\bm{X}_1)}\int\limits_{\lambda_2 \in \sigma(\bm{X}_2)}\sum_{q_1=1}^{m_{\lambda_1}-1}\sum_{q_2=1}^{m_{\lambda_2}-1}\frac{\beta^{(q_1,q_2)}(\lambda_{1},\lambda_{2})}{q_1!q_2!} \right. \nonumber \\
&&~~~~ \left. \times (\bm{X}_1-\lambda_1\bm{I})^{q_1}d\bm{E}_{\bm{X}_1}(\lambda_1)\bm{Y}(\bm{X}_2-\lambda_2\bm{I})^{q_2}d\bm{E}_{\bm{X}_2}(\lambda_2)\right\Vert_{\mbox{Part IV}}.
\end{eqnarray}

For Part I, we have 
\begin{eqnarray}\label{eq4-1: thm: GDOI norm}
\lefteqn{\underbrace{\left\Vert\int\limits_{\lambda_1 \in \sigma(\bm{X}_1)}\int\limits_{\lambda_2 \in \sigma(\bm{X}_2)}\beta(\lambda_{1}, \lambda_{2})d\bm{E}_{\bm{X}_1}(\lambda_1)\bm{Y}d\bm{E}_{\bm{X}_2}(\lambda_2)\right\Vert}_{\mbox{Part I}}}\nonumber \\
&\leq&\left[\max\limits_{\lambda_1 \in \sigma_{\bm{X}_1},\lambda_2 \in \sigma_{\bm{X}_2}} \left\vert\beta(\lambda_1,\lambda_2)\right\vert\right]\left\Vert\int\limits_{\lambda_1 \in \sigma(\bm{X}_1)}\int\limits_{\lambda_2 \in \sigma(\bm{X}_2)}d\bm{E}_{\bm{X}_1}(\lambda_1)\bm{Y}d\bm{E}_{\bm{X}_2}(\lambda_2)\right\Vert\nonumber \\
&=&\left[\max\limits_{\lambda_1 \in \sigma_{\bm{X}_1},\lambda_2 \in \sigma_{\bm{X}_2}} \left\vert\beta(\lambda_1,\lambda_2)\right\vert\right]\left\Vert\bm{Y}\right\Vert,
\end{eqnarray}
where we use the fact that $\int\limits_{\lambda_1 \in \sigma(\bm{X}_1)}d\bm{E}_{\bm{X}_1}(\lambda_1)= \bm{I}$ and $\int\limits_{\lambda_2\in \sigma(\bm{X}_2)}d\bm{E}_{\bm{X}_2}(\lambda_2)= \bm{I}$ at the last equality. 

For Part II, we have 
\begin{eqnarray}\label{eq4-2: thm: GDOI norm}
\lefteqn{\underbrace{\left\Vert\int\limits_{\lambda_1 \in \sigma(\bm{X}_1)}\int\limits_{\lambda_2 \in \sigma(\bm{X}_2)}\sum_{q_2=1}^{m_{\lambda_2}-1}\frac{\beta^{(-,q_2)}(\lambda_{1},\lambda_{2})}{q_2!}d\bm{E}_{\bm{X}_1}(\lambda_1) \bm{Y} (\bm{X}_2-\lambda_2\bm{I})^{q_2}d\bm{E}_{\bm{X}_2}(\lambda_2)\right\Vert}_{\mbox{Part II}}}\nonumber \\
&\leq&\int\limits_{\lambda_2 \in \sigma(\bm{X}_2)}\sum_{q_2=1}^{m_{\lambda_2}-1}\left[\max\limits_{\lambda_1 \in \sigma_{\bm{X}_1},\lambda_2 \in \sigma_{\bm{X}_2}} \left\vert\frac{\beta^{(-,q_2)}(\lambda_{1},\lambda_{2})}{q_2!}\right\vert\right]\left\Vert\int\limits_{\lambda_1 \in \sigma(\bm{X}_1)} d\bm{E}_{\bm{X}_1}(\lambda_1)\bm{Y} (\bm{X}_2-\lambda_2\bm{I})^{q_2}d\bm{E}_{\bm{X}_2}(\lambda_2)\right\Vert\nonumber \\
&\leq_1&\left\{\int\limits_{\lambda_2 \in \sigma(\bm{X}_2)}\sum_{q_2=1}^{m_{\lambda_2}-1}\left[\max\limits_{\lambda_1 \in \sigma_{\bm{X}_1},\lambda_2 \in \sigma_{\bm{X}_2}}\left\vert\frac{\beta^{(-,q_2)}(\lambda_{1},\lambda_{2})}{q_2!}\right\vert\right]\left\Vert (\bm{X}_2-\lambda_2\bm{I})^{q_2}\right\Vert d\bm{E}_{\bm{X}_2}(\lambda_2) \right\}\left\Vert\bm{Y}\right\Vert,
\end{eqnarray}
where we use the fact that $\int\limits_{\lambda_1 \in \sigma(\bm{X}_1)}d\bm{E}_{\bm{X}_1}(\lambda_1)=\bm{I}$ and the norm multiplicative inequality at the inequality $\leq_1$.

For Part III, we have 
\begin{eqnarray}\label{eq4-3: thm: GDOI norm}
\lefteqn{\underbrace{\left\Vert\int\limits_{\lambda_1 \in \sigma(\bm{X}_1)}\int\limits_{\lambda_2 \in \sigma(\bm{X}_2)}\sum_{q_1=1}^{m_{\lambda_1}-1}\frac{\beta^{(q_1,-)}(\lambda_{1},\lambda_{2})}{q_1!}(\bm{X}_1-\lambda_1\bm{I})^{q_1}d\bm{E}_{\bm{X}_1}(\lambda_1)\bm{Y}d\bm{E}_{\bm{X}_2}(\lambda_2)\right\Vert}_{\mbox{Part III}}}\nonumber \\
&\leq&\int\limits_{\lambda_1 \in \sigma(\bm{X}_1)}\sum_{q_1=1}^{m_{\lambda_1}-1}\left[\max\limits_{\lambda_1 \in \sigma_{\bm{X}_1},\lambda_2 \in \sigma_{\bm{X}_2}} \left\vert\frac{\beta^{(q_1,-)}(\lambda_{1},\lambda_{2})}{q_1!}\right\vert\right]\left\Vert\int\limits_{\lambda_2 \in \sigma(\bm{X}_2)}(\bm{X}_1-\lambda_1\bm{I})^{q_1}d\bm{E}_{\bm{X}_1}(\lambda_1)\bm{Y}d\bm{E}_{\bm{X}_2}(\lambda_2)\right\Vert\nonumber \\
&\leq_1&\left\{\int\limits_{\lambda_1 \in \sigma(\bm{X}_1)}\sum_{q_1=1}^{m_{\lambda_1}-1}\left[\max\limits_{\lambda_1 \in \sigma_{\bm{X}_1},\lambda_2 \in \sigma_{\bm{X}_2}}\left\vert\frac{\beta^{(q_1,-)}(\lambda_{1},\lambda_{2})}{q_1!}\right\vert\right]\left\Vert(\bm{X}_1-\lambda_1\bm{I})^{q_1}\right\Vert d\bm{E}_{\bm{X}_1}(\lambda_1)\right\}\left\Vert\bm{Y}\right\Vert,
\end{eqnarray}
where we use the fact that $\int\limits_{\lambda_2 \in \sigma(\bm{X}_2)}d\bm{E}_{\bm{X}_2}(\lambda_2)=\bm{I}$ and the norm multiplicative inequality at the  inequality $\leq_1$.

For Part IV, we have 
\begin{eqnarray}\label{eq4-4: thm: GDOI norm}
\lefteqn{\underbrace{\left\Vert\int\limits_{\lambda_1 \in \sigma(\bm{X}_1)}\int\limits_{\lambda_2 \in \sigma(\bm{X}_2)}\sum_{q_1=1}^{m_{\lambda_1}-1}\sum_{q_2=1}^{m_{\lambda_2}-1}\frac{\beta^{(q_1,q_2)}(\lambda_{1},\lambda_{2})}{q_1!q_2!}(\bm{X}_1-\lambda_1\bm{I})^{q_1}d\bm{E}_{\bm{X}_1}(\lambda_1)\bm{Y} (\bm{X}_2-\lambda_2\bm{I})^{q_2}d\bm{E}_{\bm{X}_2}(\lambda_2) \right\Vert}_{\mbox{Part IV}}}\nonumber \\
&\leq& \int\limits_{\lambda_1 \in \sigma(\bm{X}_1)}\int\limits_{\lambda_2 \in \sigma(\bm{X}_2)}\sum_{q_1=1}^{m_{\lambda_1}-1}\sum_{q_2=1}^{m_{\lambda_2}-1}\left[\max\limits_{\lambda_1 \in \sigma_{\bm{X}_1},\lambda_2 \in \sigma_{\bm{X}_2}}\left\vert\frac{\beta^{(q_1,q_2)}(\lambda_{1},\lambda_{2})}{q_1!q_2!}\right\vert\right]\nonumber \\
&& \times \left\Vert(\bm{X}_1-\lambda_1\bm{I})^{q_1}\bm{Y} (\bm{X}_2-\lambda_2\bm{I})^{q_2}\right\Vert d\bm{E}_{\bm{X}_1}(\lambda_1)d\bm{E}_{\bm{X}_2}(\lambda_2)\nonumber \\
&\leq_1& \Bigg\{\int\limits_{\lambda_1 \in \sigma(\bm{X}_1)}\int\limits_{\lambda_2 \in \sigma(\bm{X}_2)}\sum_{q_1=1}^{m_{\lambda_1}-1}\sum_{q_2=1}^{m_{\lambda_2}-1}\left[\max\limits_{\lambda_1 \in \sigma_{\bm{X}_1},\lambda_2 \in \sigma_{\bm{X}_2}}\left\vert\frac{\beta^{(q_1,q_2)}(\lambda_{1},\lambda_{2})}{q_1!q_2!}\right\vert\right]\nonumber\\
&&~~~\times \left\Vert(\bm{X}_1-\lambda_1\bm{I})^{q_1} \right\Vert \left\Vert (\bm{X}_2-\lambda_2\bm{I})^{q_2}\right\Vert d\bm{E}_{\bm{X}_1}(\lambda_1)d\bm{E}_{\bm{X}_2}(\lambda_2)\Bigg\}\left\Vert\bm{Y}\right\Vert.
\end{eqnarray}
where we use the norm multiplicative inequality again at the inequality $\leq_1$. By combining previous four parts from Eq.~\eqref{eq4-1: thm: GDOI norm} to Eq.~\eqref{eq4-4: thm: GDOI norm}, we have the desired upper bound given by Eq.~\eqref{eq1: thm: GDOI norm}

For the lower bound of $\left\Vert T_{\beta}^{\bm{X}_1,\bm{X}_2}(\bm{Y})\right\Vert$, we have Eq.~\eqref{eq2: thm: GDOI norm} immediatedly from Lemma~\ref{lma:conv triangle for Frob norm}.

If we have $\left[\min\limits_{\lambda_1 \in \sigma_{\bm{X}_1},\lambda_2 \in \sigma_{\bm{X}_2}} \beta(\lambda_1,\lambda_2)\right]\left\Vert\bm{Y}\right\Vert\geq \left\Vert\bm{A}_{2}\right\Vert+\left\Vert\bm{A}_{3}\right\Vert + \left\Vert\bm{A}_{4}\right\Vert$ and Lemma~\ref{lma:conv triangle for Frob norm}, we have
\begin{eqnarray}\label{eq5: thm: GDOI norm}
\left\Vert T_{\beta}^{\bm{X}_1,\bm{X}_2}(\bm{Y})\right\Vert&\geq&\left\Vert\bm{A}_1\right\Vert-(\left\Vert\bm{A}_{2}\right\Vert+\left\Vert\bm{A}_{3}\right\Vert + \left\Vert\bm{A}_{4}\right\Vert) \nonumber \\
&\geq&\left[\min\limits_{\lambda_1 \in \sigma_{\bm{X}_1},\lambda_2 \in \sigma_{\bm{X}_2}} \beta(\lambda_1,\lambda_2)\right]\left\Vert\int\limits_{\lambda_1 \in \sigma(\bm{X}_1)}\int\limits_{\lambda_2 \in \sigma(\bm{X}_2)}d\bm{E}_{\bm{X}_1}(\lambda_1)\bm{Y}d\bm{E}_{\bm{X}_2}(\lambda_2) \right\Vert\nonumber \\
& &-  (\left\Vert\bm{A}_{2}\right\Vert+\left\Vert\bm{A}_{3}\right\Vert + \left\Vert\bm{A}_{4}\right\Vert) \nonumber \\
&=&\left[\min\limits_{\lambda_1 \in \sigma_{\bm{X}_1},\lambda_2 \in \sigma_{\bm{X}_2}} \beta(\lambda_1,\lambda_2)\right]\left\Vert\bm{Y}\right\Vert - (\left\Vert\bm{A}_{2}\right\Vert+\left\Vert\bm{A}_{3}\right\Vert + \left\Vert\bm{A}_{4}\right\Vert),
\end{eqnarray}
which is the lower bound of $\left\Vert T_{\beta}^{\bm{X}_1,\bm{X}_2}(\bm{Y})\right\Vert$ given by Eq.~\eqref{eq3: thm: GDOI norm}.  
$\hfill\Box$

\subsection{Lipschitz Estimations}\label{sec:Lipschitz Estimations}

In this section, Theorem~\ref{thm:Lipschitz Estimations} is given to provide the lower and the upper bounds for Lipschitz Estimations.
\begin{theorem}\label{thm:Lipschitz Estimations}
Given an analytic function $f(z)$ within the domain for $|z| < R$, the first operator $\bm{X}_1$ is expressed by
\begin{eqnarray}\label{eq1-1:thm:Lipschitz Estimations}
\bm{X}_1&=&\int\limits_{\lambda_1 \in \sigma(\bm{X}_1)}\lambda_1 d\bm{E}_{\bm{X}_1}(\lambda_1)+
\int\limits_{\lambda_1 \in \sigma(\bm{X}_1)}\left(\bm{X}_1-\lambda_1\bm{I}\right)d\bm{E}_{\bm{X}_1}(\lambda_1),
\end{eqnarray}
where $\left\vert\lambda_{1}\right\vert<R$, and second operator $\bm{X}_2$ is expressed by 
\begin{eqnarray}\label{eq1-2:thm:Lipschitz Estimations}
\bm{X}_2&=&\int\limits_{\lambda_2 \in \sigma(\bm{X}_2)}\lambda_2 d\bm{E}_{\bm{X}_2}(\lambda_2)+
\int\limits_{\lambda_2 \in \sigma(\bm{X}_2)}\left(\bm{X}_2-\lambda_2\bm{I}\right)d\bm{E}_{\bm{X}_2}(\lambda_2),
\end{eqnarray}
where $\left\vert\lambda_{2}\right\vert<R$. We also assume that $\lambda_{1} \neq \lambda_{2}$ for any $\lambda_1 \in \sigma(\bm{X}_1)$ and  $\lambda_2 \in \sigma(\bm{X}_2)$. We define $f^{[1]}(\lambda_1, \lambda_2) \define \frac{f(\lambda_1) - f(\lambda_2)}{\lambda_1 - \lambda_2}$. 

Then, we have the following upper bound for Lipschitz estimation:
\begin{eqnarray}\label{eq2:thm:Lipschitz Estimations}
\lefteqn{\left\Vert f(\bm{X}_1) - f(\bm{X}_2) \right\Vert \leq}\nonumber \\
&& \left[\max\limits_{\lambda_1 \in \sigma(\bm{X}_1),\lambda_2 \in \sigma(\bm{X}_2)} \left\vert f^{[1]}(\lambda_1,\lambda_2)\right\vert\right]\left\Vert\bm{X}_1-\bm{X}_2\right\Vert \nonumber \\
&&+\left\{\int\limits_{\lambda_2 \in \sigma(\bm{X}_2)}\sum_{q_2=1}^{m_{\lambda_2}-1}\left[\max\limits_{\lambda_1 \in \sigma(\bm{X}_1),\lambda_2 \in \sigma(\bm{X}_2)}\left\vert\frac{(f^{[1]})^{(-,q_2)}(\lambda_{1},\lambda_{2})}{q_2!}\right\vert\right]\left\Vert(\bm{X}_2-\lambda_2\bm{I})^{q_2}\right\Vert d\bm{E}_{\bm{X}_2}(\lambda_2)\right\}\left\Vert\bm{X}_1-\bm{X}_2\right\Vert\nonumber \\
&&+\left\{\int\limits_{\lambda_1 \in \sigma(\bm{X}_1)}\sum_{q_1=1}^{m_{\lambda_1}-1}\left[\max\limits_{\lambda_1 \in \sigma(\bm{X}_1),\lambda_2 \in \sigma(\bm{X}_2)}\left\vert\frac{(f^{[1]})^{(q_1,-)}(\lambda_{1},\lambda_{2})}{q_1!}\right\vert\right]\left\Vert(\bm{X}_1-\lambda_1\bm{I})^{q_1}\right\Vert d\bm{E}_{\bm{X}_1}(\lambda_1)\right\} \left\Vert\bm{X}_1-\bm{X}_2\right\Vert  \nonumber \\
&&+ \left\{\int\limits_{\lambda_1 \in \sigma(\bm{X}_1)}\int\limits_{\lambda_2 \in \sigma(\bm{X}_2)}\sum_{q_1=1}^{m_{\lambda_1}-1}\sum_{q_2=1}^{m_{\lambda_2}-1}\left[\max\limits_{\lambda_1 \in \sigma(\bm{X}_1),\lambda_2 \in \sigma(\bm{X}_2)}\left\vert\frac{(f^{[1]})^{(q_1,q_2)}(\lambda_{1},\lambda_{2})}{q_1!q_2!}\right\vert\right] \right.\nonumber\\
&&~~~\left. \times\left\Vert(\bm{X}_1-\lambda_1\bm{I})^{q_1} \right\Vert \left\Vert (\bm{X}_2-\lambda_2\bm{I})^{q_2}\right\Vert d\bm{E}_{\bm{X}_1}(\lambda_1) d\bm{E}_{\bm{X}_2}(\lambda_2)\right\}\left\Vert\bm{X}_1-\bm{X}_2\right\Vert.
\end{eqnarray}

On the other hand, let us define the following operators
\begin{eqnarray}\label{eq3:thm:Lipschitz Estimations}
\bm{A}_1&\define&\int\limits_{\lambda_1 \in \sigma(\bm{X}_1)}\int\limits_{\lambda_2 \in \sigma(\bm{X}_2)}f^{[1]}(\lambda_{1}, \lambda_{2})d\bm{E}_{\bm{X}_1}(\lambda_1)(\bm{X}_1-\bm{X}_2)d\bm{E}_{\bm{X}_2}(\lambda_2)\nonumber \\
\bm{A}_2&\define&\int\limits_{\lambda_1 \in \sigma(\bm{X}_1)}\int\limits_{\lambda_2 \in \sigma(\bm{X}_2)}f^{[1]}(\lambda_{1}, \lambda_{2})\sum_{q_2=1}^{m_{\lambda_2}-1}\frac{(f^{[1]})^{(-,q_2)}(\lambda_{1},\lambda_{2})}{q_2!}\nonumber \\
&& \times d\bm{E}_{\bm{X}_1}(\lambda_1)(\bm{X}_1-\bm{X}_2)(\bm{X}_2-\lambda_2\bm{I})^{q_2}d\bm{E}_{\bm{X}_2}(\lambda_2) \nonumber \\
\bm{A}_3&\define&\int\limits_{\lambda_1 \in \sigma(\bm{X}_1)}\int\limits_{\lambda_2 \in \sigma(\bm{X}_2)}f^{[1]}(\lambda_{1}, \lambda_{2})\sum_{q_1=1}^{m_{\lambda_1}-1}\frac{(f^{[1]})^{(q_1,-)}(\lambda_{1},\lambda_{2})}{q_1!}\nonumber \\
&&\times (\bm{X}_1-\lambda_1\bm{I})^{q_1}d\bm{E}_{\bm{X}_1}(\lambda_1) (\bm{X}_1-\bm{X}_2)d\bm{E}_{\bm{X}_2}(\lambda_2)\nonumber \\
\bm{A}_4&\define&\int\limits_{\lambda_1 \in \sigma(\bm{X}_1)}\int\limits_{\lambda_2 \in \sigma(\bm{X}_2)}f^{[1]}(\lambda_{1}, \lambda_{2})\sum_{q_1=1}^{m_{\lambda_1}-1}\sum_{q_2=1}^{m_{\lambda_2}-1}\frac{(f^{[1]})^{(q_1,q_2)}(\lambda_{1},\lambda_{2})}{q_1!q_2!}\nonumber \\
&& \times (\bm{X}_1-\lambda_1\bm{I})^{q_1}d\bm{E}_{\bm{X}_1}(\lambda_1)(\bm{X}_1-\bm{X}_2)(\bm{X}_2-\lambda_2\bm{I})^{q_2}d\bm{E}_{\bm{X}_2}(\lambda_2).
\end{eqnarray}
then, we have the lower bound for the operator norm of $f(\bm{X}_1) - f(\bm{X}_2)$, which is given by
\begin{eqnarray}\label{eq4:thm:Lipschitz Estimations}
\left\Vert f(\bm{X}_1) - f(\bm{X}_2) \right\Vert\geq \max\left[0, \left\Vert\bm{A}_{\sigma(1)}\right\Vert - \sum\limits_{i=2}^4 \left\Vert\bm{A}_{\sigma(i)}\right\Vert\right],
\end{eqnarray}
where $\sigma$ is the permutation of operators $\bm{A}_i$ for $i=1,2,3,4$ such that $\left\Vert\bm{A}_{\sigma(1)}\right\Vert \geq \left\Vert\bm{A}_{\sigma(2)}\right\Vert \geq\left\Vert\bm{A}_{\sigma(3)}\right\Vert \geq\left\Vert\bm{A}_{\sigma(4)}\right\Vert \geq$. 

Further, if we have $\left[\min\limits_{\lambda_1 \in \sigma(\bm{X}_1),\lambda_2 \in \sigma(\bm{X}_2)} f^{[1]}(\lambda_1,\lambda_2)\right]\left\Vert\bm{X}_1 - \bm{X}_2\right\Vert\geq \left\Vert\bm{A}_{2}\right\Vert+\left\Vert\bm{A}_{3}\right\Vert + \left\Vert\bm{A}_{4}\right\Vert$,  the lower bound for the operator norm of $f(\bm{X}_1) - f(\bm{X}_2)$ can be expressed by
\begin{eqnarray}\label{eq5:thm:Lipschitz Estimations}
\left\Vert f(\bm{X}_1) - f(\bm{X}_2)\right\Vert\geq \left[\min\limits_{\lambda_1 \in \sigma(\bm{X}_1),\lambda_2 \in \sigma(\bm{X}_2)} f^{[1]}(\lambda_1,\lambda_2)\right]\left\Vert\bm{X}_1-\bm{X}_2\right\Vert - \sum\limits_{i=2}^4 \left\Vert\bm{A}_{i}\right\Vert.
\end{eqnarray}
\end{theorem}
\textbf{Proof:}
From Theorem~\ref{thm:pert formula},  we have
\begin{eqnarray}
f(\bm{X}_1)-f(\bm{X}_2)&=&T_{\frac{f(x_1) - f(x_2)}{x_1 -x_2}}^{\bm{X}_{1},\bm{X}_{2}}(\bm{X}_{1}-\bm{X}_{2}). 
\end{eqnarray}
This theorem is obtained by applying in $\beta = f^{[1]}$ and $\bm{Y}=\bm{X}_1 - \bm{X}_2$ in Theorem~\ref{thm: GDOI norm}.
$\hfill\Box$

\section{Continuity of GDOI for Continuous Spectrum Operators}\label{sec:Continuity of GDOI for Continuous Spectrum Operators}

In this section, we prove the continuity property of the GDOI $T_{\beta}^{\bm{X}_1,\bm{X}_2}(\bm{Y})$.


Follow the basic spirit of GDOI, we will define Generalized Triple Operator Integral (GTOI) for continuous spectrum operators in Eq.~\eqref{eq1:  GTOI def}. We have
\begin{eqnarray}\label{eq1:  GTOI def}
\lefteqn{T_{\beta}^{\bm{X}_1,\bm{X}_2,\bm{X}_3}(\bm{Y}_1, \bm{Y}_2)\define}\nonumber\\
&&\int\limits_{\lambda_1 \in \sigma(\bm{X}_1)}\int\limits_{\lambda_2 \in \sigma(\bm{X}_2)}\int\limits_{\lambda_3 \in \sigma(\bm{X}_3)}\beta(\lambda_{1}, \lambda_{2}, \lambda_{3})d\bm{E}_{\bm{X}_1}(\lambda_1) \bm{Y}_1 d\bm{E}_{\bm{X}_2}(\lambda_2) \bm{Y}_2 d\bm{E}_{\bm{X}_3}(\lambda_3)  \nonumber \\
&&+\int\limits_{\lambda_1 \in \sigma(\bm{X}_1)}\int\limits_{\lambda_2 \in \sigma(\bm{X}_2)}\int\limits_{\lambda_3 \in \sigma(\bm{X}_3)}\sum_{q_3=1}^{m_{\lambda_3}-1}\frac{\beta^{(-,-,q_3)}(\lambda_{1}, \lambda_{2}, \lambda_{3})}{q_3!}\nonumber \\
&&~~ \times d\bm{E}_{\bm{X}_1}(\lambda_1)\bm{Y}_1d\bm{E}_{\bm{X}_2}(\lambda_2)\bm{Y}_2(\bm{X}_3-\lambda_3\bm{I})^{q_3}d\bm{E}_{\bm{X}_3}(\lambda_3)\nonumber \\
&&+\int\limits_{\lambda_1 \in \sigma(\bm{X}_1)}\int\limits_{\lambda_2 \in \sigma(\bm{X}_2)}\int\limits_{\lambda_3 \in \sigma(\bm{X}_3)}\sum_{q_2=1}^{m_{\lambda_2}-1}\frac{\beta^{(-,q_2,-)}(\lambda_{1}, \lambda_{2}, \lambda_{3})}{q_2!}\nonumber \\
&&~~ \times d\bm{E}_{\bm{X}_1}(\lambda_1) \bm{Y}_1 (\bm{X}_2-\lambda_2\bm{I})^{q_2}d\bm{E}_{\bm{X}_2}(\lambda_2) \bm{Y}_2 d\bm{E}_{\bm{X}_3}(\lambda_3)  \nonumber \\
&&+\int\limits_{\lambda_1 \in \sigma(\bm{X}_1)}\int\limits_{\lambda_2 \in \sigma(\bm{X}_2)}\int\limits_{\lambda_3 \in \sigma(\bm{X}_3)}\sum_{q_1=1}^{m_{\lambda_1}-1}\frac{\beta^{(q_1,-,-)}(\lambda_{1}, \lambda_{2}, \lambda_{3})}{q_1!}\nonumber \\
&&~~ \times (\bm{X}_1-\lambda_1\bm{I})^{q_1}d\bm{E}_{\bm{X}_1}(\lambda_1) \bm{Y}_1 d\bm{E}_{\bm{X}_2}(\lambda_2) \bm{Y}_2 d\bm{E}_{\bm{X}_3}(\lambda_3)  \nonumber \\
&&+\int\limits_{\lambda_1 \in \sigma(\bm{X}_1)}\int\limits_{\lambda_2 \in \sigma(\bm{X}_2)}\int\limits_{\lambda_3 \in \sigma(\bm{X}_3)}\sum_{q_2=1}^{m_{\lambda_2}-1}\sum_{q_3=1}^{m_{\lambda_3}-1}\frac{\beta^{(-,q_2,q_3)}(\lambda_{1}, \lambda_{2}, \lambda_{3})}{q_2! q_3!} \nonumber \\
&&~~ \times d\bm{E}_{\bm{X}_1}(\lambda_1) \bm{Y}_1  (\bm{X}_2-\lambda_2\bm{I})^{q_2}d\bm{E}_{\bm{X}_2}(\lambda_2) \bm{Y}_2  (\bm{X}_3-\lambda_3\bm{I})^{q_3}d\bm{E}_{\bm{X}_3}(\lambda_3) \nonumber \\
&&+\int\limits_{\lambda_1 \in \sigma(\bm{X}_1)}\int\limits_{\lambda_2 \in \sigma(\bm{X}_2)}\int\limits_{\lambda_3 \in \sigma(\bm{X}_3)}\sum_{q_1=1}^{m_{\lambda_1}-1}\sum_{q_3=1}^{m_{\lambda_3}-1}\frac{\beta^{(q_1,-,q_3)}(\lambda_{1}, \lambda_{2}, \lambda_{3})}{q_1!q_3!} \nonumber \\
&& ~~ \times (\bm{X}_1-\lambda_1\bm{I})^{q_1}d\bm{E}_{\bm{X}_1}(\lambda_1)  \bm{Y}_1 d\bm{E}_{\bm{X}_2}(\lambda_2)  \bm{Y}_2  (\bm{X}_3-\lambda_3\bm{I})^{q_3}d\bm{E}_{\bm{X}_3}(\lambda_3)  \nonumber \\
&&+\int\limits_{\lambda_1 \in \sigma(\bm{X}_1)}\int\limits_{\lambda_2 \in \sigma(\bm{X}_2)}\int\limits_{\lambda_3 \in \sigma(\bm{X}_3)}\sum_{q_1=1}^{m_{\lambda_1}-1}\sum_{q_2=1}^{m_{\lambda_2}-1}\frac{\beta^{(q_1,q_2,-)}(\lambda_{1}, \lambda_{2}, \lambda_{3})}{q_1!q_2!}  \nonumber \\
&& ~~ \times (\bm{X}_1-\lambda_1\bm{I})^{q_1}d\bm{E}_{\bm{X}_1}(\lambda_1) \bm{Y}_1 (\bm{X}_2-\lambda_2\bm{I})^{q_2}d\bm{E}_{\bm{X}_2}(\lambda_2) \bm{Y}_2 d\bm{E}_{\bm{X}_3}(\lambda_3) \nonumber \\
&&+\int\limits_{\lambda_1 \in \sigma(\bm{X}_1)}\int\limits_{\lambda_2 \in \sigma(\bm{X}_2)}\int\limits_{\lambda_3 \in \sigma(\bm{X}_3)}\sum_{q_1=1}^{m_{\lambda_1}-1}\sum_{q_2=1}^{m_{\lambda_2}-1}\sum_{q_3=1}^{m_{\lambda_3}-1}\frac{\beta^{(q_1,q_2,q_3)}(\lambda_{1}, \lambda_{2}, \lambda_{3})}{q_1!q_2!q_3!}\nonumber \\
&&~~\times (\bm{X}_1-\lambda_1\bm{I})^{q_1}d\bm{E}_{\bm{X}_1}(\lambda_1) \bm{Y}_1 (\bm{X}_2-\lambda_2\bm{I})^{q_2}d\bm{E}_{\bm{X}_2}(\lambda_2) \bm{Y}_2 (\bm{X}_3-\lambda_3\bm{I})^{q_3}d\bm{E}_{\bm{X}_3}(\lambda_3) .
\end{eqnarray}

From the definition of GTOI given by Eq.~\eqref{eq1:  GTOI def}, we have the following Theorem~\ref{thm: GTOI norm} about the upper and the lower bounds for GTOI.
\begin{theorem}\label{thm: GTOI norm}
We have the upper bound for the norm of $T_{\beta}^{\bm{X}_1,\bm{X}_2,\bm{X}_3}(\bm{Y}_1,\bm{Y}_2)$, which is given by
\begin{eqnarray}\label{eq1: thm: GTOI norm}
\lefteqn{\left\Vert T_{\beta}^{\bm{X}_1,\bm{X}_2,\bm{X}_3}(\bm{Y}_1,\bm{Y}_2)\right\Vert\leq}\nonumber \\
&&\left[\max\limits_{\lambda_1 \in \sigma(\bm{X}_1),\lambda_2 \in \sigma(\bm{X}_2),\lambda_3 \in \Lambda_{\bm{X}_3}} \left\vert\beta(\lambda_1,\lambda_2,\lambda_3)\right\vert\right]\left\Vert\bm{Y}_1\right\Vert\left\Vert\bm{Y}_2\right\Vert\nonumber \\
&&+\Bigg\{\int\limits_{\lambda_3 \in \sigma(\bm{X}_3)}\sum_{q_3=1}^{m_{\lambda_3}-1}\left[\max\limits_{\lambda_1 \in \sigma(\bm{X}_1),\lambda_2 \in \sigma(\bm{X}_2)\lambda_3 \in \Lambda_{\bm{X}_3}}\left\vert\frac{\beta^{(-,-,q_3)}(\lambda_{1},\lambda_{2},\lambda_3)}{q_3!}\right\vert\right] \nonumber \\
&& \times \left\Vert (\bm{X}_3-\lambda_3\bm{I})^{q_3}\right\Vert d\bm{E}_{\bm{X}_3}(\lambda_3)\Bigg\} \left\Vert\bm{Y}_1\right\Vert\left\Vert\bm{Y}_2\right\Vert \nonumber \\
&&+\Bigg\{\int\limits_{\lambda_2 \in \sigma(\bm{X}_2)}\sum_{q_2=1}^{m_{\lambda_2}-1}\left[\max\limits_{\lambda_1 \in \sigma(\bm{X}_1),\lambda_2 \in \sigma(\bm{X}_2)\lambda_3 \in \Lambda_{\bm{X}_3}}\left\vert\frac{\beta^{(-,q_2,-)}(\lambda_{1},\lambda_{2},\lambda_3)}{q_2!}\right\vert\right] \nonumber \\
&& \times \left\Vert (\bm{X}_2-\lambda_2\bm{I})^{q_2}\right\Vert d\bm{E}_{\bm{X}_2}(\lambda_2)\Bigg\} \left\Vert\bm{Y}_1\right\Vert\left\Vert\bm{Y}_2\right\Vert \nonumber \\
&&+\Bigg\{\int\limits_{\lambda_1 \in \sigma(\bm{X}_1)}\sum_{q_1=1}^{m_{\lambda_1}-1}\left[\max\limits_{\lambda_1 \in \sigma(\bm{X}_1),\lambda_2 \in \sigma(\bm{X}_2)\lambda_1 \in \sigma(\bm{X}_1)}\left\vert\frac{\beta^{(q_1,-,-)}(\lambda_{1},\lambda_{2},\lambda_3)}{q_!}\right\vert\right] \nonumber \\
&& \times \left\Vert (\bm{X}_1-\lambda_1\bm{I})^{q_1}\right\Vert d\bm{E}_{\bm{X}_1}(\lambda_1)\Bigg\} \left\Vert\bm{Y}_1\right\Vert\left\Vert\bm{Y}_2\right\Vert\nonumber \\
&&+\Bigg\{\int\limits_{\lambda_2 \in \sigma(\bm{X}_2)}\int\limits_{\lambda_3 \in \sigma(\bm{X}_3)}\sum_{q_2=1}^{m_{\lambda_2}-1}\sum_{q_3=1}^{m_{\lambda_3}-1}\left[\max\limits_{\lambda_1 \in \sigma(\bm{X}_1),\lambda_2 \in \sigma(\bm{X}_2),\lambda_3 \in \Lambda_{\bm{X}_3}}\left\vert\frac{\beta^{(-,q_2,q_3)}(\lambda_{1},\lambda_{2},\lambda_3)}{q_2!q_3!}\right\vert\right]\nonumber\\
&&~~~\times \left\Vert (\bm{X}_2-\lambda_2\bm{I})^{q_2} \right\Vert  \left\Vert (\bm{X}_3-\lambda_3\bm{I})^{q_3} \right\Vert d\bm{E}_{\bm{X}_2}(\lambda_2)d\bm{E}_{\bm{X}_3}(\lambda_3)\Bigg\}\left\Vert\bm{Y}_1\right\Vert\left\Vert\bm{Y}_2\right\Vert\nonumber \\
&&+\Bigg\{\int\limits_{\lambda_1 \in \sigma(\bm{X}_1)}\int\limits_{\lambda_3 \in \sigma(\bm{X}_3)}\sum_{q_1=1}^{m_{\lambda_1}-1}\sum_{q_3=1}^{m_{\lambda_3}-1}\left[\max\limits_{\lambda_1 \in \sigma(\bm{X}_1),\lambda_2 \in \sigma(\bm{X}_2),\lambda_3 \in \Lambda_{\bm{X}_3}}\left\vert\frac{\beta^{(q_1,-,q_3)}(\lambda_{1},\lambda_{2},\lambda_3)}{q_1!q_3!}\right\vert\right]\nonumber\\
&&~~~\times \left\Vert (\bm{X}_1-\lambda_1\bm{I})^{q_1} \right\Vert  \left\Vert (\bm{X}_3-\lambda_3\bm{I})^{q_3} \right\Vert d\bm{E}_{\bm{X}_1}(\lambda_1)d\bm{E}_{\bm{X}_3}(\lambda_3)\Bigg\}\left\Vert\bm{Y}_1\right\Vert\left\Vert\bm{Y}_2\right\Vert \nonumber \\
&&+\Bigg\{\int\limits_{\lambda_1 \in \sigma(\bm{X}_1)}\int\limits_{\lambda_2 \in \sigma(\bm{X}_2)}\sum_{q_1=1}^{m_{\lambda_1}-1}\sum_{q_2=1}^{m_{\lambda_2}-1}\left[\max\limits_{\lambda_1 \in \sigma(\bm{X}_1),\lambda_2 \in \sigma(\bm{X}_2),\lambda_3 \in \Lambda_{\bm{X}_3}}\left\vert\frac{\beta^{(q_1,q_2,-)}(\lambda_{1},\lambda_{2},\lambda_3)}{q_1!q_2!}\right\vert\right]\nonumber\\
&&~~~\times \left\Vert (\bm{X}_1-\lambda_1\bm{I})^{q_1} \right\Vert  \left\Vert (\bm{X}_2-\lambda_2\bm{I})^{q_2} \right\Vert d\bm{E}_{\bm{X}_1}(\lambda_1)d\bm{E}_{\bm{X}_2}(\lambda_2)\Bigg\}\left\Vert\bm{Y}_1\right\Vert\left\Vert\bm{Y}_2\right\Vert  \nonumber \\
&&+\Bigg\{\int\limits_{\lambda_1 \in \sigma(\bm{X}_1)}\int\limits_{\lambda_2 \in \sigma(\bm{X}_2)}\int\limits_{\lambda_3 \in \sigma(\bm{X}_3)}\sum_{q_1=1}^{m_{\lambda_1}-1}\sum_{q_2=1}^{m_{\lambda_2}-1} \sum_{q_3=1}^{m_{\lambda_3}-1}\left[\max\limits_{\lambda_1 \in \sigma(\bm{X}_1),\lambda_2 \in \sigma(\bm{X}_2),\lambda_3 \in \Lambda_{\bm{X}_3}}\left\vert\frac{\beta^{(q_1,q_2,q_3)}(\lambda_{1},\lambda_{2},\lambda_3)}{q_1!q_2!q_3!}\right\vert\right]\nonumber\\
&&~~~\times \left\Vert (\bm{X}_1-\lambda_1\bm{I})^{q_1} \right\Vert  \left\Vert (\bm{X}_2-\lambda_2\bm{I})^{q_2} \right\Vert  \left\Vert (\bm{X}_3-\lambda_2\bm{I})^{q_3} \right\Vert d\bm{E}_{\bm{X}_1}(\lambda_1) d\bm{E}_{\bm{X}_2}(\lambda_2) d\bm{E}_{\bm{X}_3}(\lambda_3) \Bigg\}\left\Vert\bm{Y}_1\right\Vert\left\Vert\bm{Y}_2\right\Vert, \nonumber \\
\end{eqnarray}
where $\sigma(\bm{X}_1), \sigma(\bm{X}_2)$ and $\sigma(\bm{X}_3)$ are spectrums of operators $\bm{X}_1, \bm{X}_2$ and $\bm{X}_3$, respectively.

On the other hand, let us define the following operators
\begin{eqnarray}\label{eq1-1: thm: GTOI norm}
\bm{A}_1&\define&\int\limits_{\lambda_1 \in \sigma(\bm{X}_1)}\int\limits_{\lambda_2 \in \sigma(\bm{X}_2)}\int\limits_{\lambda_3 \in \sigma(\bm{X}_3)}\beta(\lambda_{1}, \lambda_{2}, \lambda_{3})d\bm{E}_{\bm{X}_1}(\lambda_1) \bm{Y}_1 d\bm{E}_{\bm{X}_2}(\lambda_2) \bm{Y}_2 d\bm{E}_{\bm{X}_3}(\lambda_3) \nonumber \\
\bm{A}_2&\define&\int\limits_{\lambda_1 \in \sigma(\bm{X}_1)}\int\limits_{\lambda_2 \in \sigma(\bm{X}_2)}\int\limits_{\lambda_3 \in \sigma(\bm{X}_3)}\sum_{q_3=1}^{m_{\lambda_3}-1}\frac{\beta^{(-,-,q_3)}(\lambda_{1}, \lambda_{2}, \lambda_{3})}{q_3!} \nonumber \\
&& \times d\bm{E}_{\bm{X}_1}(\lambda_1)\bm{Y}_1 d\bm{E}_{\bm{X}_2}(\lambda_2) \bm{Y}_2 (\bm{X}_3-\lambda_3\bm{I})^{q_3} d\bm{E}_{\bm{X}_3}(\lambda_3) \nonumber \\
\bm{A}_3&\define&\int\limits_{\lambda_1 \in \sigma(\bm{X}_1)}\int\limits_{\lambda_2 \in \sigma(\bm{X}_2)}\int\limits_{\lambda_3 \in \sigma(\bm{X}_3)}\sum_{q_2=1}^{m_{\lambda_2}-1}\frac{\beta^{(-,q_2,-)}(\lambda_{1}, \lambda_{2}, \lambda_{3})}{q_2!}\nonumber \\
&& \times d\bm{E}_{\bm{X}_1}(\lambda_1) \bm{Y}_1 (\bm{X}_2-\lambda_2\bm{I})^{q_2} d\bm{E}_{\bm{X}_2}(\lambda_2)  \bm{Y}_2 d\bm{E}_{\bm{X}_3}(\lambda_3)\nonumber \\
\bm{A}_4&\define&\int\limits_{\lambda_1 \in \sigma(\bm{X}_1)}\int\limits_{\lambda_2 \in \sigma(\bm{X}_2)}\int\limits_{\lambda_3 \in \sigma(\bm{X}_3)}\sum_{q_1=1}^{m_{\lambda_1}-1}\frac{\beta^{(q_1,-,-)}(\lambda_{1}, \lambda_{2}, \lambda_{3})}{q_1!}\nonumber \\
&& \times  (\bm{X}_1-\lambda_1\bm{I})^{q_1} d\bm{E}_{\bm{X}_1}(\lambda_1) \bm{Y}_1 d\bm{E}_{\bm{X}_2}(\lambda_2) \bm{Y}_2 d\bm{E}_{\bm{X}_3}(\lambda_3) \nonumber \\
\bm{A}_5&\define&\int\limits_{\lambda_1 \in \sigma(\bm{X}_1)}\int\limits_{\lambda_2 \in \sigma(\bm{X}_2)}\int\limits_{\lambda_3 \in \sigma(\bm{X}_3)}\sum_{q_2=1}^{m_{\lambda_2}-1}\sum_{q_3=1}^{m_{\lambda_3}-1}\frac{\beta^{(-,q_2,q_3)}(\lambda_{1}, \lambda_{2}, \lambda_{3})}{q_2! q_3!}\nonumber \\
&& \times d\bm{E}_{\bm{X}_1}(\lambda_1) \bm{Y}_1 (\bm{X}_2-\lambda_2\bm{I})^{q_2} d\bm{E}_{\bm{X}_2}(\lambda_2)\bm{Y}_2 (\bm{X}_3-\lambda_3\bm{I})^{q_3} d\bm{E}_{\bm{X}_3}(\lambda_3) \nonumber \\
\bm{A}_6&\define&\int\limits_{\lambda_1 \in \sigma(\bm{X}_1)}\int\limits_{\lambda_2 \in \sigma(\bm{X}_2)}\int\limits_{\lambda_3 \in \sigma(\bm{X}_3)}\sum_{q_1=1}^{m_{\lambda_1}-1}\sum_{q_3=1}^{m_{\lambda_3}-1}\frac{\beta^{(q_1,-,q_3)}(\lambda_{1}, \lambda_{2}, \lambda_{3})}{q_1!q_3!}\nonumber \\
&& \times (\bm{X}_1-\lambda_1\bm{I})^{q_1} d\bm{E}_{\bm{X}_1}(\lambda_1) \bm{Y}_1 d\bm{E}_{\bm{X}_2}(\lambda_2) \bm{Y}_2 (\bm{X}_3-\lambda_3\bm{I})^{q_3} d\bm{E}_{\bm{X}_3}(\lambda_3)  \nonumber \\
\bm{A}_7&\define&\int\limits_{\lambda_1 \in \sigma(\bm{X}_1)}\int\limits_{\lambda_2 \in \sigma(\bm{X}_2)}\int\limits_{\lambda_3 \in \sigma(\bm{X}_3)}\sum_{q_1=1}^{m_{\lambda_1}-1}\sum_{q_2=1}^{m_{\lambda_2}-1}\frac{\beta^{(q_1,q_2,-)}(\lambda_{1}, \lambda_{2}, \lambda_{3})}{q_1!q_2!}\nonumber \\
&& \times  (\bm{X}_1-\lambda_1\bm{I})^{q_1} d\bm{E}_{\bm{X}_1}(\lambda_1) \bm{Y}_1  (\bm{X}_2-\lambda_2\bm{I})^{q_2} d\bm{E}_{\bm{X}_2}(\lambda_2) \bm{Y}_2 d\bm{E}_{\bm{X}_3}(\lambda_3) \nonumber \\
\bm{A}_8&\define&\int\limits_{\lambda_1 \in \sigma(\bm{X}_1)}\int\limits_{\lambda_2 \in \sigma(\bm{X}_2)}\int\limits_{\lambda_3 \in \sigma(\bm{X}_3)}\sum_{q_1=1}^{m_{\lambda_1}-1}\sum_{q_2=1}^{m_{\lambda_2}-1}\sum_{q_3=1}^{m_{k_3,i3}-1}\frac{\beta^{(q_1,q_2,q_3)}(\lambda_{1}, \lambda_{2}, \lambda_{3})}{q_1!q_2!q_3!}\nonumber \\
&&~~\times(\bm{X}_1-\lambda_1\bm{I})^{q_1} d\bm{E}_{\bm{X}_1}(\lambda_1) \bm{Y}_1 (\bm{X}_2-\lambda_2\bm{I})^{q_2} d\bm{E}_{\bm{X}_2}(\lambda_2) \bm{Y}_2 (\bm{X}_3-\lambda_3\bm{I})^{q_3} d\bm{E}_{\bm{X}_3}(\lambda_3),
\end{eqnarray}
then, we have the lower bound for the operator norm of $T_{\beta}^{\bm{X}_1,\bm{X}_2,\bm{X}_3}(\bm{Y}_1,\bm{Y}_2)$, which is given by
\begin{eqnarray}\label{eq2: thm: GTOI norm}
\left\Vert T_{\beta}^{\bm{X}_1,\bm{X}_2,\bm{X}_3}(\bm{Y}_1,\bm{Y}_2) \right\Vert\geq \max\left[0, \left\Vert\bm{A}_{\sigma(1)}\right\Vert - \sum\limits_{i=2}^8 \left\Vert\bm{A}_{\sigma(i)}\right\Vert\right],
\end{eqnarray}
where $\sigma$ is the permutation of matrices $\bm{A}_i$ for $i=1,2,3,4,5,6,7,8$ such that $\left\Vert\bm{A}_{\sigma(1)}\right\Vert \geq \ldots \geq \left\Vert\bm{A}_{\sigma(8)}\right\Vert \geq$. 

Further, if we have $\left[\min\limits_{\lambda_1 \in \sigma(\bm{X}_1),\lambda_2 \in \sigma(\bm{X}_2),\lambda_3 \in \Lambda_{\bm{X}_3}} \beta(\lambda_1,\lambda_2,\lambda_3)\right]\left\Vert\bm{Y}_1\right\Vert \left\Vert\bm{Y}_2\right\Vert \geq \sum\limits_{i=2}^8\left\Vert\bm{A}_{i}\right\Vert$,  the lower bound for the norm of $T_{\beta}^{\bm{X}_1,\bm{X}_2,\bm{X}_3}(\bm{Y}_1,\bm{Y}_2)$ can be expressed by
\begin{eqnarray}\label{eq3: thm: GTOI norm}
\left\Vert T_{\beta}^{\bm{X}_1,\bm{X}_2,\bm{X}_3}(\bm{Y}_1,\bm{Y}_2) \right\Vert\geq \left[\min\limits_{\lambda_1 \in \sigma(\bm{X}_1),\lambda_2 \in \sigma(\bm{X}_2),\lambda_3 \in \Lambda_{\bm{X}_3}} \beta(\lambda_1,\lambda_2,\lambda_3)\right]\left\Vert\bm{Y}_1\right\Vert\left\Vert\bm{Y}_2\right\Vert - \sum\limits_{i=2}^8 \left\Vert\bm{A}_{i}\right\Vert.
\end{eqnarray}
\end{theorem}
\textbf{Proof:}
The proof is similar to the proof in Theorem~\ref{thm: GDOI norm}.
$\hfill\Box$

In the following Lemma~\ref{lma: telescope propt}, we establish the telescope property for GTOI in the setting of continuous spectrum operators, thereby generalizing the classical MOI framework, which has been primarily developed for Hermitian or self-adjoint parameter matrices~\cite{skripka2019multilinear}.
\begin{lemma}\label{lma: telescope propt}
We have
\begin{eqnarray}\label{eq1:lma: telescope propt}
T_{f^{[1]}}^{\bm{A},\bm{X}}(\bm{Y}) -  T_{f^{[1]}}^{\bm{B},\bm{X}}(\bm{Y})&=&
T_{f^{[2]}}^{\bm{A},\bm{B},\bm{X}}(\bm{A}-\bm{B},\bm{Y}),
\end{eqnarray}
where $f^{[1]}(x_0,x_1)$ and $^{[2]}(x_0,x_1,x_2)$ are first and second divide differences. 
\end{lemma}
\textbf{Proof:}
We have 
\begin{eqnarray}\label{eq2:lma: telescope propt}
T_{f^{[2]}}^{\bm{A},\bm{B},\bm{X}}(\bm{A}-\bm{B},\bm{Y})&=&T_{f^{[2]}}^{\bm{A},\bm{B},\bm{X}}(\bm{A},\bm{Y}) - T_{f^{[2]}}^{\bm{A},\bm{B},\bm{X}}(\bm{B},\bm{Y})\nonumber \\
&=&T_{x_0 f^{[2]}}^{\bm{A},\bm{B},\bm{X}}(\bm{I},\bm{Y}) - T_{x_1 f^{[2]}}^{\bm{A},\bm{B},\bm{X}}(\bm{I},\bm{Y})\nonumber \\
&=&T_{x_0 f^{[2]}-x_1 f^{[2]}}^{\bm{A},\bm{B},\bm{X}}(\bm{I},\bm{Y})\nonumber \\
&=_1&T_{f^{[1]}(x_0,x_2)}^{\bm{A},\bm{X}}(\bm{Y}) - T_{f^{[1]}(x_1,x_2)}^{\bm{B},\bm{X}}(\bm{Y})\nonumber \\
&=&T_{f^{[1]}}^{\bm{A},\bm{X}}(\bm{Y}) -  T_{f^{[1]}}^{\bm{B},\bm{X}}(\bm{Y}),
\end{eqnarray}
where we have $x_0 f^{[2]}-x_1 f^{[2]}=f^{[1]}(x_0,x_2)-f^{[1]}(x_1,x_2)$ in $=_1$.
$\hfill\Box$

The following Theorem~\ref{thm:GDOI continuity dd} will show the continuity property of the GDOI $T_{f^{[1]}}^{\bm{X}_1,\bm{X}_2}(\bm{Y})$  in the setting of continuous spectrum operators. 
\begin{theorem}\label{thm:GDOI continuity dd}
Given two sequence of operators $\bm{X}_{1,\ell_1}$ and $\bm{X}_{2,\ell_2}$ satisfying $\bm{X}_{1,\ell_1} \rightarrow \bm{X}_{1}$ and $\bm{X}_{2,\ell_2} \rightarrow \bm{X}_{2}$, respectively. Moreover, all spectrums of operators $\bm{X}_{1,\ell_1}, \bm{X}_1, \bm{X}_{2,\ell_2}$ and $\bm{X}_{2}$ are assumed bounded. If we have $f^{[2]}(\lambda_1, \lambda_2, \lambda_3)$ satisfying:
\begin{eqnarray}\label{eq0:thm:GDOI continuity dd}
\forall  \bm{\vartheta} = (\vartheta_1, \vartheta_2, \vartheta_3) \in \mathbb{N}_0^3, \quad \left| \frac{\partial^{|\bm{\vartheta}|}}{\partial \lambda_1^{\vartheta_1} \partial \lambda_2^{\vartheta_2} \partial \lambda_3^{\vartheta_3}} f^{[2]}(\lambda_1, \lambda_2, \lambda_3) \right| < \infty
\end{eqnarray}
where $\bm{\vartheta}$ is a multi-index  $(\vartheta_1, \vartheta_2, \vartheta_3)$, $|\bm{\vartheta}| = \vartheta_1 + \vartheta_2 + \vartheta_3$ is the otal order of derivatives, and $\mathbb{N}_0$ is the set of non-negative integers.

Then, we have 
\begin{eqnarray}\label{eq1:thm:GDOI continuity dd}
T_{f^{[1]}}^{\bm{X}_{1,\ell_1},\bm{X}_{2,\ell_2}}(\bm{Y}) \rightarrow T_{f^{[1]}}^{\bm{X}_{1},\bm{X}_{2}}(\bm{Y}),
\end{eqnarray}
where $\rightarrow$ is in the sense of operator norm, i.e., 
\begin{eqnarray}\label{eq2:thm:GDOI continuity dd}
\lim\limits_{\ell_1, \ell_2 \rightarrow \infty}\left\Vert T_{f^{[1]}}^{\bm{X}_{1,\ell_1},\bm{X}_{2,\ell_2}}(\bm{Y}) - T_{f^{[1]}}^{\bm{X}_{1},\bm{X}_{2}}(\bm{Y})\right\Vert=0
\end{eqnarray}
\end{theorem}
\textbf{Proof:}
Because we have
\begin{eqnarray}\label{eq3:thm:GDOI continuity dd}
\lefteqn{\left\Vert T_{f^{[1]}}^{\bm{X}_{1,\ell_1},\bm{X}_{2,\ell_2}}(\bm{Y}) - T_{f^{[1]}}^{\bm{X}_{1},\bm{X}_{2}}(\bm{Y})\right\Vert}\nonumber \\
&\leq&\left\Vert T_{f^{[1]}}^{\bm{X}_{1,\ell_1},\bm{X}_{2,\ell_2}}(\bm{Y}) - T_{f^{[1]}}^{\bm{X}_{1},\bm{X}_{2,\ell_2}}(\bm{Y})\right\Vert + \left\Vert T_{f^{[1]}}^{\bm{X}_{1},\bm{X}_{2,\ell_2}}(\bm{Y}) - T_{f^{[1]}}^{\bm{X}_{1},\bm{X}_{2}}(\bm{Y})\right\Vert\nonumber \\
&=_1&\left\Vert T_{f^{[2]}}^{\bm{X}_{1,\ell_1},\bm{X}_{1},\bm{X}_{2,\ell_2}}(\bm{X}_{1,\ell_1}-\bm{X}_1,\bm{Y})\right\Vert + \left\Vert T_{f^{[2]}}^{\bm{X}_{1},\bm{X}_{2,\ell_2},\bm{X}_2}(\bm{Y}, \bm{X}_{2,\ell_2}-\bm{X}_2)\right\Vert\nonumber \\
&\leq_2& \epsilon/2 +\epsilon/2 = \epsilon
\end{eqnarray}
where we apply Lemma~\ref{lma: telescope propt} at step $=1$. Since the sequences $\bm{X}{1,\ell_1}$ and $\bm{X}_{2,\ell_2}$ converge to $\bm{X}_1$ and $\bm{X}_2$, respectively, and the function $f^{[2]}$ satisfies Eq.\eqref{eq0:thm:GDOI continuity dd}, together with the boundedness of the spectra of $\bm{X}_{1,\ell_1}$, $\bm{X}1$, $\bm{X}{2,\ell_2}$, and $\bm{X}_2$, we can invoke Theorem\ref{thm: GTOI norm} to obtain the inequality $\leq_2$.
$\hfill\Box$

Theorem~\ref{thm:GDOI continuity dd} is restrictive to the underlying function as divided difference only. We wish to establish the continuity property for a more general function rather than just $f^{[1]}$. To this end, we must prepare the following Lemma~\ref{lma:continuity ii and iii} concerning the continuity of the variable matrix $\bm{Y}$ and the underlying function $\beta$ in the GDOI $T_{\beta}^{\bm{X}_1, \bm{X}_2}(\bm{Y})$, within the setting of continuous spectrum operators.
\begin{lemma}\label{lma:continuity ii and iii}
(i) Given a sequence of operators $\bm{Y}_\ell$ such that $\bm{Y}_\ell \rightarrow \bm{Y}$, we assume that 
those terms involving $\max\limits_{\lambda_1 \in \sigma(\bm{X}_1),\lambda_2 \in \sigma(\bm{X}_2)}$ in Eq.~\eqref{eq1: thm: GDOI norm} are finite with bounded spectrums $\sigma(\bm{X}_1)$ and $\sigma(\bm{X}_2)$. If we have $\beta(\lambda_1, \lambda_2)$ satisfying:
\begin{eqnarray}\label{eq0:lma:continuity ii and iii}
\forall  \bm{\vartheta} = (\vartheta_1, \vartheta_2) \in \mathbb{N}_0^2, \quad \left| \frac{\partial^{|\bm{\vartheta}|}}{\partial \lambda_1^{\vartheta_1} \partial \lambda_2^{\vartheta_2}}\beta(\lambda_1, \lambda_2) \right| < \infty
\end{eqnarray}
where $\bm{\vartheta}$ is a multi-index  $(\vartheta_1, \vartheta_2)$, $|\bm{\vartheta}| = \vartheta_1 + \vartheta_2$ is the otal order of derivatives, and $\mathbb{N}_0$ is the set of non-negative integers.

Then, we have
\begin{eqnarray}\label{eq1:lma:continuity ii and iii}
T_{\beta}^{\bm{X}_1, \bm{X}_2}(\bm{Y}_\ell)&\rightarrow&T_{\beta}^{\bm{X}_1, \bm{X}_2}(\bm{Y}).
\end{eqnarray}

(ii) Given a sequence of function $\beta_\ell$ such that ${\beta_\ell}$ converges to $\beta$ in the $C^\infty$-topology on $\sigma(\bm{X}_1) \times \sigma(\bm{X}_2)$, i.e., 
\begin{eqnarray}\label{eq1.1:lma:continuity ii and iii}
\lim\limits_{\ell \rightarrow \infty} \left\Vert \beta_\ell - \beta \right\Vert_{C^k} \rightarrow 0,~~{\mbox{for all $k \geq 0$}.}
\end{eqnarray}
Then, we have
\begin{eqnarray}\label{eq2:lma:continuity ii and iii}
T_{\beta_\ell}^{\bm{X}_1, \bm{X}_2}(\bm{Y})&\rightarrow&T_{\beta}^{\bm{X}_1, \bm{X}_2}(\bm{Y}).
\end{eqnarray}
\end{lemma}
\textbf{Proof:}
Since we have
\begin{eqnarray}\label{eq3:lma:continuity ii and iii}
\left\Vert T_{\beta}^{\bm{X}_1, \bm{X}_2}(\bm{Y}_\ell) - T_{\beta}^{\bm{X}_1, \bm{X}_2}(\bm{Y}) \right\Vert
&=&\left\Vert T_{\beta}^{\bm{X}_1, \bm{X}_2}(\bm{Y}_\ell - \bm{Y}) \right\Vert \nonumber \\
&\leq& \epsilon,
\end{eqnarray}
where $\epsilon$ is any positive number and this inequality comes from the condition provided by Eq.~\eqref{eq0:lma:continuity ii and iii} and Theorem~\ref{thm: GDOI norm} with $\bm{Y}_\ell \rightarrow \bm{Y}$. This proves Part (i).

For Part (ii), we have
\begin{eqnarray}\label{eq4:lma:continuity ii and iii}
\left\Vert T_{\beta_\ell}^{\bm{X}_1, \bm{X}_2}(\bm{Y}) - T_{\beta}^{\bm{X}_1, \bm{X}_2}(\bm{Y}) \right\Vert
&=&\left\Vert T_{\beta_\ell-\beta}^{\bm{X}_1, \bm{X}_2}( \bm{Y}) \right\Vert \nonumber \\
&\leq& \epsilon,
\end{eqnarray}
where $\epsilon$ is any positive number and this inequality comes from the condition provided by Eq.~\eqref{eq1.1:lma:continuity ii and iii} and Theorem~\ref{thm: GDOI norm}. 
$\hfill\Box$

We are ready to present the following GDOI continuity theorem for continuous spectrum operators with the general function $\beta$.
\begin{theorem}\label{thm:GDOI continuity}
Given two sequence of operators $\bm{X}_{1,\ell_1}$ and $\bm{X}_{2,\ell_2}$ satisfying $\bm{X}_{1,\ell_1} \rightarrow \bm{X}_{1}$ and $\bm{X}_{2,\ell_2} \rightarrow \bm{X}_{2}$, respectively. We assume that all spectrums of operators $\bm{X}_{1,\ell_1}, \bm{X}_{2,\ell_2}, \bm{X}_{1}$ and $\bm{X}_{2}$ are within a bounded, polynomially convex compact set. If we also have holomorphic $\beta(\lambda_1, \lambda_2)$ satisfying:
\begin{eqnarray}\label{eq0:thm:GDOI continuity}
\forall  \bm{\vartheta} = (\vartheta_1, \vartheta_2) \in \mathbb{N}_0^2, \quad \left| \frac{\partial^{|\bm{\vartheta}|}}{\partial \lambda_1^{\vartheta_1} \partial \lambda_2^{\vartheta_2}}\beta(\lambda_1, \lambda_2) \right| < \infty
\end{eqnarray}
where $\bm{\vartheta}$ is a multi-index  $(\vartheta_1, \vartheta_2)$, $|\bm{\vartheta}| = \vartheta_1 + \vartheta_2$ is the otal order of derivatives, and $\mathbb{N}_0$ is the set of non-negative integers.

Then, we have 
\begin{eqnarray}\label{eq1:thm:GDOI continuity}
T_{\beta}^{\bm{X}_{1,\ell_1},\bm{X}_{2,\ell_2}}(\bm{Y}) \rightarrow T_{\beta}^{\bm{X}_{1},\bm{X}_{2}}(\bm{Y}),
\end{eqnarray}
where $\rightarrow$ is in the sense of the operator norm, i.e., 
\begin{eqnarray}\label{eq2:thm:GDOI continuity}
\lim\limits_{\ell_1, \ell_2 \rightarrow \infty}\left\Vert T_{\beta}^{\bm{X}_{1,\ell_1},\bm{X}_{2,\ell_2}}(\bm{Y}) - T_{\beta}^{\bm{X}_{1},\bm{X}_{2}}(\bm{Y})\right\Vert=0.
\end{eqnarray}
\end{theorem}
\textbf{Proof:}
For any bivariate polynomial function of degree $\ell$, we can exactly represent it using a linear combination of the divided differences of monomials and the second variable. We have
\begin{eqnarray}\label{eq2-1:thm:GDOI continuity}
x_1^{k_1}x_2^{k_2}&=&\left[f^{[1]}_{k_1 + 1}(x_1,x_2) - x_2 f^{[1]}_{k_1}(x_1,x_2)\right]x_2^{k_2}\nonumber \\
&=&f^{[1]}_{k_1 + 1}(x_1,x_2)x_2^{k_2}-  f^{[1]}_{k_1}(x_1,x_2)x_2^{k_2 + 1},
\end{eqnarray}
where $f^{[1]}_{k_1}(x_1,x_2) \define \frac{x_1^{k_1} - x_2^{k_2}}{x_1 - x_2}$. Then, we have
\begin{eqnarray}\label{eq2-2:thm:GDOI continuity}
\sum\limits_{k_1=0,k_2=0}^{k_1+k_2=\ell} c_{k_1,k_2}x_1^{k_1}x_2^{k_2}&=&\sum\limits_{k_1=0,k_2=0}^{k_1+k_2=\ell}c_{k_1,k_2}\left[f^{[1]}_{k_1 + 1}(x_1,x_2)x_2^{k_2}-  f^{[1]}_{k_1}(x_1,x_2)x_2^{k_2 + 1}\right]\nonumber \\
&=&\sum\limits_{k'_1=1,k'_2=0}^{k'_1+k'_2=\ell+1} d_{k'_1,k'_2}f^{[1]}_{k'_1 }(x_1,x_2)x_2^{k'_2},
\end{eqnarray}
which shows that any bivariate polynomial can be expressed as linear combination of divided difference and its product with the variable $x_2$.

By Oka-Weil approximation theorem, let us use $\tilde{\beta}\define\sum\limits_{k_1,k_2} c_{k_1,k_2}x_1^{k_1}x_2^{k_2}$ as a bivariate polynomial function with complex arguments to approximate $\beta$ uniformly in function values and the derivatives, i.e., $\tilde{\beta}\rightarrow\beta$. Therefore, we have
\begin{eqnarray}\label{eq3:thm:GDOI continuity}
\tilde{\beta}&=&\sum\limits_{k'_1=1,k'_2=0}^{k'_1+k'_2=\ell+1} d_{k'_1,k'_2}f^{[1]}_{k'_1 }(x_1,x_2)x_2^{k'_2},
\end{eqnarray}
where $d_{k'_1,k'_2}$ are complex scalars.  

Then, we have
\begin{eqnarray}\label{eq4:thm:GDOI continuity}
\left\Vert T_{\beta}^{\bm{X}_{1,\ell_1},\bm{X}_{2,\ell_2}}(\bm{Y}) - T_{\beta}^{\bm{X}_{1},\bm{X}_{2}}(\bm{Y})\right\Vert&\leq&\underbrace{\left\Vert T_{\beta}^{\bm{X}_{1,\ell_1},\bm{X}_{2,\ell_2}}(\bm{Y}) - T_{\tilde{\beta}}^{\bm{X}_{1,\ell_1},\bm{X}_{2,\ell_2}}(\bm{Y})\right\Vert}_{\mbox{Part I}}\nonumber \\
&&+\underbrace{\left\Vert T_{\tilde{\beta}}^{\bm{X}_{1,\ell_1},\bm{X}_{2,\ell_2}}(\bm{Y}) - T_{\tilde{\beta}}^{\bm{X}_{1},\bm{X}_{2}}(\bm{Y})\right\Vert}_{\mbox{Part II}}\nonumber \\
&&+\underbrace{\left\Vert T_{\tilde{\beta}}^{\bm{X}_{1},\bm{X}_{2}}(\bm{Y}) - T_{\beta}^{\bm{X}_{1},\bm{X}_{2}}(\bm{Y})\right\Vert}_{\mbox{Part III}}.
\end{eqnarray}

Let us analyze each term in the R.H.S. of Eq.~\eqref{eq4:thm:GDOI continuity}. For Part I, we have
\begin{eqnarray}\label{eq5:thm:GDOI continuity}
\left\Vert T_{\beta}^{\bm{X}_{1,\ell_1},\bm{X}_{2,\ell_2}}(\bm{Y}) - T_{\tilde{\beta}}^{\bm{X}_{1,\ell_1},\bm{X}_{2,\ell_2}}(\bm{Y})\right\Vert \leq \epsilon/3,
\end{eqnarray}
by part (ii) in Lemma~\ref{lma:continuity ii and iii}.  For Part III, we also have
\begin{eqnarray}\label{eq6:thm:GDOI continuity}
\left\Vert T_{\tilde{\beta}}^{\bm{X}_{1},\bm{X}_{2}}(\bm{Y}) - T_{\beta}^{\bm{X}_{1},\bm{X}_{2}}(\bm{Y})\right\Vert\leq \epsilon/3,
\end{eqnarray}
by part (ii) in Lemma~\ref{lma:continuity ii and iii} again. For Part II, we have
\begin{eqnarray}\label{eq7:thm:GDOI continuity}
\lefteqn{\left\Vert T_{\tilde{\beta}}^{\bm{X}_{1,\ell_1},\bm{X}_{2,\ell_2}}(\bm{Y}) - T_{\tilde{\beta}}^{\bm{X}_{1},\bm{X}_{2}}(\bm{Y})\right\Vert}\nonumber \\
&=&\left\Vert T_{\sum\limits_{k'_1=1,k'_2=0}^{k'_1+k'_2=\ell+1} d_{k'_1,k'_2}f^{[1]}_{k'_1 }(x_1,x_2)x_2^{k'_2}}^{\bm{X}_{1,\ell_1},\bm{X}_{2,\ell_2}}(\bm{Y}) - T_{\sum\limits_{k'_1=1,k'_2=0}^{k'_1+k'_2=\ell+1} d_{k'_1,k'_2}f^{[1]}_{k'_1 }(x_1,x_2)x_2^{k'_2}}^{\bm{X}_{1},\bm{X}_{2}}(\bm{Y})\right\Vert\nonumber \\
&=_1&\left\Vert \sum\limits_{k'_1=1,k'_2=0}^{k'_1+k'_2=\ell+1} d_{k'_1,k'_2} T_{ f^{[1]}_{k'_1}(x_1,x_2)}^{\bm{X}_{1,\ell_1},\bm{X}_{2,\ell_2}}(T_{x_2^{k'_2}}^{\bm{X}_{1,\ell_1},\bm{X}_{2,\ell_2}}(\bm{Y})) -\sum\limits_{k'_1=1,k'_2=0}^{k'_1+k'_2=\ell+1} d_{k'_1,k'_2} T_{f^{[1]}_{k'_1}(x_1,x_2)}^{\bm{X}_{1},\bm{X}_{2}}(T_{x_2^{k'_2}}^{\bm{X}_{1},\bm{X}_{2}}(\bm{Y}))\right\Vert\nonumber \\
&\leq&\left\Vert \sum\limits_{k'_1=1,k'_2=0}^{k'_1+k'_2=\ell+1} d_{k'_1,k'_2} T_{ f^{[1]}_{k'_1}(x_1,x_2)}^{\bm{X}_{1,\ell_1},\bm{X}_{2,\ell_2}}(T_{x_2^{k'_2}}^{\bm{X}_{1,\ell_1},\bm{X}_{2,\ell_2}}(\bm{Y})) -\sum\limits_{k'_1=1,k'_2=0}^{k'_1+k'_2=\ell+1} d_{k'_1,k'_2} T_{f^{[1]}_{k'_1}(x_1,x_2)}^{\bm{X}_{1,\ell_1},\bm{X}_{2,\ell_2}}(T_{x_2^{k'_2}}^{\bm{X}_{1},\bm{X}_{2}}(\bm{Y}))\right\Vert \nonumber \\
&&+\left\Vert \sum\limits_{k'_1=1,k'_2=0}^{k'_1+k'_2=\ell+1} d_{k'_1,k'_2} T_{ f^{[1]}_{k'_1}(x_1,x_2)}^{\bm{X}_{1,\ell_1},\bm{X}_{2,\ell_2}}(T_{x_2^{k'_2}}^{\bm{X}_{1},\bm{X}_{2}}(\bm{Y})) -\sum\limits_{k'_1=1,k'_2=0}^{k'_1+k'_2=\ell+1} d_{k'_1,k'_2} T_{f^{[1]}_{k'_1}(x_1,x_2)}^{\bm{X}_{1},\bm{X}_{2}}(T_{x_2^{k'_2}}^{\bm{X}_{1},\bm{X}_{2}}(\bm{Y}))\right\Vert\nonumber \\
&\leq&\sum\limits_{k'_1=1,k'_2=0}^{k'_1+k'_2=\ell+1} |d_{k'_1,k'_2}|\left\Vert T_{ f^{[1]}_{k'_1}(x_1,x_2)}^{\bm{X}_{1,\ell_1},\bm{X}_{2,\ell_2}}(T_{x_2^{k'_2}}^{\bm{X}_{1,\ell_1},\bm{X}_{2,\ell_2}}(\bm{Y})) -T_{f^{[1]}_{k'_1}(x_1,x_2)}^{\bm{X}_{1,\ell_1},\bm{X}_{2,\ell_2}}(T_{x_2^{k'_2}}^{\bm{X}_{1},\bm{X}_{2}}(\bm{Y}))\right\Vert \nonumber \\
&&+\sum\limits_{k'_1=1,k'_2=0}^{k'_1+k'_2=\ell+1}|d_{k'_1,k'_2}|\left\Vert T_{ f^{[1]}_{k'_1}(x_1,x_2)}^{\bm{X}_{1,\ell_1},\bm{X}_{2,\ell_2}}(T_{x_2^{k'_2}}^{\bm{X}_{1},\bm{X}_{2}}(\bm{Y})) -  T_{f^{[1]}_{k'_1}(x_1,x_2)}^{\bm{X}_{1},\bm{X}_{2}}(T_{x_2^{k'_2}}^{\bm{X}_{1},\bm{X}_{2}}(\bm{Y}))\right\Vert \nonumber \\
&\leq_2&\epsilon/6 + \epsilon/6 =\epsilon/3,
\end{eqnarray}
where we apply Lemma~\ref{lma: GDOI linear homomorphism} in $=_1$, apply (i) from Lemma~\ref{lma:continuity ii and iii} to obtain the first $\epsilon/6$ in $\leq_2$, and apply Theorem~\ref{thm:GDOI continuity dd} to obtain the second $\epsilon/6$ in $\leq_2$. Therefore,  this theorem is proved by combining Eq.~\eqref{eq5:thm:GDOI continuity}, Eq.~\eqref{eq6:thm:GDOI continuity}, and Eq.~\eqref{eq7:thm:GDOI continuity}. 
$\hfill\Box$

\section{Differentiation of Functions with Arguments as Continuous Spectrum Operators}\label{sec:Differentiation of Functions with Arguments as Continuous Spectrum Operators}

We consider a family of operators indexed by $t$ as $\bm{X}(t)$ such that $t \rightarrow \bm{X}(t)$ is a smooth $C^\infty$ function with norm topology. We have the following Lemma~\ref{lma:d x power ell} about expressing the derivative of the $\bm{X}^\ell(t)$, where $\ell$ is a natural number, with respect to $t$ in terms of the GDOI.
\begin{lemma}\label{lma:d x power ell}
For any natural number $\ell$, we have
\begin{eqnarray}\label{eq1:lma:d x power ell}
\frac{d\bm{X}^\ell(t)}{dt}&=&T_{\sum\limits_{i=0}^{\ell-1}\lambda_1(t)^{\ell-i-1}\lambda_2(t)^i}^{\bm{X}(t),\bm{X}(t)}\left(\frac{d \bm{X}(t)}{dt}\right),
\end{eqnarray}
where $\lambda_1(t) \in \sigma(\bm{X}(t))$ and $\lambda_2(t) \in \sigma(\bm{X}(t))$.
\end{lemma}
\textbf{Proof:}
We will establish this lemma by induction.  For $\ell=1$, it is trivial to have Eq.~\eqref{eq1:lma:d x power ell}.

For $\ell=2$, the L.H.S. of Eq.~\eqref{eq1:lma:d x power ell} is
\begin{eqnarray}\label{eq2:lma:d x power ell}
\frac{d\bm{X}^2(t)}{dt}&=&\bm{X}(t) \frac{d \bm{X}(t)}{dt} + \frac{d \bm{X}(t)}{dt} \bm{X}(t).
\end{eqnarray}
Since we have 
\begin{eqnarray}\label{eq3:lma:d x power ell}
\bm{X}(t)&=&\int\limits_{\lambda_1 \in \sigma(\bm{X}(t))}\lambda_1(t) d\bm{E}_{\bm{X}(t)}(\lambda_1(t))+
\int\limits_{\lambda_1(t) \in \sigma(\bm{X}(t))}\left(\bm{X}(t)-\lambda_1(t)\bm{I}\right)d\bm{E}_{\bm{X}(t)}(\lambda_1(t)),
\end{eqnarray}
and 
\begin{eqnarray}\label{eq3.1:lma:d x power ell}
\bm{X}(t)&=&\int\limits_{\lambda_2 \in \sigma(\bm{X}(t))}\lambda_2(t) d\bm{E}_{\bm{X}(t)}(\lambda_2(t))+
\int\limits_{\lambda_2(t) \in \sigma(\bm{X}(t))}\left(\bm{X}(t)-\lambda_2(t)\bm{I}\right)d\bm{E}_{\bm{X}(t)}(\lambda_2(t)),
\end{eqnarray}
by applying Eq.~\eqref{eq3:lma:d x power ell} and Eq.~\eqref{eq3.1:lma:d x power ell} to Eq.~\eqref{eq2:lma:d x power ell},  we obtain
\begin{eqnarray}\label{eq4:lma:d x power ell}
\lefteqn{\frac{d\bm{X}^2(t)}{dt}=}\nonumber \\
&& \left(\int\limits_{\lambda_1(t) \in \sigma(\bm{X}(t))}\lambda_1(t) d\bm{E}_{\bm{X}(t)}(\lambda_1(t))+
\int\limits_{\lambda_1(t) \in \sigma(\bm{X}(t))}\left(\bm{X}(t)-\lambda_1(t)\bm{I}\right)d\bm{E}_{\bm{X}(t)}(\lambda_1(t))\right) \frac{d \bm{X}(t)}{dt}\nonumber \\
&& + \frac{d \bm{X}(t)}{dt}\left(\int\limits_{\lambda_2(t) \in \sigma(\bm{X}(t))}\lambda_2(t) d\bm{E}_{\bm{X}(t)}(\lambda_2(t))+
\int\limits_{\lambda_2(t) \in \sigma(\bm{X}(t))}\left(\bm{X}(t)-\lambda_2(t)\bm{I}\right)d\bm{E}_{\bm{X}(t)}(\lambda_2(t))\right).
\end{eqnarray}

For $\ell=2$, the R.H.S. of Eq.~\eqref{eq1:lma:d x power ell} is obtained by applying GDOI definition given by Eq.~\eqref{eq1:  GDOI def}
\begin{eqnarray}\label{eq5:lma:d x power ell}
\lefteqn{T_{\lambda_1(t)+\lambda_2(t)}^{\bm{X}(t),\bm{X}(t)}\left(\frac{d \bm{X}(t)}{dt}\right)=}\nonumber \\
&& \int\limits_{\lambda_1(t) \in \sigma(\bm{X}(t))}\int\limits_{\lambda_2(t) \in \sigma(\bm{X}(t))}(\lambda_1(t)+\lambda_2(t))d\bm{E}_{\bm{X}(t)}(\lambda_1(t)) \frac{d \bm{X}(t)}{dt} d\bm{E}_{\bm{X}(t)}(\lambda_2(t)) \nonumber \\
&&+ \int\limits_{\lambda_1(t) \in \sigma(\bm{X}(t))}\int\limits_{\lambda_2(t) \in \sigma(\bm{X}(t))}d\bm{E}_{\bm{X}(t)}(\lambda_1(t)) \frac{d \bm{X}(t)}{dt}\left(\bm{X}(t)-\lambda_2(t)\bm{I}\right)d\bm{E}_{\bm{X}(t)}(\lambda_2(t))  \nonumber \\
&&+ \int\limits_{\lambda_1(t) \in \sigma(\bm{X}(t))}\int\limits_{\lambda_2(t) \in \sigma(\bm{X}(t))} \left(\bm{X}(t)-\lambda_1(t)\bm{I}\right) d\bm{E}_{\bm{X}(t)}(\lambda_1(t)) \frac{d \bm{X}(t)}{dt}d\bm{E}_{\bm{X}(t)}(\lambda_2(t)) \nonumber \\
&=_1&\left(\int\limits_{\lambda_1(t) \in \sigma(\bm{X}(t))}\lambda_1(t)d\bm{E}_{\bm{X}(t)}(\lambda_1(t))\right)\frac{d \bm{X}(t)}{dt} + \frac{d \bm{X}(t)}{dt}  \left(\int\limits_{\lambda_2(t) \in \sigma(\bm{X}(t))}\lambda_2(t) d\bm{E}_{\bm{X}(t)}(\lambda_2(t))\right)\nonumber \\
&&+ \frac{d \bm{X}(t)}{dt}\left(\int\limits_{\lambda_2(t) \in \sigma(\bm{X}(t))}\left(\bm{X}(t)-\lambda_2(t)\bm{I}\right)d\bm{E}_{\bm{X}(t)}(\lambda_2(t)) \right)\nonumber \\
&&+\left(\int\limits_{\lambda_1(t) \in \sigma(\bm{X}(t))} \left(\bm{X}(t)-\lambda_1(t)\bm{I}\right) d\bm{E}_{\bm{X}(t)}(\lambda_1(t))\right)\frac{d \bm{X}(t)}{dt},
\end{eqnarray}
where we apply $\int\limits_{\lambda_1(t) \in \sigma(\bm{X}(t))}d\bm{E}_{\bm{X}(t)}(\lambda_1(t))=\bm{I}$ and $\int\limits_{\lambda_2(t) \in \sigma(\bm{X}(t))}d\bm{E}_{\bm{X}(t)}(\lambda_2(t))=\bm{I}$ in $=_1$.  We can find that Eq.~\eqref{eq4:lma:d x power ell} agress with Eq.~\eqref{eq5:lma:d x power ell}. Hence, we  have Eq.~\eqref{eq1:lma:d x power ell} at $\ell = 2$. 

Suppose, we have Eq.~\eqref{eq1:lma:d x power ell} at $\ell = m$. Then, for the case with $\ell=m+1$, we have
\begin{eqnarray}\label{eq6:lma:d x power ell}
\frac{d\bm{X}^{m+1}(t)}{dt}&=& \frac{d\bm{X}(t)}{dt}\bm{X}^{m}(t)+\bm{X}(t)\frac{d\bm{X}^{m}(t)}{dt}\nonumber \\
&=_1& T_{\lambda_2(t)^m}^{\bm{X}(t),\bm{X}(t)}\left(\frac{d \bm{X}(t)}{dt}\right)+
T_{\lambda_1(t)\sum\limits_{i=0}^{m-1}\lambda_1(t)^{m-1-i}\lambda_2(t)^i}^{\bm{X}(t),\bm{X}(t)}\left(\frac{d \bm{X}(t)}{dt}\right) \nonumber \\
&=_2& T_{\sum\limits_{i=0}^{m}\lambda_1(t)^{m-i}\lambda_2(t)^i}^{\bm{X}(t),\bm{X}(t)}\left(\frac{d \bm{X}(t)}{dt}\right),
\end{eqnarray}
where we apply induction assumption for $\ell=m$ in $=_1$and Lemma~\ref{lma: GDOI linear homomorphism} in $=_2$. Therefore, this lemma is proved by induction. 
$\hfill\Box$

Then, we have the following Theorem~\ref{thm:d f X dt}.
\begin{theorem}\label{thm:d f X dt}
Let $f$ be a complex-valued and differentiable function, we have 
\begin{eqnarray}\label{eq1:thm:d f X dt}
\frac{df(\bm{X}(t))}{dt}&=&T_{f^{[1]}}^{\bm{X}(t),\bm{X}(t)}\left(\frac{d \bm{X}(t)}{dt}\right),
\end{eqnarray}
where $T_{f^{[1]}}^{\bm{X}(t),\bm{X}(t)}\left(\frac{d \bm{X}(t)}{dt}\right)$ is the GDOI with continuous spectrum operators. We assume that $\sigma(\bm{X}_t)$ within a compact set for all $t$. 
\end{theorem}
\textbf{Proof:}
From Lemma~\ref{lma: GDOI linear homomorphism} and Lemma~\ref{lma:d x power ell}, we have
\begin{eqnarray}\label{eq2:thm:d f X dt}
\frac{d p(\bm{X}(t))}{dt}&=&T_{p^{[1]}}^{\bm{X}(t),\bm{X}(t)}\left(\frac{d \bm{X}(t)}{dt}\right),
\end{eqnarray}
where $p(\bm{X}(t))$ is a polynomial with respect to the variable operator $\bm{X}(t)$.
 
Since we have
\begin{eqnarray}\label{eq3:thm:d f X dt}
\lefteqn{\left\Vert \frac{df(\bm{X}(t))}{dt} -  T_{f^{[1]}}^{\bm{X}(t),\bm{X}(t)}\left(\frac{d \bm{X}(t)}{dt}\right) \right\Vert =}\nonumber \\
&& 
\left\Vert \frac{df(\bm{X}(t))}{dt} -  \frac{d p(\bm{X}(t))}{dt} + \frac{d p(\bm{X}(t))}{dt} - T_{p^{[1]}}^{\bm{X}(t),\bm{X}(t)}\left(\frac{d \bm{X}(t)}{dt}\right) \right. \nonumber \\
&& \left.  + T_{p^{[1]}}^{\bm{X}(t),\bm{X}(t)}\left(\frac{d \bm{X}(t)}{dt}\right) -  T_{f^{[1]}}^{\bm{X}(t),\bm{X}(t)}\left(\frac{d \bm{X}(t)}{dt}\right) \right\Vert \nonumber \\
&\leq & \left\Vert \frac{df(\bm{X}(t))}{dt} -  \frac{d p(\bm{X}(t))}{dt} \right\Vert + \left\Vert \frac{d p(\bm{X}(t))}{dt} - T_{p^{[1]}}^{\bm{X}(t),\bm{X}(t)}\left(\frac{d \bm{X}(t)}{dt}\right) \right\Vert \nonumber \\
&&+ \left\Vert  T_{p^{[1]}}^{\bm{X}(t),\bm{X}(t)}\left(\frac{d \bm{X}(t)}{dt}\right) -  T_{f^{[1]}}^{\bm{X}(t),\bm{X}(t)}\left(\frac{d \bm{X}(t)}{dt}\right) \right\Vert \nonumber \\
&\leq_1& \epsilon/2+0+\epsilon/2= \epsilon,
\end{eqnarray}
where the inequality $\leq_1$ comes from Oka-Weil approximation theorem ($p \rightarrow f$), Eq.~\eqref{eq2:thm:d f X dt}, and the part (ii) in Lemma~\ref{lma:continuity ii and iii}. This theorem is proved accordingly.  
$\hfill\Box$

\bibliographystyle{IEEETran}
\bibliography{SpecialCase_and_DOI_Bib}

\end{document}